\documentclass[11pt]{amsart}

\usepackage{mathtools}
\usepackage{amssymb}
\usepackage{fullpage}
\usepackage{amsthm}
\usepackage{tikz}
\usepackage{tikz-cd}
\usetikzlibrary{positioning}
\usepackage{ytableau}
\usepackage[hidelinks]{hyperref}

\usepackage{biblatex}
\newtheorem{thm}{Theorem}[section]
\newtheorem{lemma}[thm]{Lemma}
\newtheorem{cor}[thm]{Corollary}
\newtheorem{prop}[thm]{Proposition}

\theoremstyle{definition}
\newtheorem{defn}[thm]{Definition}
\newtheorem{example}[thm]{Example}

\DeclareMathOperator{\GL}{GL}
\DeclareMathOperator{\SL}{SL}
\DeclareMathOperator{\Sor}{SO}
\DeclareMathOperator{\Or}{O}
\DeclareMathOperator{\Un}{U}
\DeclareMathOperator{\Sp}{Sp}

\newcommand{\g}{\mathfrak{g}}
\newcommand{\kk}{\mathfrak{k}}
\newcommand{\p}{\mathfrak{p}}
\newcommand{\h}{\mathfrak{h}}
\newcommand{\gl}{\mathfrak{gl}}
\newcommand{\sll}{\mathfrak{sl}}
\newcommand{\sor}{\mathfrak{so}}

\DeclareMathOperator{\Hom}{Hom}
\DeclareMathOperator{\Ind}{Ind}
\DeclareMathOperator{\Res}{Res}
\DeclareMathOperator{\RR}{RR}
\DeclareMathOperator{\chr}{char}
\DeclareMathOperator{\wt}{wt}
\DeclareMathOperator{\mwt}{M-wt}
\DeclareMathOperator{\tr}{tr}
\DeclareMathOperator{\HC}{HC}
\DeclareMathOperator{\LC}{LC}
\DeclareMathOperator{\HL}{HL}
\DeclareMathOperator{\LL}{LL}
\DeclareMathOperator{\LR}{LR}
\DeclareMathOperator{\Rect}{Rect}
\DeclareMathOperator{\SST}{SST}
\DeclareMathOperator{\JK}{JK}
\DeclareMathOperator{\pRJK}{pRJK}
\DeclareMathOperator{\RJK}{RJK}
\DeclareMathOperator{\SM}{SM}
\DeclareMathOperator{\AM}{AM}
\DeclareMathOperator{\M}{M}
\DeclareMathOperator{\Lie}{Lie}
\DeclareMathOperator{\Spec}{Spec}
\DeclareMathOperator{\mult}{mult}
\DeclareMathOperator{\Mult}{Mult}
\DeclareMathOperator{\gr}{gr}

\newcommand{\Z}{\mathbb{Z}}
\newcommand{\R}{\mathbb{R}}
\newcommand{\C}{\mathbb{C}}
\newcommand{\Qu}{\mathbb{H}}

\newcommand{\Ha}{\mathcal{H}}
\newcommand{\B}{\mathcal{B}}
\newcommand{\Par}{\mathcal{P}}

\title{Graded multiplicities in the Kostant-Rallis setting}

\date{}

\author{Andrew Frohmader}

\begin{document}

\begin{abstract}
   This paper contains two main results. First, we provide combinatorial branching rules for $\GL_n \downarrow \Or_n$ and $\GL_{2n} \downarrow \Sp_{2n}$ extending the Littlewood restriction rules. Second, we use these branching rules and the combinatorics of $\GL_n$-crystals to derive a formula for the graded multiplicity of a $K$-type in the regular functions on the $K$-nilpotent cone for $\GL(n, \R)$, $\GL(n, \C)$ and $\GL(n, \Qu)$. For $\GL(n,\R)$, work of Schmid and Vilonen identifies these graded multiplicities with the Hodge $K$-character of the tempered spherical principal series with real infinitesimal character.
\end{abstract}

\maketitle
\tableofcontents

\section{Introduction}

\subsection{Kostant-Rallis setting}
Let $G$ be a connected complex reductive algebraic group, $G_{\R}$ a real form with maximal compact subgroup $K_{\R}$, and $\g_{\R} = \kk_{\R} \oplus \p_{\R}$ the Cartan decomposition of its Lie algebra. Let $\theta: G \to G$ be the corresponding complexified Cartan involution, $K = G^\theta$ the group of $\theta$-fixed points, $\g$ the Lie algebra of $G$, $\kk = \{ X \in \g : d\theta(X) = X \}$ the Lie algebra of $K$, and $\p = \{X \in \g : d\theta(X) = -X \}$, so $\g = \kk \oplus \p$. 

The adjoint representation of $G$ restricts to give a representation of $K$ on $\p$. Let $\mathfrak{a}_{\R}$ be a maximal abelian subalgebra of $\p_{\R}$ and let $\mathfrak{a} \subset \p$ be its complexification. Let $M$ be the centralizer of $\mathfrak{a}$ in $K$. Let $\C[\p]$ denote the ring of polynomial functions on $\p$, $\C[\p]^K$ the ring of $K$-invariant polynomials, and $\C[\p]^K_+$ the ideal in $\C[\p]$ generated by all elements of $\C[\p]^K$ without a constant term. Let $\mathcal{N} \subset \g$ be the nilpotent cone in $\g$ and $\mathcal{N}_\theta = \mathcal{N} \cap \p$ be the $K$-nilpotent cone. We have

\begin{thm}[{\cite[Theorem 14]{KostantRallis1971}}]
    The ideal $\C[\p]^K_+$ is the vanishing ideal of $\mathcal{N}_\theta$, i.e. $\C[\p]^K_+ = \{f \in \C[\p] : f(x) = 0 \text{ for all }  x \in \mathcal{N}_\theta \}$.
\end{thm}

Similarly, let $S(\p)$ denote the symmetric algebra on $\p$ regarded as a $K$-module, $S(\p)^K$ the $K$-invariant symmetric tensors, and $S(\p)^K_+$ the ideal in $S(\p)^K$ generated by all elements without a constant term. Each element  $v \in S(\p)$ is naturally associated with a constant-coefficient differential operator $\partial_v$ on $\p$ and a polynomial $f \in \C[\p]$ is called \textit{$K$-harmonic} if $\partial_vf = 0$ for all $v \in S(\p)^K_+$. Let $\Ha \subset \C[\p]$ be the subspace of all harmonic polynomials on $\p$. Crucially for us, the spaces $\C[\p]$, $\C[\p]^K$, and $\Ha$ are all naturally graded by degree. 

\begin{thm}[{\cite[Theorems 15, 16, 17]{KostantRallis1971}}]\label{thm:Kostant-Rallis}
    The space $\C[\p]$ is a free module over the ring of $K$-invariant polynomials, with
    $$\C[\p] = \C[\p]^K \otimes \Ha.$$

    \noindent Furthermore, as a representation of $K$, $\Ha$ is isomorphic to the representation algebraically induced from the trivial representation of $M$ to $K$,
    $$\Ha \cong \Ind_M^K 1.$$
\end{thm}

We call this setup the \textit{Kostant-Rallis setting}. We draw a couple of conclusions from the two theorems above. 

\begin{cor} \label{cor:H=nilpotent cone}
    As a graded $K$-representation, $\Ha$ is isomorphic to $\C[\mathcal{N}_\theta]$. Moreover, as an ungraded $K$-representation, $\Ha$ is isomorphic to the $K$-finite vectors of the spherical principal series representations of $G_{\R}$.
\end{cor}

 We can immediately compute the multiplicity of a $K$-type in $\Ha$ by applying Frobenius reciprocity to $\Ind_M^K 1$. Let $\widehat{K}$ denote the set of equivalence classes of finite-dimensional irreducible complex representations of $K$ with rational matrix coefficients. For $\nu \in \widehat{K}$, let $\pi^\nu_K$ be a model for $\nu$. 

\begin{cor}\label{cor:K-type mult in princ series}
    The multiplicity of the irreducible $K$-representation $\pi^\nu_K$ in $\Ha$ is the dimension of the subspace of $M$-invariant vectors in $\pi^\nu_K$, i.e. $\dim (\pi^\nu_K)^M$.
\end{cor}

It is more difficult to describe the graded multiplicities. For any $K$-module $V$ and any $\nu \in \widehat{K}$, set $\mult(\pi_K^\nu, V) \coloneqq \dim \Hom_K(\pi^\nu_K, V)$. If $V = \bigoplus_{d \geq 0} V_d$ is graded and $q$ is an indeterminate, define
\[
\mult_q(\pi^\nu_K,V) \coloneqq \sum_{d \geq 0} \mult(\pi^\nu_K,V_d)q^d.
\]
We call this the \textit{graded multiplicity} of $\pi_K^\nu$ in $V$. Specializing $q=1$, $\mult_1(\pi_K^\nu, V) = \mult(\pi_K^\nu, V)$. Define
$$m^{\nu,0}_{(G,K)}(q) \coloneqq \mult_q(\pi^\nu_K, \Ha).$$
\noindent These polynomials are the primary objects of this paper. The parameter $\nu$ denotes the $K$-type and the $0$ indicates that we are considering the trivial representation of $M$. We will only consider the trivial $M$-type here, but it would be interesting to consider other $M$-types in connection with $K$-types in non-spherical principal series, as discussed in Section \ref{sec:motivation-real-groups} and Section \ref{sec:principal-series-M-polynomials}.

A large body of literature in combinatorial representation theory explores the special case of the Kostant-Rallis setting corresponding to complex Lie groups. Consider $G_{\R} = \GL(n, \C)$, $\Or(n, \C)$, and $ \Sp(2n, \C)$ regarded as real groups with $K_{\R} = \Un(n)$, $\Or(n)$, and $\Sp(n)$ respectively. As we will primarily work with complex algebraic groups, we denote $\GL(n, \C), \Or(n, \C), \Sp(2n, \C)$ by $\GL_n, \Or_n, \Sp_{2n}$. Complexifying, we have $(G,K) = (K^2, K) = (\GL_n^2, \GL_n)$, $(\Or_n^2, \Or_n)$, and $(\Sp_{2n}^2, \Sp_{2n})$. Here $K^2 = K \times K$ and more generally, $K^m = K \times \dots \times K$.  The Cartan involution exchanges the two factors of $K$. In this case, $\p$ is naturally identified with $\kk$, $M$ is the maximal torus of $K$, and the representation of $K$ on $\p$ is the adjoint representation. Since $M$ is a maximal torus, the subspace of $M$-invariant vectors in Corollary \ref{cor:K-type mult in princ series} is the zero weight subspace. Hence, in this setting, the multiplicity of a $K$-type $\pi^\nu_K$ in $\Ha$ is the dimension of the zero weight space of $\pi^\nu_K$. We call this the \textit{Kostant setting}, since it corresponds to \cite{Kostant1963}.

In the Kostant setting, a lot of work has gone into describing the graded multiplicity $m^{\nu,0}_{(K^2,K)}(q)$, which Kostant called \textit{generalized exponents}, in terms of elements of the zero weight subspace. For example, a theorem of Hesselink establishes that $m^{\nu,0}_{(K^2,K)}(q)$ equals the Lusztig $q$-analogue associated to the zero weight multiplicity in $\pi^\nu_K$, \cite{Hesselink1980}. For more background, see \cite{NelsenRam2003, LecouveyLenart2020, JosephLetzterZelikson2000}.

In this paper, we consider the real groups $\GL(n, \R)$, $\GL(n, \C)$, and $\GL(n, \Qu)$. The data we need is summarized in the table below, see \cite[Section 12.3.2]{GoodmanWallach2009}. Let $\M_{p,q}$ denote the space of complex $p \times q$ matrices, $\M_n \coloneqq \M_{n,n}$, $\SM_n$ the $n \times n$ symmetric matrices, and $\AM_n$ the $n \times n$ antisymmetric matrices. 

\begin{center}
\begin{tabular}{c|c|c|c}
    $G_{\R}, K_{\R}$ & $G, K$ & $\p$ & $M$ \\
    \hline
    $\GL(n, \R), \Or(n)$ & $\GL_n, \Or_n$ &  $\SM_n$ & $\Or_1^{n}$\\
    $\GL(n, \C), \Un(n)$ & $\GL_n^2, \GL_n$ & $\M_n$ & $\GL_1^{n}$\\
    $\GL(n, \Qu), \Sp(n)$ & $\GL_{2n}, \Sp_{2n}$ & $\AM_{2n}$ & $\Sp_2^{n}$  \\
\end{tabular}
\end{center}

For these cases, our main result is
\begin{thm}[Theorem \ref{thm:graded mults} below]\label{thm:intro-gradedmults}
    $$m^{\nu,0}_{(G,K)}(q) = \sum_{T \in (\mathcal{T}_K^\nu)_0} q^{d(T)},$$
    \noindent where $(\mathcal{T}^\nu_K)_0$ is a subset of semistandard tableaux indexing a basis for the space of $M$-invariant vectors in the irreducible representation $\pi^\nu_K$ and $d(T)$ is a statistic computed using crystal operators.
\end{thm}

The $(\GL_n^2, \GL_n)$ case falls within the Kostant setting and $d(T)$ can be seen as a version of the charge statistic on rational $\GL_n$-tableaux, \cite{LascouxSchutzenberger1978}. The other two cases are outside the Kostant setting. Although crystals play an important role in proving this theorem, the resulting statistic can be computed in a purely combinatorial manner via the signature rule. We include many explicit computations in Section \ref{sec:graded-mults}. 

The literature on graded multiplicities in the Kostant-Rallis setting is less developed than in the Kostant setting. Kostant and Rallis describe the graded multiplicity of an irreducible representation $\pi^\nu_K$ in $\Ha$ in terms of the eigenvalues of a certain element of $\kk$, see \cite[Theorem 21]{KostantRallis1971}. In \cite{WallachWillenbring2000}, Wallach and Willenbring obtain Hesselink-type formulas for some examples including $(\GL_{2n}, \Sp_{2n})$, $(\Sor_{2n+2},\Sor_{2n+1})$, and $(E_6, F_4)$. Wallach and Willenbring also worked out the example of $(\SL_4, \Sor_4)$ explicitly, and other results in special cases have appeared; see \cite{JohnsonWallach1977, WillenbringVanGroningen2014}. Mason-Brown obtained a Hesselink-type formula for $\C[\mathcal{N}_\theta]$ in the case where $G_{\R}$ is split modulo center \cite{Mason-Brown2022}. In general, prior to the present work, explicit combinatorial formulas were only available in the stable range governed by the classical restriction rules of Littlewood \cite{Littlewood1945, Littlewood1944, Willenbring2002, HoweTanWillenbring2005, HoweTanWillenbring2008}.

\subsection{Branching rules}
In order to prove the graded multiplicity formula above, we establish generalizations of the Littlewood restriction rules valid for all parameters; in particular there are no stable-range constraints. These results build on \cite{Kwon2018, KwonJang2021, LecouveyLenart2020} and may be of independent interest. For $\lambda \in \widehat{G}$ and $\nu \in \widehat{K}$, let $b^\lambda_\nu \coloneqq \mult(\pi_K^\nu, \Res^G_K(\pi^\lambda_G))$ denote the \textit{branching multiplicity}. We will drop $\Res^G_K$ from the notation and simply write $\mult(\pi_K^\nu, \pi^\lambda_G)$ since the first argument makes clear that the second representation is being restricted to $K$.

\begin{thm}[Theorem \ref{thm:main_branching} below]~
\begin{enumerate}
    \item For $\lambda \in \widehat{\GL}_{2n}$ corresponding to a polynomial representation and $\nu \in \widehat{\Sp}_{2n}$,
    $$\mult(\pi^\nu_{\Sp_{2n}}, \pi^\lambda_{\GL_{2n}}) = \sum_{\mu \in \Par^{(1,1)}_{2n}} c^\lambda_{\mu \nu }(\Sp_{2n}).$$
    \item For $\lambda \in \widehat{\GL}_n$ corresponding to a polynomial representation and $\nu \in \widehat{\Or}_n$,
    $$\mult(\pi^\nu_{\Or_{n}}, \pi^\lambda_{\GL_{n}}) =  \sum_{\mu \in \Par_n^{(2)}} c^\lambda_{\mu \nu }(\Or_n).$$
\end{enumerate}
\end{thm}

The integers $c^\lambda_{\mu \nu}(K)$ are a generalization of the standard Littlewood-Richardson coefficients and can be computed using mild flag conditions on the usual Littlewood-Richardson tableaux. $\Par^{(2)}_n$ denotes the set of partitions with even rows and length at most $n$. $\Par^{(1,1)}_{2n}$ denotes the set of partitions with length at most $2n$ and even columns. There are a number of examples in Section \ref{sec:branching rules} demonstrating explicit computations, which in our opinion are no worse than the computation of stable branching via the Littlewood restriction rules.

\subsection{Crystal bases}

The main combinatorial tool used throughout this paper is the theory of crystal bases; our conventions follow \cite{BumpSchilling2017}. We review and develop the necessary material in Sections \ref{sec:crystals-ofreductive} and \ref{sect:gln-crystals}. Two features of our exposition are particularly important for the applications below. First, we work with crystals for complex reductive Lie algebras, rather than only semisimple Lie algebras, with $\GL_n$ as the principal example. Second, we treat highest- and lowest-weight constructions symmetrically. In particular, we give lowest-weight counterparts to familiar highest-weight constructions, including Yamanouchi words and the various tableau formulations of the Littlewood--Richardson rule.

The principal tool for passing between highest- and lowest-weight constructions is a crystal anti-isomorphism, which allows the source and target crystals to be different, providing a slight generalization of the usual notion of a crystal involution. We also describe how crystal anti-isomorphisms interact with tensor products. For $\GL_n$-crystals of tableaux, the two maps used here are the Sch\"utzenberger involution
\[
S:\mathcal{B}_n^\nu \longrightarrow \mathcal{B}_n^\nu
\]
and the rotation--complement anti-isomorphism
\[
S^\pi:\mathcal{B}_n^{\lambda/\mu}
\longrightarrow
\mathcal{B}_n^{\lambda^\pi/\mu^\pi}.
\]
Using these maps, we relate several equivalent highest- and lowest-weight forms of the Littlewood--Richardson rule.

This framework allows us to translate the branching formulas of
\cite{Kwon2018,KwonJang2021,LecouveyLenart2020} into a common language of lowest-weight companion tableaux. For example, the rule of Jang and Kwon is stated using highest-weight companion tableaux on the rotated shape $\nu^\pi$; applying $S^\pi$ gives an equivalent description using lowest-weight companion tableaux on $\nu$. These reformulations lead to the simpler and uniform branching rules stated in Theorem \ref{thm:main_branching}.

The lowest-weight formulation is important not only for the statement of the branching rules. The resulting companion tableaux belong to the same tableau sets $\mathcal{T}_K^\nu$ used to parametrize the relevant $M$-types in $\pi_K^\nu$. Consequently, in the proof of the graded multiplicity formula, the branching sums can be reorganized tableau by tableau. The crystal statistics $\varphi_i(T)$ then determine the smallest allowable term in the branching sum and hence the degree $d(T)$ appearing in Theorem \ref{thm:graded mults}.

\subsection{Motivation from real groups} \label{sec:motivation-real-groups}
Next, we discuss the connection between the graded multiplicities studied here and the unitary dual problem for real reductive groups. Recently, it has been shown that the unitarity of an irreducible Harish-Chandra module $\pi$ with real infinitesimal character can be tested using the associated graded of a certain filtration on $\pi$, called the Hodge filtration; see \cite{AdamsTrapaVogan, DavisVilonen2025}. We will refer to the graded $K$-character of $\operatorname{gr}(\pi)$ as the \textit{Hodge $K$-character}, $\chi^H$. This Hodge filtration arises geometrically by realizing $\pi$ as the global sections of a $K$-equivariant $\mathcal{D}$-module on the flag variety via Beilinson-Bernstein localization. Saito's theory of mixed Hodge modules can then be applied to this localized object. The base case in this setting is that of the tempered spherical principal series with real infinitesimal character for $G_{\R}$ split. We have the following,

\begin{thm}\label{thm:intro Hodge filtration}
    Let $G_{\R}$ be split and $\pi$ be the tempered spherical principal series of $G_{\R}$ with real infinitesimal character. Then, as a graded $K$-representation, the associated graded of the Hodge filtration  on $\pi$ is isomorphic to $\Ha$.
\end{thm}

\begin{proof}
    Using Corollary \ref{cor:H=nilpotent cone} this is immediate from a result of Schmid and Vilonen. See \cite[Proposition 13.1]{AdamsTrapaVogan} where the result is attributed to Schmid and Vilonen and \cite[Theorem 1.10]{DavisVilonen2025} where a proof of a slightly more general fact is provided.
\end{proof}

Let $I(\delta)$ denote the tempered principal series with real infinitesimal character induced from $\pi^\delta_M$ and $\gr I(\delta)$ its associated graded with respect to the Hodge filtration. Then Theorem \ref{thm:intro Hodge filtration} says that $\gr I(0) \cong \mathcal{H}$ when $G_{\R}$ is split. Consequently, its Hodge $K$-character is 
\[ \chi^H(I(0)) = \sum_{\nu \in \widehat{K}} m_{(G,K)}^{\nu,0}(q)\chi_K^\nu. \]
Among the three families considered in this paper, this applies directly to $G_{\R}=\GL(n,\R)$. Thus, in the split case, our tableau formulas give an explicit combinatorial description of the Hodge $K$-character of the tempered spherical principal series. We hope that this calculation may shed additional light on Hodge filtrations and unitarity more generally.

It would be interesting to consider whether the statistic $d(T)$, or a suitable extension of it, applied to tableaux of nontrivial $M$-type, could provide information about the Hodge filtrations of non-spherical principal series representations.  To this end, we define the \textit{$M$-polynomial}, $M^\nu_{(G,K)}$ in Section \ref{sec:principal-series-M-polynomials} and give a tableau-theoretic method for computing it. This polynomial simultaneously records the multiplicity of $\pi^\nu_K$ in all principal series representations of $G_{\R}$. For connections to Hodge filtrations, what we want is a $q$-analogue of the $M$-polynomial, where $q$ records the degree of $\pi^\nu_K$ in $\gr I(\delta)$. This is also defined in Section \ref{sec:principal-series-M-polynomials}. The computation of $m^{\nu,0}_{(G,K)}$ contained in this paper is related to $\gr I(0)$ only. Extending the connection with Hodge filtrations to nonsplit groups, computing the coefficients corresponding to nontrivial $M$-types, and extending the results to the remaining classical symmetric pairs are natural directions for future work.

\subsection{Outline}
Finally, we briefly outline our approach. For any $G$-module $V$, let $\Spec(V) \coloneqq \{ \lambda \in \widehat{G} : \mult(\pi_G^\lambda, V) \neq 0\}$ denote the \textit{spectrum} of $V$. In the cases we consider, there is a multiplicity-free decomposition of $\C[\p]$ under the action of $G$, Theorem \ref{thm:mult-free spaces} below. Hence, we can write 
$$\chr_q(\C[\p]) = \sum_{\lambda \in \Spec(\C[\p])} q^{d_G^\lambda} \chi^\lambda_G$$

\noindent where $d_G^\lambda$ is the degree in which $\pi^\lambda_G$ appears in the multiplicity-free decomposition of $\C[\p]$ and $\chi^\lambda_G$ is the character of $\pi^\lambda_G$. Branching from $G$ down to  $K$,
\begin{align*}
    \chr_q(\C[\p]) &= \sum_{\lambda \in \Spec(\C[\p])} q^{d_G^\lambda} \sum_{\nu \in \widehat{K}} b^\lambda_\nu \cdot \chi^\nu_K\\
    &= \sum_{\nu \in \widehat{K}} \chi^\nu_K \sum_{\lambda \in \Spec(\C[\p])} q^{d_G^\lambda}b^\lambda_\nu,
\end{align*}
which we recognize as the decomposition $$\chr_q(\C[\p]) = \sum_{\nu \in \widehat{K}} \chi^\nu_K \mult_q(\pi^\nu_K, \C[\p]).$$
From this, we see 
$$\mult_q(\pi^\nu_K, \C[\p]) = \sum_{\lambda \in \Spec(\C[\p])} q^{d_G^\lambda}b^\lambda_\nu.$$

By Theorem \ref{thm:Kostant-Rallis},
$$\chr_q(\C[\p]) = \chr_q(\C[\p]^K) \cdot \chr_q(\Ha).$$

\noindent Our goal is to compute $\chr_q(\Ha)$. Using the previous identity, we can write
\begin{align*}
    \chr_q(\Ha) &= \frac{1}{\chr_q(\C[\p]^K)} \cdot \chr_q(\C[\p])\\
    &= \frac{1}{\chr_q(\C[\p]^K)} \cdot \sum_{\nu \in \widehat{K}} \chi^\nu_K \sum_{\lambda \in \Spec(\C[\p])} q^{d_G^\lambda}b^\lambda_\nu.
\end{align*}

\noindent So,
$$m^{\nu,0}_{(G,K)}(q) = \mult_q(\pi^\nu_K, \Ha) = \frac{1}{\chr_q(\C[\p]^K)} \sum_{\lambda \in \Spec(\C[\p])} q^{d_G^\lambda}b^\lambda_\nu.$$

But we would like to realize $m^{\nu, 0}_{(G,K)}(q)$ more explicitly in terms of a basis for the $M$-invariants, $(\pi^\nu_K)^M$, as in Theorem \ref{thm:intro-gradedmults} above. To achieve this, we take the following approach:

\begin{enumerate}
    \item Recall the structure of $\C[\p]^K$.
    \item Understand $\C[\p]$ as a multiplicity-free $G$-representation.
    \item Branch from $G$ to $K$.
    \item Factor out the invariants.
\end{enumerate}

Items (1) and (2) are treated in the next section. (3) is handled in Sections 3-5 and (4) is completed in Section 6.

\subsection*{Acknowledgments} I would like to thank Jeb Willenbring for helpful discussions and Jeff Adams for explaining the connection with Hodge filtrations.

\section{Representation theory preliminaries}
We begin by establishing some notation for reductive Lie algebras and Lie groups that we will need in our treatment of crystals. Our notation mainly follows \cite{GoodmanWallach2009}. As we will only need $\GL_n$-crystals for the development of the branching rules, in Section \ref{sec:conventions for gln} we set conventions for $\GL_n$. In Section \ref{section:types_of_tableaux}, we recall parameterizations of the irreducible representations of $\GL_n$, $\Or_n$, and $\Sp_{2n}$ by Young diagrams. We discuss $M$-polynomials and connections with non-spherical principal series in Section \ref{sec:principal-series-M-polynomials}. In Section \ref{sec:ring-of-invariants} we discuss the ring of invariants of $\C[\p]$ and in Section \ref{sec:mult-free-spaces} we review some multiplicity-free actions of $G$ on $\C[\p]$.

\subsection{Reductive Lie algebras} \label{sec:reductive Lie algebras}
Let $H \subset G$ be a maximal torus, $n = \dim H$ the rank of $G$, and $\chi(H)$ the group of rational characters of $H$. Since $\g =\Lie(G)$ is reductive, $\g = \mathfrak{z}(\g) \oplus\g'$, where $\mathfrak{z}(\g)$ is the center of $\g$ and $\g' = [\g, \g]$ is semisimple. Fix the Cartan subalgebra $\h = \Lie(H)$ of $\g$ and set $\h_0 = \h \cap \g'$ which is a Cartan subalgebra of $\g'$. Then we have $\h =  \h_0 \oplus \mathfrak{z}(\g)$. Let $\Phi \subset\h_0^*$ be the root system of $\g'$ relative to $\h_0$, $\Phi^+$ a choice of positive roots and $\Delta = \{\alpha_1, \dots , \alpha_r \}$ the corresponding set of simple roots. Let $r = |\Delta|$ be the semisimple rank of $\g$. For each $\alpha \in \Phi$, choose an $\sll(2,\C)$-triple  $\{e_\alpha, f_\alpha, h_\alpha \} \subset \g'$ with $h_\alpha$ the coroot to $\alpha$ so that $\alpha(h_\alpha) = 2$. For each simple root $\alpha_i$, write $\{e_i, f_i, h_i \}$ for the corresponding triple. Let $X = \{e_i, f_i, h_i \}_{i=1}^r$ be the set of standard generators for the semisimple Lie algebra $\g'$ and extend $X$ to a set of generators for $\g$ by adding a basis $\{h_i \}_{i = r+1}^n$ for $\mathfrak{z}(\g)$.

Let $P(\g)=\{\mu\in\h^*:\mu(h_\alpha)\in\Z\text{ for all }\alpha\in\Phi\}$ be the weight lattice of $\g$. When $\g$ is not semisimple, this is not literally a lattice in the usual sense, since $P(\g)=\mathfrak z(\g)^*\oplus P(\g')$, but we retain the terminology in \cite{GoodmanWallach2009}. Let $\{\varpi_1, \dots , \varpi_r \} \subset \h^*$ be the fundamental weights relative to $\Delta$, extended by zero on $\mathfrak z(\g)$, so that $\varpi_i(h_{j}) = \delta_{ij}$ for $i \in [r]$ and $j \in [n]$. We use $[m]$ to denote the set $\{1, 2, \dots , m \}$. Then, $P(\g') = \{n_1 \varpi_1 + \cdots + n_r \varpi_r : n_i \in \Z \}$ is the semisimple weight lattice. When $\g = \g'$ is semisimple, $P(\g)$ is a free abelian group of rank $r$ having the fundamental weights as basis. If $\mathfrak{z}(\g) \neq 0$, we define additional $\varpi_i$ for $i = r+1, \dots , n$ by the same formula, $\varpi_i(h_j) = \delta_{ij}$ for $j \in [n]$. The $\varpi_i$ are not fundamental weights for $i > r$. The point is that the basis $\{\varpi_1,\dots,\varpi_n\}\subset \h^*$ is dual to the extended basis $\{h_1,\dots,h_n\}\subset \h$. The indices $1,\dots,r$ are
for the simple coroots and fundamental weights of the semisimple Lie algebra $\g'$, while the indices $r+1,\dots,n$ are for $\mathfrak z(\g)$. For $\mu =m_1 \varpi_1 + \cdots + m_n \varpi_n \in P(\g)$ we define the \textit{semisimple weight map}, $ss(\mu) = m_1 \varpi_1 + \cdots + m_r \varpi_r$ projecting $\mu$ onto $P(\g')$.

Define the dominant integral weights to be $P_{++}(\g) = \{\lambda \in \h^* : \lambda(h_i) \in \Z_{\geq 0} \text{ for } i \in [r] \}$. With the above setup, we have 
\[P_{++}(\g) = P_{++}(\g') \oplus \mathfrak{z}(\g)^* = (\Z_{\geq 0} \varpi_1 + \cdots + \Z_{\geq 0} \varpi_r) + (\C \varpi_{r+1} + \cdots + \C \varpi_n).\]
Let $P(G) = \{d \psi : \psi \in \chi(H) \} \subset \h^*$ be the weight lattice of $G$. Then $P(G) \subset P(\g)$ and we set $P_{++}(G) = P(G) \cap P_{++}(\g)$ to be the dominant weights for $G$.

Define a partial order $\geq$ on $P(\g')$ by $\lambda \geq \mu$ if $\lambda - \mu = \sum_{i=1}^r a_i \varpi_i$ where the coefficients $a_i$ are nonnegative for all $i \in [r]$, i.e. $\lambda - \mu \in P_{++}(\g')$. Notice, this is different from the usual partial order $\succeq$ on $P(\g')$ where $\lambda \succeq \mu$ if $\lambda - \mu = \sum_{i=1}^r b_i \alpha_i$ where the coefficients $b_i$ are nonnegative for all $i \in [r]$. 

For each $\alpha \in \Phi$, let $s_\alpha$ be the corresponding reflection of $\h^*$, i.e $s_\alpha(\beta) = \beta - \beta(h_\alpha)\alpha$ for $\beta \in \h^*$ and $W$ the Weyl group generated by the $s_\alpha$. Note that $W$ acts trivially on $\mathfrak z(\g)^*$. Denote the longest Weyl group element $w_0$. Recall $w_0 \Phi^+ = - \Phi^+$ and $w_0$ takes the highest weight of an irreducible finite-dimensional representation of $\g$ to its lowest weight.

\subsection{Conventions for \texorpdfstring{$\GL_n$}{GLn}} \label{sec:conventions for gln}

As discussed in the introduction, $\GL_n$ denotes the complex general linear group. Similarly, $\SL_n$ is the complex special linear group and $\gl_n, \sll_n$ are their Lie algebras. We take $H \subset \GL_n$ to be the subgroup of diagonal matrices. $\g = \gl_n \cong \C \oplus \sll_n$ is reductive, not semisimple. $\h = \h_0 \oplus \C I$ where $\h_0$ consists of diagonal trace zero matrices and $I$ is the $n \times n$ identity element generating the center of $\gl_n$. Let $\epsilon_i \in \h^*$ be the functional $\epsilon_i(\text{diag}(a_1, \dots a_n)) = a_i$, where $\text{diag}(a_1, \dots , a_n)$ is the $n \times n$ diagonal matrix with entries $a_1, \dots , a_n$. Let $E_{i,j}$ denote the $n \times n$ matrix with one in position $(i,j)$ and zeros elsewhere. For $1 \leq i < n$, set $h_i = E_{i,i} - E_{i+1,i+1}$ and $\varpi_i = \epsilon_1 + \cdots + \epsilon_i - \frac{i}{n}(\epsilon_1 + \cdots + \epsilon_n)$. These are the coroots and fundamental weights of $\sll_n$. Define $h_n = I$ and $\varpi_n = \frac{1}{n}(\epsilon_1 + \cdots + \epsilon_n)$. The dominant integral weights of $\gl_n$ are $P_{++}(\gl_n) = P_{++}(\sll_n) \oplus \mathfrak{z}(\gl_n)^* = \Z_{\geq 0} \varpi_1 + \cdots + \Z_{\geq 0} \varpi_{n-1} + \C \varpi_n $. From this, it is easy to see that $P_{++}(\gl_n)$ consists of all weights $m_1\epsilon_1+\cdots+m_n\epsilon_n$ such that $m_i-m_{i+1}\in \Z_{\geq 0}$ for $ 1\leq i<n.$

The weight lattice for $\GL_n$ is 
\[P(\GL_n) = \Z \epsilon_1 + \cdots + \Z \epsilon_n\] 
and from the characterization of $P_{++}(\gl_n)$ above, we see $P_{++}(\GL_n) = P(\GL_n) \cap P_{++}(\gl_n)$ consists of all weights 
$$\text{$\mu =m_1 \epsilon_1 + \cdots + m_n \epsilon_n$ with $m_1 \geq m_2 \geq \cdots \geq m_n$, $m_i \in \Z$}.$$

Define the dominant weights $\gamma_i = \epsilon_1 + \cdots + \epsilon_i$ for $i \in [n]$. The restriction of $\gamma_i$ to $\h_0$ is $\varpi_i$ for $i\in[n-1]$, while $\gamma_n|_{\h_0}=0$. The restriction of $\gamma_i$ to $\mathfrak{z}(\gl_n) = \C I$ is $i \cdot \varpi_n$ for $i \in [n]$. Hence, in terms of our extended fundamental weights basis, $\gamma_i = \varpi_i + i \cdot \varpi_n$ for $i \in [n-1]$ and $\gamma_n = n \cdot \varpi_n$. 

The importance of the $\gamma_i$ is that the elements of $P_{++}(\GL_n)$ can be written uniquely as 
$$\mu = k_1 \gamma_1 + \cdots + k_n \gamma_n, \text{ with $k_i \in \Z_{\geq 0}$ for $i \in [n-1]$ and $k_n \in \Z$}.$$
The connection between the coefficients for $\mu$ in terms of $\gamma_i$ and $\epsilon_i$ is $k_i= m_i - m_{i+1}$, where we define $m_{n+1}= 0$. The restriction of $\mu$ to $\h_0$ is 
$$\mu_0 = k_1 \varpi_1 + \cdots + k_{n-1} \varpi_{n-1}$$
and the restriction of $\mu$ to $\mathfrak{z}(\gl_n)$ is $(m_1 + \cdots + m_n)\varpi_n$. The restriction $\mu \mapsto \mu_0$ is just the semisimple weight map discussed above, and we will also write $ss(\mu) = \mu_0$.

Notice that $\varpi_i|_{\h_0} = \epsilon_1|_{\h_0} + \cdots + \epsilon_i|_{\h_0}$, since $-\frac{i}{n}(\epsilon_1 + \cdots + \epsilon_n)$ is zero on $\h_0$ and we can also take $\{ \epsilon_1|_{\h_0}, \dots , \epsilon_{n-1}|_{\h_0} \}$ as $\Z$-basis for $P(\sll_n)$.

The Weyl group $W$  is the symmetric group  $S_n$ and the longest Weyl group element $w_0$ acts on $P(\GL_n)$ by $w_0(\epsilon_i) =\epsilon_{n+1-i}$. Hence $w_0(\gamma_i) = \gamma_n - \gamma_{n-i}$ for $i \in [n-1]$ and $w_0(\gamma_n) = \gamma_n$. Moreover, $w_0(\varpi_i)=-\varpi_{n-i}$ for $i\in[n-1]$ and $w_0(\varpi_n) = \varpi_n$.

\begin{thm} [{\cite[Theorem 5.5.22, Theorem 3.2.13]{GoodmanWallach2009}}] \label{thm:GLn_highest_lowest} 
    For $\mu \in P_{++}(\GL_n)$ and $\mu_0$ as above, there exists a unique irreducible rational representation $\pi^\mu_{\GL_n}$ of highest weight $\mu$ such that 
    \begin{enumerate}
        \item The restriction of $\pi^\mu_{\GL_n}$ to $\SL_n$ is $\pi^{\mu_0}_{\SL_n}$, i.e. it has highest weight $\mu_0$.
        \item The central elements $zI$ of $\GL_n$ act by multiplication by $z^{m_1+\cdots+m_n}$.
        \item The lowest weight of $\pi^\mu_{\GL_n}$ is $w_0(\mu)$ and the lowest weight of $\pi^{\mu_0}_{\SL_n}$ is $w_0(\mu_0)$.
        \item The highest weight of the dual representation $(\pi^\mu_{\GL_n})^*$ is $-w_0(\mu)$ and the highest weight of $(\pi^{\mu_0}_{\SL_n})^*$ is $- w_0(\mu_0)$.
    \end{enumerate}
\end{thm}

Denote the highest weights $- w_0(\mu)$ and $- w_0(\mu_0)$ by $\mu^*$ and $\mu_0^*$ respectively. We state (3) and (4) in terms of the coordinates $\varpi_i$ since this will be important in our development of the Littlewood-Richardson rule below. 

\begin{lemma} \label{lemma:hw-dual-rep}
$ss$ and $w_0$ commute and for $\mu$ and $\mu_0$ as above,
$$w_0(\mu) = (\sum_{i=1}^n k_i)\gamma_n - k_1 \gamma_{n-1} - k_2\gamma_{n-2} - \cdots - k_{n-1} \gamma_1.$$
\noindent and 
$$ss(w_0(\mu)) = - k_1 \varpi_{n-1} - k_2\varpi_{n-2} - \cdots - k_{n-1} \varpi_1.$$

\noindent Hence, 
$$\mu_0^* = k_1 \varpi_{n-1} + k_2\varpi_{n-2} + \cdots + k_{n-1} \varpi_1.$$
\end{lemma}

\begin{proof}
    Use $w_0(\gamma_i) = \gamma_n - \gamma_{n-i}$ and $ss(\gamma_i) = \varpi_i$ for $i \in [n-1]$ along with $w_0(\gamma_n) = \gamma_n$ and $ss(\gamma_n) = 0$.
\end{proof}

Theorem \ref{thm:GLn_highest_lowest} gives the usual parameterization of the irreducible \textit{rational representations} of $\GL_n$. The irreducible \textit{polynomial representations} are the subset with highest weight $\mu = m_1 \epsilon_1 + \cdots + m_n \epsilon_n$ where $m_1 \geq \cdots \geq m_n \geq 0$ equivalently $\mu = k_1 \gamma_1 + \cdots + k_n \gamma_n$ with $k_n \geq 0$. For $\GL_n$, we extend $\geq$ to $P(\GL_n)$ by $\lambda \geq \mu$ if $\lambda - \mu$ parametrizes a polynomial representation, i.e. $\lambda - \mu = \sum_{i=1}^n a_i \gamma_i$ with $a_i \geq 0$ for all $i \in [n]$. Equivalently, $\lambda - \mu = \sum_{i=1}^n b_i \epsilon_i$ with $b_1 \geq b_2 \geq \cdots \geq b_n\geq 0$.

\subsection{Parameterization of representations} \label{section:types_of_tableaux}
In order to provide concrete multiplicity formulas for $K$-types, we need an explicit parametrization of the irreducible representations of $K$. We recall one way to do this for $K = \GL_n$,  $\Or_n$, and  $\Sp_{2n}$ using Young diagrams. We also give a parametrization of the irreducible representations of $M$ for the cases we are considering.

A \textit{partition} of $k \in \Z_{\geq 0}$ is a tuple of nonnegative integers $\lambda = (\lambda_1, \lambda_2, \dots, \lambda_n)$ such that $\lambda_1 \geq \lambda_2 \geq \cdots \geq \lambda_n \geq 0$ and $\sum_{i=1}^n \lambda_i = k$. We write $|\lambda| = k$ when $\lambda$ is a partition of $k$. The \textit{length} of $\lambda$ is the largest $l$ such that $\lambda_l \neq 0$. Write $\ell(\lambda) = l$. We call the total number of $\lambda_i$ the \textit{ambient length} and write $\ell_a(\lambda) = n$. Notice $\ell(\lambda) \leq \ell_a(\lambda)$. This ambient length is important for us. In notation and examples we will
usually suppress trailing zeros, but operations depending on the ambient
length, such as \(\lambda^\pi\) (defined below), are always taken with respect to the ambient \(n\) determined by context. By the discussion above, we see the irreducible polynomial representations of $\GL_n$ are parametrized by partitions of ambient length $n$. We just take the coefficients of the highest weight $\lambda$ with respect to the $\epsilon_i$ basis: $\lambda = \lambda_1 \epsilon_1 + \cdots + \lambda_n \epsilon_n  \longleftrightarrow (\lambda_1, \dots , \lambda_n)$. Similarly, the irreducible rational representations of $\SL_n$ are parametrized by dominant integral weights of the form
\(
\lambda_0=\lambda_1\epsilon_1|_{\h_0}+\cdots+\lambda_{n-1}\epsilon_{n-1}|_{\h_0},
\)
which we identify with partitions $(\lambda_1,\dots,\lambda_{n-1})$ of ambient length $n-1$.

The \textit{Young diagram} of a partition $\lambda$ is an array of boxes in the lower left or lower right quadrant of $\R^2$ in which the $i^{\text{th}}$ row has $\lambda_i$ boxes. We choose coordinates on $\R^2$ with $y$ negative above the $x$-axis and positive below it and view the $i$th row of $\lambda$ as sitting between $i-1$ and $i$, see Example \ref{ex: young diagrams} below. 

We think of partitions as Young diagrams and do not use separate notation to distinguish them. Hence irreducible polynomial representations of $\GL_n$ are in one-to-one correspondence with Young diagrams of ambient length $n$. Restricting to $\h_0$, Young diagrams of ambient length $n-1$ give our chosen parametrization of the irreducible representations of $\SL_n$.

Let $\Par_n$ denote the set of all partitions of ambient length $n$. By the above discussion, this set parametrizes the polynomial irreducible representations of $\GL_n$. On $\Par_n$, the partial order $\geq$ is given by $\lambda \geq \mu$ if $\lambda - \mu \in \Par_n$, where subtraction is defined by viewing the partition as a tuple of integers in $\Z^n$. Let $\Par^{(2)}_n$ be those partitions of ambient length $n$ with even rows and $\Par^{(1,1)}_{2n}$ be those partitions of ambient length $2n$ with even columns. Notice $\Par_n$, $\Par_n^{(2)}$, and $\Par_{2n}^{(1,1)}$ are additive monoids. 

The \textit{conjugate} diagram $\lambda'$ is obtained by flipping the diagram over its main diagonal. If the Young diagram of $\mu$ has the same ambient length as $\lambda$ and is contained in $\lambda$, we write $\lambda \supseteq \mu$. The \textit{skew diagram} $\lambda/\mu$ is the diagram obtained by removing $\mu$ from $\lambda$. A Young diagram $\lambda$ can be identified with the skew diagram $\lambda/ \varnothing$. Finally, let $\lambda^\pi$ be the skew diagram obtained by rotating $\lambda$ by $180^\circ$ about the origin and then translating vertically by $\ell_a(\lambda)+1=n+1$ (downward translation in our chosen coordinates). This has the effect of reversing the order of the parts of $\lambda$, i.e. $\lambda^\pi = (\lambda_n, \lambda_{n-1}, \dots , \lambda_1)$, and in general $\lambda^\pi$ is no longer a partition since $\lambda_n \leq \lambda_{n-1} \leq \cdots \leq \lambda_1$.
\begin{example}\label{ex: young diagrams}
    If $\lambda = (5,3,3,1,0) = \begin{ytableau}
    \none[1] &  &  & & &\\
    \none[2] &  &  & \\
    \none[3] &  &  & \\
    \none[4] &  \\
    \none[5]\\
\end{ytableau}$, then $|\lambda| = 12$, $\ell(\lambda) = 4$, $\ell_a(\lambda) =5$ and 
$\lambda^\pi = \begin{ytableau}
    \none & \none & \none & \none & \none & \none[1]\\
    \none & \none & \none & \none & &\none[2]\\
    \none& \none& &&&\none[3]\\
    \none& \none & &&& \none[4] \\
    &&&&& \none[5]\\
\end{ytableau} = (0, 1, 3, 3, 5)$.
\end{example}

\noindent Notice that $\lambda^\pi$ is the diagram corresponding to the unique lowest weight of $\pi^\lambda_{\GL_n}$. As mentioned above, we will often suppress trailing zeros if ambient length is clear from context, e.g. $\lambda = (5,3,3,1) \in \Par_5$.

\begin{lemma}
    $\lambda^\pi = w_0(\lambda).$ Hence the $\lambda^\pi$ with ambient length $n$ provide another parametrization of the irreducible polynomial representations of $\GL_n$.
\end{lemma}

\begin{proof}
    Write $\lambda = \lambda_1 \epsilon_1 + \dots + \lambda_n \epsilon_n$. $w_0(\lambda) = \lambda_1 \epsilon_n + \lambda_2 \epsilon_{n-1} + \cdots + \lambda_n \epsilon_1$ which corresponds to  $(\lambda_n, \lambda_{n-1}, \cdots , \lambda_1)$.
\end{proof}

Outside of the polynomial representations of $\GL_n$, we do not need a lowest weight parametrization and use the common parameterization of $\widehat{K}$ given in the following table:
\begin{center}
\begin{tabular}{c|c}
    $K$ & $\widehat{K}$\\
    \hline
    $\Or_n$ & \{$\lambda$ : first two columns of $\lambda$ contain $\leq n$ boxes\} \\
    $\GL_n$ & \{ $\lambda = (\lambda^+, \lambda^-)$ : $\ell(\lambda^+) + \ell(\lambda^-) \leq n \}$\\
    $\Sp_{2n}$ & \{ $\lambda$ : $\ell(\lambda) \leq n$ \}\\
\end{tabular}
\end{center}

If $K=\GL_n$ and $\lambda^+=(\lambda_1^+,\dots,\lambda_i^+)$, $\lambda^-=(\lambda_1^-,\dots,\lambda_j^-)$, then $\pi^\lambda_{\GL_n}$ is the irreducible representation with highest weight
\(
(\lambda_1^+,\dots,\lambda_i^+,0,\dots,0,-\lambda_j^-,\dots,-\lambda_1^-),
\)
see \cite[Proposition 2.1]{STEMBRIDGE1987}. In particular, the polynomial representations are given by $\lambda = (\lambda^+, \varnothing)$. In this case, we will suppress the $\varnothing$ and think of $\lambda$ as a Young diagram as above. For $\Sp_{2n}$, $\pi^\lambda_{\Sp_{2n}}$ is just the irreducible representation with highest weight $\lambda$ in standard coordinates. For $K = \Or_n$, which is not connected, highest weights do not directly apply, see \cite[p.~438-439]{GoodmanWallach2009}. Recall, the irreducible representations of $\Or_n$ and $\Sp_{2n}$ are self-dual, while the irreducible representations of $\GL_n$ are not. For $\lambda = (\lambda^+, \lambda^-)$, the diagram associated to $(\pi^\lambda_{\GL_n})^*$ is $\lambda^* = (\lambda^-, \lambda^+)$. Notice, for a polynomial representation given by $\lambda = (\lambda^+, \varnothing)$, $\lambda^* = (\varnothing, \lambda^+)$ is not polynomial unless $\lambda^+ = \varnothing$.

Next we want to fill these diagrams in a way that will allow us to easily compute the $M$-type multiplicities in $\Res^K_M(\pi^\lambda_K)$. The usual way to do this for polynomial irreducible representations of $\GL_n$ in the context of the symmetric pair $(\GL_n^2,\GL_n)$ is by semistandard Young tableaux. Recall a \textit{(semistandard Young) tableau} is a filling of $\lambda$ that is
\begin{enumerate}
    \item weakly increasing across each row,
    \item strictly increasing down each column.
\end{enumerate}
More generally, a \textit{skew tableau} on $\lambda/\mu$ is a filling satisfying (1) and (2). Denote the set of skew tableaux on shape $\lambda/\mu$ with fillings in the alphabet $[n]$ by $\SST^{\lambda / \mu}_n$ and the Young tableaux on $\lambda$ with fillings in $[n]$ by $\SST^\lambda_n$. For every $\lambda$ with ambient length $n$, let $s^\lambda_n$ be the associated \textit{Schur polynomial}, $$s^\lambda_n(x_1, \dots , x_n) = \sum_{T \in \SST^\lambda_n} x^{\wt(T)}.$$
Here $\wt(T) = w_1 \epsilon_1 + \cdots + w_n \epsilon_n$ where $w_i = \text{number of $i$'s in $T$}$, and $x^{\wt(T)} = \prod_{i=1}^n x_i^{w_i}$. The Schur polynomial $s^\lambda_n$ is the character of $\pi^\lambda_{\GL_n}$. We continue to denote the diagonal maximal torus in $\GL_n$ by $H$. Notice that the coefficient of a monomial $x^\eta$ in $s_n^\lambda$ is the multiplicity of the torus representation $\pi_H^{\eta}$ in $\Res_H^{\GL_n}(\pi^\lambda_{\GL_n})$, i.e. the multiplicity of the $\eta$-weight space. For the symmetric pair $(\GL_n^2, \GL_n)$, $M= H$ and the semistandard tableaux are exactly capturing the multiplicities in the branching problem from $K$ to $M$ for polynomial representations.
Definition \ref{def:K-tableaux} characterizes the fillings that have the same property in the three cases we consider. Before presenting the definition, we give a parametrization of the $\widehat{M}$ in the relevant cases. 

The irreducible representations $\pi^{\delta}_M$ of $M$ are labeled by the following $n$-tuples $\delta$:
\begingroup
\renewcommand*{\arraystretch}{1.5}
\begin{equation*}
    \begin{array}{l|l|l}
        K & M & \delta \in \widehat{M} \\ \hline

        \Or_n & \Or_1^n & \delta \in \{0,1\}^n\\
        \GL_n & \GL_1^n & \delta \in \mathbb{Z}^n\\
        \Sp_{2n} & \Sp_2^n & \delta \in \mathbb{Z}_{\geq 0}^n     
    \end{array}
\end{equation*}
\endgroup

In the case of $\Or_1$, 0 denotes the trivial representation and 1 denotes the alternating representation. $\GL_1^n$ is just the diagonal maximal torus in $\GL_n$ and $\delta = (\delta_1, \dots \delta_n) \in \Z^n$ denotes the character with weight $\delta_1 \epsilon_1 + \cdots + \delta_n \epsilon_n$. For $\Sp_2$, an integer $m \in \Z_{\geq 0}$ parametrizes the irreducible representation with highest weight $m\varpi_1 = m\epsilon_1$.

\begin{defn}\label{def:K-tableaux}
$K$-tableaux and their $M$-weights
    
    \begin{enumerate}
        \item A tableau $T \in \SST^\lambda_n$ is an \emph{$\Or_n$-tableau} if, for all $1 \leq i \leq n$, we have
        \[
        \#\{ \text{boxes in first two columns of $T$ with entry $\leq i$}\} \leq i.
        \]
        Let $\mathcal{T}_{\Or_n}^\lambda$ denote the set of $\Or_n$-tableaux on shape $\lambda$. We define $\mwt_{\Or_n}(T)$ to be the vector $(w_1, \ldots, w_n) \in \widehat{M} = \{0,1\}^n$ such that 
        \[
        w_i = \#\{\text{boxes in $T$ with entry $i$}\} \: {\rm mod} \: 2.
        \]
        
        \item A pair $(T^+, T^-) \in \SST^{\lambda^+}_n \times \SST^{\lambda^-}_n$ is a \emph{$\GL_n$-tableau} if, for all $1 \leq i \leq n$, we have
        \[
        \#\{\text{boxes in first column of $T^+$ or $T^-$ with entry $\leq i$}\} \leq i.
        \]
        Let $\mathcal{T}_{\GL_n}^\lambda$ denote the set of $\GL_n$-tableaux on shape $\lambda$. We define $\mwt_{\GL_n}(T^+, T^-)$ to be the vector $(w_1, \ldots, w_n) \in \widehat{M} = \mathbb{Z}^n$ such that
        \[
        w_i = \#\{\text{boxes in $T^+$ with entry $i$}\} - \#\{\text{boxes in $T^-$ with entry $i$}\}.
        \]

        \item A tableau $T \in \SST^\lambda_{2n}$ is an \emph{$\Sp_{2n}$-tableau} if, for all $1 \leq i \leq n$, we have
        \[
        \#\{\text{boxes in first column of $T$ with entry $\leq 2i$}\} \leq i.
        \]
        Moreover, we say that $T$ is an \emph{$\Sp_{2n}$-ballot tableau} if, as one reads the word of $T$ right to left (word is defined in Section \ref{subsec:tableaux and knuth equiv} below), the number of $2i$'s never exceeds the number of $2i-1$'s, for each $1 \leq i \leq n$.
        
        Let $\mathcal{T}_{\Sp_{2n}}^{\lambda, H}$ denote the set of $\Sp_{2n}$-tableaux on shape $\lambda$ and $\mathcal{T}_{\Sp_{2n}}^{\lambda}$ be the set of $\Sp_{2n}$-ballot tableaux on $\lambda$. If $T$ is an $\Sp_{2n}$-ballot tableau, then we define $\mwt_{\Sp_{2n}}(T)$ to be the vector $(w_1, \ldots, w_n) \in \widehat{M} = \mathbb{Z}_{\geq 0}^n$ such that
        \[
        w_i = \#\{\text{boxes in $T$ with entry $2i-1$}\} - \#\{\text{boxes in $T$ with entry $2i$}\}.
        \]
    \end{enumerate}
\end{defn}

We call $\mwt_K(T)$ the \textit{$M$-weight} of $T$, meaning the parameter in $\widehat M$ of the $M$-type associated to $T$; in the torus case this is literally a weight. 

Definition~\ref{def:K-tableaux} can be restated as follows: $T$ is a $K$-tableau if and only if, upon restricting $T$ to those entries $\leq i$ (or $\leq 2i$ for $\Sp_{2n}$), the resulting shape lies in the set $\widehat{\Or}_i$ or $\widehat{\GL}_i$ or $\widehat{\Sp}_{2i}$ for all $i$. Call the shape obtained by retaining only those boxes with entries $\leq i$ the shape \emph{supported} on $[i]$. Then Definition~\ref{def:K-tableaux} says that $T$ is a $K$-tableau if and only if, for every $i$, the shape supported on $[i]$ (or on $[2i]$ in the symplectic case) again satisfies the defining condition for the corresponding smaller-rank group. Also, observe that the set $\mathcal{T}^\lambda_K$ is empty unless $\lambda \in \widehat{K}$. Finally, an equivalent characterization of the $\Sp_{2n}$-tableaux is that the entries in row $i$ are at least $2i-1$ for $i \in [n]$.

Note, the order preserving bijection $i \mapsto 2i-1$, $\overline{\iota} \mapsto 2i$ gives a bijection between the $\Sp_{2n}$-tableaux used in \cite{CEFW2025} and those used here.

The tableaux in $\mathcal{T}^\lambda_K$ parameterize the $M$-types occurring in $\Res^K_M(\pi^\lambda_K)$, with multiplicity. In the case of $\Or_n$ and $\GL_n$, $M$ is abelian so all irreducible representations are one-dimensional and $|\mathcal{T}^\lambda_K| = \dim \pi^\lambda_K$. For $\Sp_{2n}$, $M$ is a product of $\Sp_2$. Here the  $M$ irreducible representations are not one-dimensional unless $\delta = (0, \dots , 0)$, and thus $|\mathcal{T}^\lambda_{\Sp_{2n}}| \neq \dim \pi^\lambda_{\Sp_{2n}}$ in general. The full set of $\Sp_{2n}$-tableaux, $\mathcal{T}^{\lambda,H}_{\Sp_{2n}}$, comes from branching down to a maximal torus $H$, hence the notation. They parameterize a weight basis for $\pi^\lambda_{\Sp_{2n}}$. The ballot tableaux come from branching to $M$. They parameterize the highest weights of the $M = \Sp_{2}^n$ irreducible representations occurring in $\Res^{\Sp_{2n}}_{\Sp_{2}^n} (\pi^\lambda_{\Sp_{2n}})$.

\begin{example}
    We illustrate the difference between $\Sp_{4}$-tableaux and $\Sp_4$-ballot tableaux for $\lambda = (2,1) = \ydiagram{2,1}$.

    There are 4 $\Sp_4$-ballot tableaux and 16 $\Sp_4$-tableaux. The ballot tableaux are:
    $$\ytableaushort{1 1,3}, \ytableaushort{1 3, 3}, \ytableaushort{1 3 ,4}, \ytableaushort{3 3, 4}.$$

    These are the highest weights for four $\Sp_2^2$ irreducible representations with highest weights $2 \varpi_1 + \varpi_2$, $\varpi_1+ 2 \varpi_2$, $\varpi_1$, and $\varpi_2$. The first two are 6-dimensional, and the last two are 2-dimensional. It is easy to find the remaining $\Sp_4$-tableaux.
\end{example}

\begin{defn} \label{def:M-weight delta}
    Denote the subset of $\mathcal{T}^\lambda_K$ with $M$-weight $\delta$ by $(\mathcal{T}^\lambda_K)_\delta$. 
\end{defn}

The following result of \cite{CEFW2025} formalizes the discussion above. Recall that for $\lambda\in\widehat K$ and $\delta\in\widehat M$,
\[
b^\lambda_\delta=\mult(\pi^\delta_M,\pi^\lambda_K).
\]
\begin{thm}[{\cite[Theorem 4.2]{CEFW2025}}]
    \label{thm:K to M}

    Let $\lambda \in \widehat{K}$ and $\delta \in \widehat{M}$.
\[
    b^\lambda_{\delta} = \# (\mathcal{T}^\lambda_K)_\delta .
\]

\end{thm}
See Section \ref{sec:graded-mults} for examples. 

\subsection{Principal series and $M$-polynomials}
\label{sec:principal-series-M-polynomials}

A useful way to organize this information is through what we call the \textit{$M$-polynomial}. For $\nu \in \widehat{K}$, define
\[M_{(G,K)}^\nu \coloneqq \sum_{T \in \mathcal{T}_K^\nu} \chi_M^{\mwt_K(T)},\]
where $\chi_M^{\mwt_K(T)}$ is the character of the irreducible $M$-representation corresponding to $\mwt_K(T)$. Since the tableaux in $\mathcal{T}_K^\nu$ index the irreducible $M$-summands of $\Res_M^K(\pi_K^\nu)$, counted with multiplicity, we have
\[M_{(G,K)}^\nu = \sum_{\delta \in \widehat{M}} b_\delta^\nu \chi_M^\delta.\]
Thus, $M_{(G,K)}^\nu$ is the $M$-character of $\Res_M^K(\pi_K^\nu)$, and the coefficient of $\chi_M^\delta$ is the branching multiplicity $b_\delta^\nu$ computed in Theorem~\ref{thm:K to M}.

By Frobenius reciprocity, $b_\delta^\nu$ is also the multiplicity of $\pi_K^\nu$ in any principal series representation induced from $\pi_M^\delta$. Therefore, $M_{(G,K)}^\nu$ simultaneously records the multiplicity of the fixed $K$-type $\pi_K^\nu$ in all principal series representations of $G_{\R}$. In the case corresponding to the symmetric pair $(\GL_n^2,\GL_n)$, we have
\[ M_{(\GL_n^2,\GL_n)}^\nu=s_n^\nu, \]
the rational Schur polynomial; see \cite{STEMBRIDGE1987}. In this case, $M$ is a torus, so its irreducible characters may be written as monomials $\chi_M^\delta=x^\delta$. After a determinant shift, the coefficient of $x^\delta$ is an ordinary Kostka number.

Summing over all $K$-types gives the formal \textit{$M$-series}
\[M_{(G,K)} \coloneqq  \sum_{\nu \in \widehat{K}} M_{(G,K)}^\nu \chi_K^\nu.\]
Equivalently,
\[M_{(G,K)} = \sum_{\delta \in \widehat{M}} \chi_M^\delta \sum_{\nu \in \widehat{K}} b_\delta^\nu \chi_K^\nu.\]
The coefficient of $\chi_K^\nu$ is the $M$-polynomial $M_{(G,K)}^\nu$, while the coefficient of $\chi_M^\delta$ is the ungraded $K$-character of a principal series representation induced from $\pi_M^\delta$. Thus, this single formal character simultaneously packages the compact pictures of all principal series representations.

This interpretation can be organized using the regular functions $\C[K]$. Twisting the usual left regular action, as in \cite[Section~4.2.3 and Corollary~5.6.8]{GoodmanWallach2009}, the algebraic Peter-Weyl theorem takes the form
\[\C[K] \cong \bigoplus_{\nu \in \widehat{K}} \pi_K^\nu \boxtimes \pi_K^\nu.\]
After restricting the right $K$-action to $M$, this becomes
\[\C[K] \cong \bigoplus_{\nu \in \widehat{K}} \pi_K^\nu \boxtimes \Res_M^K(\pi_K^\nu).\]
Reorganizing the right-hand side by $M$-type and applying Frobenius reciprocity gives
\[\C[K] \cong \bigoplus_{\delta \in \widehat{M}} \Ind_M^K(\pi_M^\delta) \boxtimes \pi_M^\delta.\]
Consequently, $M_{(G,K)}$ is the $K \times M$-character of $\C[K]$ for this choice of $K \times M$-action.

The preceding discussion suggests considering a \textit{graded
$M$-polynomial} $M_{(G,K)}^\nu(q)$ and a \textit{graded $M$-series}
$M_{(G,K)}(q)$. These are $q$-analogues in which the exponent of $q$
records the degree in which $\pi_K^\nu$ occurs in the associated graded
of the Hodge filtration of a principal series representation. At the
level of formal characters, this amounts to seeking a graded
enhancement of the $K \times M$-character of $\C[K]$.

Let $I(\delta)$ denote the tempered principal series induced from
$\pi_M^\delta$ with real infinitesimal character and $\gr I(\delta)$ the associated graded of the Hodge filtration. Define
\[h_{(G,K)}^{\nu,\delta}(q) \coloneqq \mult_q \left( \pi_K^\nu, \gr I(\delta) \right).\]
Then
\[h_{(G,K)}^{\nu,\delta}(1)=b_\delta^\nu,\]
so $h_{(G,K)}^{\nu,\delta}(q)$ is a $q$-analogue of the branching multiplicity $b_\delta^\nu$. We define
\[M_{(G,K)}^\nu(q) \coloneqq \sum_{\delta \in \widehat{M}} h_{(G,K)}^{\nu,\delta}(q)\chi_M^\delta\]
and
\begin{align*}
    M_{(G,K)}(q) & \coloneqq \sum_{\nu \in \widehat{K}} M^\nu_{(G,K)}(q) \chi^\nu_K\\
& =\sum_{\delta \in \widehat{M}}
\chi_M^\delta
\sum_{\nu \in \widehat{K}}
h_{(G,K)}^{\nu,\delta}(q)\chi_K^\nu.
\end{align*}
At $q=1$, these specialize to the ungraded $M$-polynomial and $M$-series. The coefficient of $\chi^\nu_K$ is the graded $M$-polynomial, while the coefficient of $\chi^\delta_M$ is the Hodge $K$-character of $I(\delta)$.

Theorem~\ref{thm:intro Hodge filtration} implies that, for $G_{\R}=\GL(n,\R)$,
\[h_{(G,K)}^{\nu,0}(q) = m_{(G,K)}^{\nu,0}(q).\]
One of the main results of this paper computes this polynomial by defining a grading statistic on the tableaux in $(\mathcal{T}_K^\nu)_0$.

More generally, since the tableaux in $\mathcal{T}_K^\nu$ index the irreducible $M$-summands of $\Res_M^K(\pi_K^\nu)$, counted with multiplicity, it is natural to ask whether there is a statistic $d(T)$, or a suitable extension of the statistic defined in this paper, such that
\[h^{\nu, \delta}_{(G,K)}(q) = \sum_{T \in (\mathcal{T}_K^\nu)_\delta} q^{d(T)}.\]
Equivalently, for each fixed $\nu$, one would have
\[M_{(G,K)}^\nu(q) = \sum_{T \in \mathcal{T}_K^\nu} q^{d(T)} \chi_M^{\mwt_K(T)}. \]

From this perspective, the present paper's computation of $m^{\nu,0}_{(G,K)}(q)$ for the three families is aimed at determining the coefficient of the trivial $M$-character $\chi_M^0=1$ in the graded $M$-series. For $G_{\R}=\GL(n,\R)$, this coefficient has the Hodge-theoretic interpretation given above. For the nonsplit families, it would be interesting to determine whether the same combinatorics admits a corresponding Hodge-theoretic interpretation. It would also be interesting to compute the remaining coefficients $h_{(G,K)}^{\nu,\delta}(q)$ in the tempered real-infinitesimal-character case and, more generally, to study the corresponding Hodge-graded multiplicities as the continuous parameter varies.

\subsection{The ring of invariants \texorpdfstring{$\C[\p]^K$}{C[p]K}} \label{sec:ring-of-invariants}
The well-known work of Chevalley establishes the structure of $\C[\p]^K$. In all cases, this is a polynomial algebra, i.e. generated by algebraically independent elements. In the three cases we consider, $\C[\p]^K \cong \C[u_1, \dots , u_n]$ where $u_i(X) = \tr(X^i)$. Hence, for $(\GL_n, \Or_n)$, $(\GL_n^2, \GL_n)$, and $(\GL_{2n}, \Sp_{2n})$ the graded character of the invariants is,
$$\chr_q(\C[\p]^K) = \frac{1}{\prod_{i=1}^n (1 - q^i)}.$$
For more detail, see \cite[Section 12.4.2]{GoodmanWallach2009}. 

\subsection{\texorpdfstring{$\C[\p]$}{C[p]} as a multiplicity-free space} \label{sec:mult-free-spaces}
In the cases we consider, there are well-known multiplicity-free decompositions of $\C[\p]$ as a $G$-representation whose restriction to $K$ is the given representation of $K$ on $\C[\p]$, see \cite[Section 5.1]{Howe1995}.

\begin{thm} \label{thm:mult-free spaces}
 We have the following multiplicity-free decompositions under the action of $G$,

\begin{center}
\begin{tabular}{c|l|c}
\multicolumn{1}{c}{$G$} & \multicolumn{1}{c}{$\C[\p]$} & \multicolumn{1}{c}{degree}\\
    \hline
    $\GL_n$ & $\displaystyle \C[\SM_n] = \bigoplus_{\lambda \in \Par_n^{(2)}} \pi_{\GL_n}^\lambda$ & $|\lambda|/2$\\
    $\GL_n^2$ & $\displaystyle \C[\M_n] = \bigoplus_{\lambda \in \Par_n} \pi^{[\lambda^*, \lambda]}_{\GL_n^2}$ & $|\lambda|$\\
    $\GL_{2n}$ & $\displaystyle \C[\AM_{2n}] = \bigoplus_{\lambda \in \Par_{2n}^{(1,1)}} \pi_{\GL_{2n}}^\lambda$ & $|\lambda|/2$
\end{tabular}
\end{center}

\noindent where $\pi^{[\lambda^*, \lambda]}_{\GL_n^2}$ denotes the $\GL_n^2$ irreducible representation $(\pi_{\GL_n}^\lambda)^* \boxtimes \pi_{\GL_n}^\lambda$. The representation $\pi^\lambda_G$ occurs in degree shown in the third column.
\end{thm}

\begin{proof}
    These are theorems 5.6.7, 5.7.3, and 5.7.5 in \cite{GoodmanWallach2009}. See also \cite{Howe1995}.
\end{proof}

\section{Crystals of reductive Lie algebras} \label{sec:crystals-ofreductive}

We give a short review of the facts we need about crystals. For more, see \cite{BumpSchilling2017, HongKang2002}. Our conventions follow \cite{BumpSchilling2017}. In this paper, the point of crystals is to provide a combinatorial model for finite-dimensional representations of complex reductive Lie algebras. Crystals originally arose independently as limits/crystallizations of representations of quantum groups on the one hand and out of standard monomial theory and Littelmann paths on the other. Bump and Schilling \cite{BumpSchilling2017} provide a purely combinatorial treatment relying on character theory to make the connection with Lie algebras. Much of this can be done in the generality of Kac-Moody algebras as in \cite{HongKang2002} but we will not need that here. Let $\g$ be a complex reductive Lie algebra and carry forward the notation from Section \ref{sec:reductive Lie algebras}.

In Section \ref{subsec:defs}, we define crystals and other basic notions. In Section \ref{subsec:tensor products}, we discuss tensor products of crystals and characterize highest and lowest weights in a tensor product. In Section \ref{subsec:normal crystals}, we briefly review normal crystals and how they mirror the finite-dimensional representation theory of Lie algebras. Finally, in Section \ref{subsec:crystal_anti-iso}, we introduce crystal anti-isomorphisms, a slight generalization of crystal involutions, and show how they interact with tensor products. Our exposition is for crystals of reductive Lie algebras, not just semisimple Lie algebras. We try to put highest and lowest weights on equal footing and use anti-isomorphisms to move between them.

\subsection{Definitions} \label{subsec:defs}
We begin with the definition of a finite-type crystal. Any crystal mentioned in this paper is finite-type. In fact, we will only work with seminormal crystals as defined below. 

\begin{defn}\label{def:crystals}
    A \textit{(finite-type) crystal} of $\g$ is a nonempty set $\B$ together with maps
    \begin{align*}
        e_i, f_i & : \B \to \B \sqcup \{0 \},\\
        \varepsilon_i, \varphi_i & : \B \to \Z, \\
        \wt & : \B \to P(\g),
    \end{align*}
    for all $i \in [r]$. The maps must satisfy the following conditions:
    \begin{enumerate}
        \item If $x,y \in \B$, then $e_i(x) = y$ if and only if $f_i(y) = x$. In this case,
        \begin{itemize}
            \item $\wt(e_ix) = \wt(x) + \alpha_i$ and hence $\wt(f_i y ) = \wt(y) - \alpha_i$.
            \item $\varepsilon_i(e_i x) = \varepsilon_i(x) - 1$ and hence $\varepsilon_i(f_i y) = \varepsilon_i(y) + 1$.
            \item $\varphi_i(e_ix) = \varphi_i(x) + 1$ and hence $\varphi_i(f_iy) = \varphi_i(y) - 1$.
        \end{itemize}

        \item $\varphi_i(x) = \langle \wt(x), h_i \rangle + \varepsilon_i(x)$ for all $x \in \B$ and $i \in [r]$.
    \end{enumerate}
    
\end{defn}

The map $\wt$ is the \textit{weight map} and the maps $e_i, f_i$ are the \textit{Kashiwara operators}, \textit{(raising and lowering operators)}. $\B$ is called \textit{seminormal} if 
$$\varphi_i(x) = \max \{ k \in \Z_{\geq 0} : f_i^k(x) \neq 0 \} \text{ and } \varepsilon_i(x) = \max \{k \in \Z_{\geq 0} : e_i^k(x) \neq 0 \},$$
\noindent for all $i \in [r]$. All crystals we consider will be seminormal. The set $\B$ should be thought of as a weight basis for a representation of $\g$, with $e_i$ and $f_i$ corresponding to the raising and lowering operators in the $\sll_2$-triple attached to the simple root $\alpha_i$; condition (2) is a natural relation present in any $\sll_2$-string.

For $x \in \B$ let
\begin{align*}
\varphi(x) &= \sum_{i =1}^r\varphi_i(x)\varpi_i,   &   \varepsilon(x) &= \sum_{i = 1}^r\varepsilon_i(x)\varpi_i.    
\end{align*}

Notice $\varphi(x), \varepsilon(x) \in P_{++}(\g')$ are dominant weights for the semisimple Lie algebra $\g'$. Starting in Section \ref{sect:gln-crystals}, all crystals will be $\gl_n$-crystals. At times, it will be important to separate the weight coming from the semisimple $\sll_n$ from the weight coming from the action of the center. To this end, we make the following general definitions analogous to the definitions made in Section \ref{sec:reductive Lie algebras} for reductive Lie algebras. Let $x \in \B$. Define the \textit{semisimple weight map}
$$\wt_{ss}(x) = \sum_{i = 1}^r \langle \wt(x), h_i \rangle \varpi_i$$
and the \textit{central weight map}
$$\wt_{c}(x) = \sum_{i = r+1}^n \langle \wt(x), h_i \rangle \varpi_i.$$
\begin{lemma}
    If $\B$ is a crystal for a complex reductive Lie algebra $\g$, then $\wt(x) = \wt_{ss}(x) + \wt_c(x)$ for all $x \in \B$. If $\g$ is semisimple, $\wt(x) = \wt_{ss}(x)$.
\end{lemma}

\begin{proof}
    We are simply extracting the coefficients of $\wt(x)$ with respect to the basis $\varpi_1,\dots,\varpi_n$. If $\g$ is semisimple, then $r=n$ and hence $\wt_c(x)=0$.
\end{proof}

The \textit{crystal graph} of $\B$ is the directed graph with vertex set $\B$ and labeled edges $x \xrightarrow{i} y$ whenever $f_i(x)=y$. An element $x \in \B$ such that $e_i(x) = 0$ for all $i \in [r]$ is called a \textit{highest-weight element} and in this case the weight of $x$ is called a \textit{highest weight}. Similarly, an element $x \in \B$ such that $f_i(x) = 0$ for all $i \in [r]$ is called a \textit{lowest-weight element} and in this case the weight of $x$ is called a \textit{lowest weight}.

If $\B$ and $\mathcal{C}$ are two crystals associated with the same Lie algebra $\g$, a \textit{crystal morphism} is a map $\psi:\B\to \mathcal{C}\sqcup\{0\}$ such that
\begin{enumerate}
    \item if $x\in \B$ and $\psi(x)\in \mathcal{C}$, then
    \begin{enumerate}
        \item $\wt(\psi(x))=\wt(x)$,
        \item $\varepsilon_i(\psi(x))=\varepsilon_i(x)$ for all $i\in [r]$,
        \item $\varphi_i(\psi(x))=\varphi_i(x)$ for all $i\in [r]$;
    \end{enumerate}
    \item if $x,e_ix\in \B$ and $\psi(x),\psi(e_ix)\in \mathcal{C}$, then $\psi(e_ix)=e_i\psi(x)$;
    \item if $x,f_ix\in \B$ and $\psi(x),\psi(f_ix)\in \mathcal{C}$, then $\psi(f_ix)=f_i\psi(x)$.
\end{enumerate}
A morphism $\psi$ is \textit{strict} if it commutes with $e_i$ and $f_i$ for all $i\in [r]$, and it is a \textit{crystal isomorphism} if $\psi$ is a bijection and $\psi^{-1}$ is also a crystal morphism.

\subsection{Tensor products of crystals} \label{subsec:tensor products}
If $\B$ and $\mathcal{C}$ are two crystals associated with the same Lie algebra $\g$, the \textit{tensor product} $\B \otimes \mathcal{C}$ is defined to be the Cartesian product of $\B$ and $\mathcal{C}$ as a set, with ordered pairs denoted $x \otimes y$. The crystal structure is given by:
\begin{align*}
    \wt(x \otimes y) & = \wt(x) + \wt(y),\\
    \varphi_i(x \otimes y) & = \varphi_i(x) + \max(0, \varphi_i(y) - \varepsilon_i(x)),\\
    \varepsilon_i(x \otimes y) &= \varepsilon_i(y) + \max ( 0, \varepsilon_i(x) - \varphi_i(y)), \\
    f_i(x \otimes y ) & = \begin{cases}
        f_i(x) \otimes y \quad \text{if } \varepsilon_i(x) \geq \varphi_i(y), \\
        x \otimes f_i(y) \quad \text{if } \varepsilon_i(x) < \varphi_i(y),\\
    \end{cases}\\
    e_i(x \otimes y ) & = \begin{cases}
        e_i(x) \otimes y \quad \text{if } \varepsilon_i(x) > \varphi_i(y), \\
        x \otimes e_i(y) \quad \text{if } \varepsilon_i(x) \leq \varphi_i(y).\\
    \end{cases}
\end{align*}

This is known as the \textit{tensor product rule}. We provide a characterization of highest and lowest-weight elements in tensor products. 

\begin{prop}\label{prop-highest/lowest-in-tensor-product}
    Let $\B$ and $\mathcal{C}$ be seminormal crystals.
    
    \begin{enumerate}
        \item $x \otimes y \in \B \otimes \mathcal{C}$ is a highest-weight element if and only if $y$ is a highest-weight element in $ \mathcal{C}$ and $\varepsilon(x) \leq \wt_{ss}(y)$.
        \item $x \otimes y \in \B \otimes \mathcal{C}$ is a lowest-weight element if and only if $x$ is a lowest-weight element in $\mathcal{B}$ and $-\wt_{ss}(x) \geq \varphi(y)$.
    \end{enumerate}
\end{prop}

\begin{proof}
    (1) By the seminormal assumption, $x \otimes y$ is a highest-weight element if and only if $\varepsilon_i(x \otimes y ) = 0$ for all $i \in [r]$. By the tensor product rule, this occurs if and only if $\varepsilon_i(y) =0$ and $\varepsilon_i(x) \leq \varphi_i(y)$ for all $i \in [r]$. Now $\varphi_i(y) = \langle \wt(y), h_i \rangle + \varepsilon_i(y) = \langle \wt(y), h_i \rangle$ by the crystal axioms and the fact that $y$ is a highest-weight element. Thus $x \otimes y$ is a highest-weight element if and only if $y$ is a highest-weight element and $\varepsilon_i(x) \leq \langle \wt(y), h_i \rangle$ for all $i \in [r]$. Now, $\varepsilon_i(x) \leq \langle \wt(y), h_i \rangle$ for all $i \in [r]$ if and only if $\varepsilon(x) \leq \wt_{ss}(y)$ and the result follows. 

    (2) The idea is the same as (1). $x \otimes y$ is a lowest-weight element if and only if $\varphi_i(x \otimes y) =0$ for all $i \in [r]$ if and only if $\varphi_i(x) =0$ and $\varphi_i(y) \leq \varepsilon_i(x)$. With $x$ a lowest-weight element, $\varepsilon_i(x) = - \langle \wt(x), h_i \rangle$ and the result follows.
\end{proof}

We have the following generalization of the tensor product rule to $k$-fold tensor products,

\begin{lemma}[\cite{BumpSchilling2017} Lemma 2.33] \label{lemma:multiple_tensors}
    Let $x_1 \otimes \cdots \otimes x_k \in \B_1 \otimes \cdots \otimes \B_k$ for $\g$-crystals $\B_1, \dots , \B_k$. Then,
    $$\varphi_i(x_1 \otimes \cdots \otimes x_k) = \max_{j=1}^k \bigg ( \sum_{h =1}^j \varphi_i(x_h) - \sum_{h=1}^{j-1} \varepsilon_i(x_h) \bigg ),$$
    \noindent and if $j$ is the first value where the maximum is attained, then 
    $$f_i(x_1 \otimes \cdots \otimes x_k) = x_1 \otimes \cdots \otimes f_i(x_j) \otimes \cdots \otimes x_k.$$
    \noindent Similarly, for $x_k \otimes \cdots \otimes x_1 \in \B_k \otimes \cdots \otimes \B_1$,
    $$\varepsilon_i(x_k \otimes \cdots \otimes x_1) = \max_{j=1}^k \bigg ( \sum_{h =1}^j \varepsilon_i(x_h) - \sum_{h=1}^{j-1} \varphi_i(x_h) \bigg ),$$
    \noindent and if $j$ is the first value where the maximum is attained, then 
    $$e_i(x_k \otimes \cdots \otimes x_1) = x_k \otimes \cdots \otimes e_i(x_j) \otimes \cdots \otimes x_1.$$ 
\end{lemma}

Let $\B_1$ and $\B_2$ be seminormal crystals of complex reductive Lie algebras $\g_1$ and  $\g_2$ with semisimple ranks $r_1$ and $r_2$. Let $e_i^1, f_i^1, \varphi_i^1, \varepsilon_i^1, \wt^1$ and $e_i^2, f_i^2, \varphi_i^2, \varepsilon_i^2, \wt^2$ denote their operators. Reindex the operators of $\B_2$ by $i \mapsto i+r_1$ so that $i \in \{r_1+1,\dots,r_1+r_2\}$. We continue to denote the reindexed operators by $e_i^2, f_i^2, \varphi_i^2$, and $\varepsilon_i^2$. Define the \textit{outer tensor product} $\B_1 \boxtimes \B_2$ as the Cartesian product of $\B_1$ and $\B_2$ as a set, with ordered pairs denoted $x_1 \boxtimes x_2$. We equip $\B_1 \boxtimes \B_2$ with the structure of a $\g_1 \oplus \g_2$-crystal. For $i \in [r_1 + r_2]$, the maps behave as follows:

\begin{align*}
    \wt(x_1 \boxtimes x_2) & = \wt^1(x_1) + \wt^2(x_2),\\
    \varphi_i(x_1 \boxtimes x_2) & = \begin{cases}
        \varphi^1_i(x_1)  \quad \text{if } i \leq r_1, \\
        \varphi^2_{i}(x_2)  \quad \text{if } i > r_1,\\
    \end{cases}\\
    \varepsilon_i(x_1 \boxtimes x_2) & = \begin{cases}
        \varepsilon^1_i(x_1)  \quad \text{if } i \leq r_1, \\
        \varepsilon^2_{i}(x_2)  \quad \text{if } i > r_1,\\
    \end{cases}\\
    f_i(x_1 \boxtimes x_2 ) & = \begin{cases}
        f^1_i(x_1) \boxtimes x_2 \quad \text{if } i \leq r_1, \\
        x_1 \boxtimes f_i^2(x_2) \quad \text{if } i > r_1, \\
    \end{cases}\\
    e_i(x_1 \boxtimes x_2 ) & = \begin{cases}
        e^1_i(x_1) \boxtimes x_2 \quad \text{if } i \leq r_1, \\
        x_1 \boxtimes e_i^2(x_2) \quad \text{if } i > r_1. \\
    \end{cases}\\
\end{align*}

Again, we characterize the highest and lowest-weight elements.
\begin{prop}\label{prop:exterior_tensor}
     Let $\B_1$ and $\B_2$ be crystals as above.
    
    \begin{enumerate}
        \item $x_1 \boxtimes x_2 \in \B_1 \boxtimes \B_2$ is a highest-weight element if and only if $x_1$ is a highest-weight element in $\B_1$ and $x_2$ is a highest-weight element in $\B_2$.
        \item $x_1 \boxtimes x_2 \in \B_1 \boxtimes \B_2$ is a lowest-weight element if and only if $x_1$ is a lowest-weight element in $\B_1$ and $x_2$ is a lowest-weight element in $\B_2$.
    \end{enumerate}
\end{prop}

\begin{proof}
    By definition, the operators indexed by $1,\dots,r_1$ act only on the first factor, while those indexed by $r_1+1,\dots,r_1+r_2$ act only on the second factor. Hence $x_1 \boxtimes x_2$ is highest weight (respectively lowest weight) if and only if each factor is highest weight (respectively lowest weight).
\end{proof}

\subsection{Normal crystals} \label{subsec:normal crystals}
Recall that we want to use crystals as a combinatorial model for finite-dimensional representations of complex reductive Lie algebras. The category of seminormal crystals is too large for this purpose. The connection with representation theory is obtained by restricting to normal crystals. There are a number of equivalent ways to define normal crystals and make this connection: via quantum groups \cite{Kashiwara1995}, Littelmann paths \cite{Littelmann1995}, and combinatorially \cite{BumpSchilling2017}. 

Once normal crystals are defined, the category is semisimple in the sense that crystals decompose into a disjoint union of connected subcrystals. The connected subcrystals correspond exactly to finite-dimensional irreducible representations of the associated complex reductive Lie algebra, and the disjoint union of connected crystals corresponds to the direct sum of irreducible representations. Restriction to Levi subalgebras and tensor products of crystals produce the same multiplicities as for finite-dimensional representations of the associated Lie algebra. 

The character of a crystal is defined as $\sum_{x \in \B} t^{\wt(x)}$, and for normal crystals this agrees with the character of the associated finite-dimensional Lie algebra representation. In particular, a connected normal crystal has a unique highest weight and a unique lowest weight. We will use later the fact that normal crystals are seminormal, so $\varepsilon_i$ and $\varphi_i$ can be used to detect highest and lowest weights.

\subsection{Crystal anti-isomorphisms}\label{subsec:crystal_anti-iso}

Let $\g$ be a complex reductive Lie algebra with the notation above, and let $w_0$ be the longest element of $W$. Recall that there exists a permutation $i \mapsto i'$ of $[r]$ such that
\[
w_0(\alpha_i) = -\alpha_{i'}.
\]

\noindent We say $S: \mathcal{C}_1 \to \mathcal{C}_2$ is a \textit{crystal anti-isomorphism} if $S$ is a bijection and
$$\wt(Sx) = w_0 \wt(x)$$

\begin{center}
\begin{tabular}{ c c }
 $e_{i'}(Sx) = S(f_i(x)), $ & $f_{i'}(Sx) = S(e_i(x)),$  \\ 
 $\varepsilon_{i'}(Sx) = \varphi_{i}(x)$, & $\varphi_{i'}(Sx) = \varepsilon_{i}(x)$, \\    
\end{tabular}
\end{center}
and the same holds for $S^{-1}$. This is a slight generalization of the familiar definition of a crystal involution. A \textit{crystal involution} is a crystal anti-isomorphism $S: \mathcal{C} \to \mathcal{C}$, i.e. the domain and the range are the same set. Crystal involutions are often called Lusztig involutions; in type A they are realized combinatorially by the Sch\"utzenberger involution. Notice a crystal anti-isomorphism sends $e_i$ to $f_{i'}$, $f_i$ to $e_{i'}$, and exchanges highest and lowest weights.

\begin{thm}
    A connected normal crystal has a unique crystal involution. 
\end{thm}

\begin{proof}
    See \cite[Exercises 5.1, 5.2]{BumpSchilling2017} also \cite[Proposition 21.1.2]{Lusztig1993}.
\end{proof}

Crystal anti-isomorphisms are the connection between the branching rules we use here and the rules developed in \cite{Kwon2018, KwonJang2021, LecouveyLenart2020}. To this end, we need to understand how crystal anti-isomorphisms should interact with tensor products. 

 \begin{lemma} \label{lemma:S_base_case}
     Let $\B_1$ and $\B_2$ be normal crystals of $\g$ and $S_1: \B_1 \to \B_1'$, $S_2: \B_2 \to \B_2'$ be crystal anti-isomorphisms. Define the map
    $$S : \B_1 \otimes \B_2 \to \B_2' \otimes \B_1'$$
    \noindent by $S(x_1 \otimes x_2 ) = S_2(x_2) \otimes S_1(x_1)$. Then, $S$ is a crystal anti-isomorphism. 
 \end{lemma}

 \begin{proof}
     We just check that $S$ satisfies the definition of crystal anti-isomorphism. First, $\wt(S(x_1 \otimes x_2)) = w_0(\wt(x_1 \otimes x_2))$:
     \begin{align*}
         \wt(S(x_1 \otimes x_2)) &= \wt(S_2(x_2) \otimes S_1(x_1)) \\
         &= \wt(S_2(x_2)) + \wt(S_1(x_1))\\
         &= w_0 \wt(x_2) + w_0 \wt(x_1)\\
         &= w_0(\wt(x_1 \otimes x_2)). 
     \end{align*}

     \noindent The first line is the definition of $S$. The second line is from the tensor product rule. The third line is from properties of the crystal anti-isomorphisms $S_1, S_2$. The fourth line is from the linearity of $w_0$ and the tensor product rule. Next, $e_{i'}(S(x_1 \otimes x_2)) = S(f_i(x_1 \otimes x_2))$:
\begin{align*}
    e_{i'} ( S(x_1 \otimes x_2)) &= e_{i'}(S_2(x_2) \otimes S_1(x_1))\\
    &= \begin{cases}
        e_{i'}(S_2(x_2)) \otimes S_1(x_1)\quad \text{if } \varepsilon_{i'}(S_2(x_2)) > \varphi_{i'}(S_1(x_1)), \\
        S_2(x_2) \otimes  e_{i'}(S_1(x_1)) \quad \text{if } \varepsilon_{i'}(S_2(x_2)) \leq \varphi_{i'}(S_1(x_1)), \\
    \end{cases}\\
    &= \begin{cases}
        S_2(f_i(x_2)) \otimes S_1(x_1)\quad \text{if } \varphi_{i}(x_2) > \varepsilon_{i}(x_1), \\
        S_2(x_2) \otimes S_1(f_i(x_1)) \quad \text{if } \varphi_{i}(x_2) \leq \varepsilon_{i}(x_1), \\
    \end{cases}\\
    &= \begin{cases}
        S(x_1 \otimes f_i(x_2)) \quad \text{if } \varepsilon_{i}(x_1) < \varphi_{i}(x_2), \\
        S(f_i(x_1) \otimes x_2) \quad \text{if }  \varepsilon_{i}(x_1) \geq \varphi_{i}(x_2), \\
    \end{cases}\\
    &= S(f_i(x_1 \otimes x_2)).
\end{align*}

\noindent The first line is the definition of $S$. The second line is the tensor product rule. The third line comes from properties of the crystal anti-isomorphisms $S_1, S_2$. The fourth is the definition of $S$. The fifth line is the tensor product rule. The remaining verifications are similar and left to the reader.     
\end{proof}

Next, we extend the definition of $S$ to $k$-fold tensor products.

\begin{defn}\label{def:S}
    Let $\B_1, \B_2, \dots , \B_k$ be normal crystals of $\g$ and $S_1, \dots , S_k$ be crystal anti-isomorphisms, where $S_i : \B_i \to \B_i'$. Define the map
    $$S : \B_1 \otimes \B_2 \otimes \cdots \otimes \B_k \to \B_k' \otimes \cdots \otimes \B_2' \otimes \B_1'$$

    \noindent by $S(x_1 \otimes x_2 \otimes \cdots \otimes x_k ) = S_k(x_k) \otimes \cdots \otimes S_2(x_2) \otimes S_1(x_1)$.
\end{defn}

\begin{thm}
    $S$ is a crystal anti-isomorphism.
\end{thm}

\begin{proof}
    Lemma \ref{lemma:S_base_case} is the base case in an induction. Suppose the construction of $S$ yields a crystal anti-isomorphism for a $(k-1)$-fold tensor product and consider $S : \B_1 \otimes \B_2 \otimes \cdots \otimes \B_k \to \B_k' \otimes \cdots \otimes \B_2' \otimes \B_1'$. Let $S' : \B_1 \otimes \B_2 \otimes \cdots \otimes \B_{k-1} \to \B_{k-1}' \otimes \cdots \otimes \B_2' \otimes \B_1'$ and notice 
    $$S((x_1 \otimes \cdots \otimes x_{k-1}) \otimes x_k) = S_k (x_k) \otimes S'(x_1 \otimes \cdots \otimes x_{k-1}).$$
    $S'$ is an anti-isomorphism by induction so $S$ is an anti-isomorphism by Lemma \ref{lemma:S_base_case}.
\end{proof}

\section{\texorpdfstring{$\GL_n$}{GLn}-crystals and tableaux}\label{sect:gln-crystals}

From here on the only crystals encountered will be normal $\gl_n$-crystals. In fact, we will only use $\GL_n$-crystals, i.e. normal $\gl_n$-crystals associated to representations of $\gl_n$ that integrate to representations of $\GL_n$. Recall from Section \ref{sec:conventions for gln} that this places constraints on the action of the center $\mathfrak{z}(\gl_n) = \C I$. 

Normal $\gl_n$-crystals can be characterized axiomatically as in \cite{Stembridge2003}. Since normal crystals correspond to finite-dimensional representations, and the character of a $\GL_n$-crystal agrees with the character of its associated $\GL_n$-representation, we will move freely between the two viewpoints below.

We focus on crystals associated to polynomial representations until the last two sections. The last two sections review a few facts about crystals of rational $\GL_n$-representations. Rational $\GL_n$-crystals are only used for our treatment of graded multiplicities for $(\GL_n^2, \GL_n)$ in Section \ref{sec:graded-mults}. They are not used elsewhere. In Section \ref{subsec:tableaux and knuth equiv}, we review some combinatorial constructions with tableaux. In Section \ref{subsec:Bk and the sig rule} we discuss crystals related to tensor products of the defining representation of $\GL_n$ and the signature rule. In Section \ref{subsec:GLn-crystals-of-tableaux} we recall how to define a $\GL_n$-crystal structure on semistandard tableaux. In Section \ref{subsec:plactic equiv} we discuss plactic equivalence. In Section \ref{subsec:gln anti-iso} we introduce two $\GL_n$-crystal anti-isomorphisms. In Section \ref{subsec:lR rule} we present some results related to the Littlewood-Richardson rule. Finally, in Sections \ref{subsec:rational gln crystals} and \ref{subsec:LR-symmetries} we review the facts we need about rational $\GL_n$-crystals.

\subsection{Tableaux and Knuth equivalence} \label{subsec:tableaux and knuth equiv}
As discussed in Section \ref{section:types_of_tableaux}, $\SST^\lambda_n$ parametrizes a weight basis of the polynomial $\GL_n$ irreducible representation with highest weight $\lambda$, so before discussing $\GL_n$-crystals, we briefly review some of the classical combinatorial theory of tableaux as presented in \cite{Fulton1996}. In particular, we recall how to define a product on the set of semistandard tableaux. The combinatorial approach is through Knuth equivalence which is the equivalence relation on the set of words in the alphabet $\Z_{>0}$ generated by elementary Knuth transformations. A \textit{word} $w$ is just an ordered list of positive integers. We refer to \cite[Chapter 2]{Fulton1996} for the definitions. If $w$ is Knuth equivalent to $w'$, we will write $w \equiv_K w'$.

For a skew tableau $T$ with shape $\lambda/\mu$, we define the \textit{(row) word} of $T$, $w(T)$ to be the entries of $T$  read left to right, bottom to top.

\begin{example}
    Consider $\lambda = (5, 4, 4,2)$ and $\mu = (3,1)$. Let 
    $$T = \begin{ytableau}
        *(black) & *(black) & *(black) & 2 & 5\\
        *(black) & 1 & 3 & 3\\
        2 & 3 & 4 & 4\\
        4 & 4\\
    \end{ytableau}.$$

    Then $w(T) = 44 \cdot 2344 \cdot 133 \cdot 25$. We included the dots to make clear which row the entries are coming from. They are not part of the word.
\end{example}

There are two operations on tableaux we will need: \textit{row bumping} and \textit{sliding}. See \cite[Chapter 1]{Fulton1996} for definitions. These two operations lead to different ways of defining the product of tableaux. The tableau obtained from $T$ by row bumping element $x$ is denoted $T \leftarrow x$. 

Given a skew tableau, we can perform a sequence of slides to produce a semistandard tableau. The tableau so constructed is called the \textit{rectification} of $T$, $\Rect(T)$, and is independent of any choices. The key facts are

\begin{prop} \label{prop:knuth_equiv_of_simple_ops}~
    \begin{enumerate}
        \item $w(T \leftarrow x) \equiv_K w(T) \cdot x$, here the dot is just concatenation of words.
        \item If one skew tableau can be obtained from another by a sequence of slides, their row words are Knuth equivalent. In particular $w(\Rect(T)) \equiv_K w(T).$
    \end{enumerate}
\end{prop}

\begin{proof}
    All facts can be found in \cite{Fulton1996}: (1) is Proposition 2.1, (2) is Proposition 2.2.
\end{proof}

\begin{thm}[{\cite[Theorem 2.1]{Fulton1996}}] \label{thm:knuth_equiv_tableau}
    Every word $x$ is Knuth equivalent to the word of a unique tableau $P(x)$.
\end{thm}

$P(x)$ can be constructed from $x = x_1 x_2 \cdots x_k$ by repeatedly performing row-insertion, $P(x) = (((\varnothing \leftarrow x_1) \leftarrow x_2) \leftarrow \cdots ) \leftarrow x_k$. Next, we recall three ways to define the product of two tableaux $T_1, T_2$. 
\begin{enumerate}
    \item $P(w(T_1) \cdot w(T_2))$ is constructed by concatenating the words of $T_1$ and $T_2$ then forming the unique tableau with word Knuth equivalent to $w(T_1) \cdot w(T_2)$.
    \item $T_1 \leftarrow T_2$ is constructed by row inserting $w(T_2) = x_1 x_2 \dots x_k$ into $T_1$, that is $T_1 \leftarrow T_2 = (((T_1 \leftarrow x_1) \leftarrow x_2) \leftarrow \cdots )\leftarrow x_k$.
    \item Let $T_1 * T_2$ be the skew tableau with an $m \times n$ rectangle of empty squares. Here $m$ is the number of rows of $T_2$ and $n$ is the number of columns of $T_1$. $T_1$ is placed below the rectangle and $T_2$ is placed to the right of it. The third product construction is $\Rect(T_1 * T_2)$.
\end{enumerate}

\begin{thm}[{\cite[Corollary 2.1.2]{Fulton1996}}] \label{thm:tableau_multiplication} 
    The three product constructions above yield the same tableau. 
\end{thm}
Denote this unique tableau $T_1 \cdot T_2$ and call it the \textit{product} of $T_1$ and $T_2$. Hence, the set of all words in the alphabet $\Z_{>0}$ modulo Knuth equivalence becomes a monoid, called the \textit{plactic monoid}. Each element of the monoid is represented by a unique tableau and we can understand the multiplication by any of the three constructions above.

In a few instances it will be useful to extend the definition of product to skew tableaux. We can do this by (1) or (3). The result is a semistandard tableau, in particular $T \cdot \varnothing = \Rect(T)$.

\begin{lemma}\label{lemma:product_of_skew}
    For any skew tableaux $T_1,T_2$, 
    $$w(T_1 \cdot T_2) \equiv_K w(T_1) \cdot w(T_2).$$
\end{lemma}

\begin{proof}
    By definition, $T_1 \cdot T_2 = P(w(T_1) \cdot w(T_2))$ and $w \big (P \big (w(T_1) \cdot w(T_2) \big )\big ) \equiv_K w(T_1) \cdot w(T_2)$ by Theorem \ref{thm:knuth_equiv_tableau}.
\end{proof}
 
This ends what we need about Knuth equivalence. We want to show that the product defined above in terms of Knuth equivalence is compatible with the tensor product of crystals. To do this, we recall how $\SST_n^\lambda$ is given a crystal structure. 

\subsection{\texorpdfstring{$\mathbb{B}_n^{\otimes k}$}{Bnk} and the signature rule} \label{subsec:Bk and the sig rule}
Let $\mathbb{B}_n$ denote the crystal associated to the defining representation of $\GL_n$. 

\begin{example}
The crystal graph of $\mathbb{B}_5$ is,
\tikzstyle{arrow} = [thick, ->, >=stealth]
$$\begin{tikzpicture}[squarednode/.style={rectangle, draw=black!100, minimum size=1.5em}]
    \node[squarednode] (1) {1};
    \node[squarednode] (2) [right=of 1] {2};
    \node[squarednode] (3) [right=of 2] {3};
    \node[squarednode] (4) [right=of 3] {4};
    \node[squarednode] (5) [right=of 4] {5};

    \draw[arrow] (1.east) -- (2.west) node[midway,above] {1};
    \draw[arrow] (2.east) -- (3.west) node[midway,above] {2};
    \draw[arrow] (3.east) -- (4.west) node[midway,above] {3};
    \draw[arrow] (4.east) -- (5.west) node[midway,above] {4};
\end{tikzpicture}.$$

$\wt(\ytableaushort{i}) = \epsilon_i$ and $e_i, f_i, \varepsilon_i, \varphi_i$ can be read from the graph by the seminormality of $\mathbb{B}_5$. For example, $\varphi_i(\ytableaushort{j}) = \delta_{ij}$, the Kronecker delta function. Recall, an arrow labeled $i$ indicates the action of $f_i$.
\end{example}

$\mathbb{B}_n$ is a normal crystal, so it behaves well under tensor product. Hence $\mathbb{B}_n^{\otimes k} $ has the same structure as $(\C^n)^{\otimes k}$ and is well understood by Schur-Weyl duality,

\begin{prop}[Schur-Weyl Duality]
    $$(\C^n)^{\otimes k} \cong \bigoplus_{|\lambda| = k, \ell(\lambda) \leq n} \pi^\lambda_{S_k} \otimes \pi^\lambda_{\GL_n}. $$
\end{prop}

It follows that every irreducible polynomial $\GL_n$-representation of highest weight $\lambda$ occurs in $(\C^n)^{\otimes k}$ with $k=|\lambda|$, and its multiplicity is $\dim \pi^\lambda_{S_k}$. In terms of crystals, this means that $\mathbb{B}_n^{\otimes k}$ contains a connected component isomorphic to the irreducible $\GL_n$-crystal of highest weight $\lambda$, where again $k=|\lambda|$. In the next section, we recall how to define a crystal structure on $\SST_n^\lambda$ by embedding $\SST_n^\lambda$ into $\mathbb{B}_n^{\otimes k}$ and pulling back the crystal structure from its image. Before that, we review a few facts about computing in $\mathbb{B}_n^{\otimes k}$, since computations in connected polynomial $\GL_n$-crystals can be reduced to computations in $\mathbb{B}_n^{\otimes k}$ once such an embedding is chosen.

The weight of an element $x_1 \otimes \cdots \otimes x_k \in \mathbb{B}_n^{\otimes k}$ is
\[
\wt(x_1 \otimes \cdots \otimes x_k)=\wt(x_1)+\cdots+\wt(x_k)
\]
by the tensor product rule. There is a simple combinatorial algorithm for computing the action of $e_i$ and $f_i$ on an element $x_1 \otimes x_2 \otimes \cdots \otimes x_k \in \mathbb{B}_n^{\otimes k}$, known as the \textit{signature rule}; see \cite[Section 2.4]{BumpSchilling2017}. It follows quickly from the tensor product rule and the structure of $\mathbb{B}_n$.

For fixed $i \in [n-1]$, the rule is as follows. Write a $-$ above each occurrence of $i$ and a $+$ above each occurrence of $i+1$. Then successively cancel adjacent $(+,-)$ pairs, ignoring blanks, until no such pair remains. After all cancellations, the remaining signs are of the form
\[
-\cdots-+\cdots+,
\]
possibly with blanks interspersed. Suppose there are $a$ remaining $+$'s and $b$ remaining $-$'s. Then $e_i$ acts on the tensor factor corresponding to the leftmost remaining $+$, and
\[
\varepsilon_i(x_1 \otimes \cdots \otimes x_k)=a.
\]
Similarly, $f_i$ acts on the tensor factor corresponding to the rightmost remaining $-$, and
\[
\varphi_i(x_1 \otimes \cdots \otimes x_k)=b.
\]

\begin{example} \label{ex:sig rule 1}
    Consider the following element of $\mathbb{B}^{\otimes 9}_5$,
    $$x = \begin{ytableau}
        2 & \none[\otimes] & 3 & \none[\otimes] & 5 & \none[\otimes] & 2 & \none[\otimes] & 3 & \none[\otimes] & 1 & \none[\otimes] & 2 & \none[\otimes] & 3 & \none[\otimes] & 3\\
    \end{ytableau}.$$

    The weight of this element is $\wt(x) = \epsilon_1 + 3 \epsilon_2 + 4 \epsilon_3 + \epsilon_5$. Suppose we want to compute $e_2, \varepsilon_2, f_2, \varphi_2$ for $x$. We start by putting $-$ over every $2$ and $+$ over every $3$,
    $$\begin{ytableau}
        \none[-] & \none & \none[+] & \none &  \none & \none & \none[-] & \none &  \none[+] & \none & \none & \none &  \none[-] & \none & \none[+] & \none & \none[+]\\ 
        2 & \none[\otimes] & 3 & \none[\otimes] & 5 & \none[\otimes] & 2 & \none[\otimes] & 3 & \none[\otimes] & 1 & \none[\otimes] & 2 & \none[\otimes] & 3 & \none[\otimes] & 3\\
    \end{ytableau}$$
    Then we cancel all $(+ -)$ pairs,
    $$\begin{ytableau}
        \none[-] & \none & \none & \none &  \none & \none & \none & \none &  \none & \none & \none & \none &  \none & \none & \none[+] & \none & \none[+]\\ 
        2 & \none[\otimes] & 3 & \none[\otimes] & 5 & \none[\otimes] & 2 & \none[\otimes] & 3 & \none[\otimes] & 1 & \none[\otimes] & 2 & \none[\otimes] & 3 & \none[\otimes] & 3\\
    \end{ytableau}$$
    Now $\varepsilon_2(x) =2$, the number of remaining $+$'s and $e_2$ acts on the leftmost remaining $+$, 
    $$e_2(x) = \begin{ytableau}
        2 & \none[\otimes] & 3 & \none[\otimes] & 5 & \none[\otimes] & 2 & \none[\otimes] & 3 & \none[\otimes] & 1 & \none[\otimes] & 2 & \none[\otimes] & *(yellow) 2 & \none[\otimes] & 3\\
    \end{ytableau}.$$
    Similarly, $\varphi_2(x) = 1$, the number of remaining $-$'s and $f_2$ acts on the rightmost remaining $-$,
    $$f_2(x) = \begin{ytableau}
        *(yellow) 3 & \none[\otimes] & 3 & \none[\otimes] & 5 & \none[\otimes] & 2 & \none[\otimes] & 3 & \none[\otimes] & 1 & \none[\otimes] & 2 & \none[\otimes] & 3 & \none[\otimes] & 3\\
    \end{ytableau}.$$
\end{example}

We identify an element $x = x_1 \otimes x_2 \otimes \cdots \otimes x_k$ with the word $x_1 x_2 \cdots x_k$ and write $x$ for both. We extend Knuth equivalence to $\mathbb{B}^{\otimes k}_n$ via this identification. 

A word $x_1x_2\cdots x_k$ is a \textit{Yamanouchi word} if every final segment $x_jx_{j+1}\cdots x_k$ contains at least as many $i$'s as it does $i+1$'s, for all $i\in[n-1]$. A word $x_1x_2\cdots x_k$ is an \textit{anti-Yamanouchi word} if every initial segment $x_1x_2\cdots x_j$ contains at least as many $i+1$'s as it does $i$'s, for all $i\in[n-1]$.

\begin{prop}\label{prop:yamanouchi}~
    \begin{enumerate}
        \item An element of $\mathbb{B}_n^{\otimes k}$ is a highest-weight element if and only if it is a Yamanouchi word.
        \item An element of $\mathbb{B}_n^{\otimes k}$ is a lowest-weight element if and only if it is an anti-Yamanouchi word.
    \end{enumerate}
\end{prop}

\begin{proof}
    (1) is \cite[Proposition 8.2]{BumpSchilling2017}. 
    
    (2) Since $\mathbb{B}_n^{\otimes k}$ is a seminormal crystal, $f_i(x_1 \otimes x_2 \otimes \cdots \otimes x_k)=0$ if and only if $\varphi_i(x_1 \otimes x_2 \otimes \cdots \otimes x_k)=0$. By Lemma \ref{lemma:multiple_tensors},
    \[
    \varphi_i(x_1 \otimes \cdots \otimes x_k) = \max_{j=1}^k \left( \sum_{h=1}^j \varphi_i(x_h) - \sum_{h=1}^{j-1} \varepsilon_i(x_h) \right).
    \]

    Here
    \[
    \varphi_i(x_h) =
    \begin{cases}
        1 & \text{if } x_h = i, \\
        0 & \text{otherwise,}
    \end{cases}
    \qquad
    \varepsilon_i(x_h) =
    \begin{cases}
        1 & \text{if } x_h = i+1,\\
        0 & \text{otherwise.}
    \end{cases}
    \]
    
    From this, it is immediate that $\varphi_i(x_1 \otimes \cdots \otimes x_k)=0$ if and only if every initial segment contains at least as many $i+1$'s as it does $i$'s. Requiring this for all $i\in[n-1]$ is exactly the anti-Yamanouchi word condition.
\end{proof}

\subsection{\texorpdfstring{$\GL_n$}{GLn}-crystals of tableaux} \label{subsec:GLn-crystals-of-tableaux}
Now we recall how to put a $\GL_n$-crystal structure on $\SST_n^\lambda$. More generally, we put a $\GL_n$-crystal structure on $\SST_n^{\lambda/\mu}$ for any skew diagram $\lambda/\mu$. To this end, we embed $\SST_n^{\lambda/\mu}$ into $\mathbb{B}_n^{\otimes k}$, where $k=|\lambda|-|\mu|$. One convenient choice is row reading. Define
\[
\RR : \SST_n^{\lambda/\mu} \to \mathbb{B}_n^{\otimes k}
\]
by setting $\RR(T)=x_1\otimes x_2\otimes \cdots \otimes x_k$, where $x_1x_2\cdots x_k=w(T)$, i.e. we tensor the entries of the row word of $T$ together in order. We call $\RR(T)$ the \textit{row reading} of $T$.

\begin{thm}[{\cite[Theorem 8.8]{BumpSchilling2017}}]
    The subset $\RR(\SST_n^{\lambda/\mu})$ is a subcrystal of $\mathbb{B}_n^{\otimes k}$.
\end{thm}

We give the set $\SST_n^{\lambda/\mu}$ the crystal structure it inherits from $\mathbb{B}_n^{\otimes k}$ via the embedding $\RR$. When $\SST_n^{\lambda/\mu}$ is equipped with this crystal structure, we denote it by $\B_n^{\lambda/\mu}$. The weight of an element $T \in \B_n^{\lambda/\mu}$ is computed by the tensor product rule:
\[
\wt(T)=\sum_{i=1}^n a_i\epsilon_i,
\]
where $a_i$ is the number of $i$'s in the filling of $T$. This is the same definition of weight we used for tableaux in Section \ref{section:types_of_tableaux}. We compute $e_i$, $f_i$, $\varepsilon_i$, and $\varphi_i$ using the signature rule on $\mathbb{B}_n^{\otimes k}$. For example,
\[
\varphi_i(T)=\varphi_i(\RR(T)),
\]
and $\varphi_i(\RR(T))$ is computed by the signature rule.

There are two special cases of $\B_n^{\lambda / \mu}$ that are particularly important for us: $\B_n^{\lambda}$ and $\B_n^{\lambda^\pi}$, where $\lambda$ is a partition and $\lambda^\pi$ is the $180^\circ$ rotation described in Section \ref{section:types_of_tableaux}. Here, and throughout this paper, the rotation $\lambda^\pi$ is taken with respect to the ambient rank $n$, as discussed in Section \ref{section:types_of_tableaux}. As we will see, these provide two equivalent combinatorial models for the irreducible $\GL_n$-representation of highest weight $\lambda$.

\begin{example} \label{ex:sig rule 2}
    Let $\lambda = (4, 2, 1)$ and consider the following element of $\B^\lambda_5$,
    $$T = \begin{ytableau}
        \none[1] & 1 & 3 & 3 & 4\\
        \none[2] & 3 & 4 \\
        \none[3] & 5\\
        \none[4]\\
        \none[5]\\
    \end{ytableau}$$
$\wt(T) = \epsilon_1 + 3\epsilon_3 + 2 \epsilon_4 + \epsilon_5$. It is embedded into $\mathbb{B}^{\otimes 7}_5$ via $\RR$ as 
$$\RR(T) = \begin{ytableau}
    5 & \none[\otimes] & 3 & \none[\otimes] & 4 & \none[\otimes] & 1 & \none[\otimes] & 3 & \none[\otimes] & 3 & \none[\otimes] & 4\\
\end{ytableau}.$$
Suppose we want to compute $e_3, \varepsilon_3, f_3, \varphi_3$. Applying the signature rule,
$$\begin{ytableau}
    \none & \none & \none[-] & \none & \none & \none & \none & \none & \none & \none & \none[-] & \none & \none[+]\\
    5 & \none[\otimes] & 3 & \none[\otimes] & 4 & \none[\otimes] & 1 & \none[\otimes] & 3 & \none[\otimes] & 3 & \none[\otimes] & 4\\
\end{ytableau}$$
and pulling the result back to $T$ we see 
$$\varepsilon_3(T) = 1, e_3(T) = \begin{ytableau}
        \none[1] & 1 & 3 & 3 & *(yellow) 3\\
        \none[2] & 3 & 4 \\
        \none[3] & 5\\
        \none[4]\\
        \none[5]\\
    \end{ytableau}, \varphi_3(T) = 2, \text{ and } f_3(T) = \begin{ytableau}
        \none[1] & 1 & 3 & *(yellow) 4 & 4\\
        \none[2] & 3 & 4 \\
        \none[3] & 5\\
        \none[4]\\
        \none[5]\\
    \end{ytableau}.$$
\end{example}

\begin{example} \label{ex:sig rule 3}
    Continuing with $\lambda = (4,2,1)$, consider the following element of $\B^{\lambda^\pi}_5$,
    $$T = \begin{ytableau}
        \none & \none & \none & \none & \none[1]\\
        \none & \none & \none & \none & \none[2]\\
        \none & \none & \none & 3 & \none[3]\\
        \none & \none & 1 & 4 & \none[4]\\
        3 & 3 & 4 & 5 & \none[5]\\
    \end{ytableau}$$
$\wt(T) = \epsilon_1 + 3\epsilon_3 + 2 \epsilon_4 + \epsilon_5$ and it is embedded into $\mathbb{B}^{\otimes 7}_5$ as
$$\RR(T) = \begin{ytableau}
    3 & \none[\otimes] & 3 & \none[\otimes] & 4 & \none[\otimes] & 5 & \none[\otimes] & 1 & \none[\otimes] & 4 & \none[\otimes] & 3\\
\end{ytableau}.$$
Suppose we want to compute $e_3, \varepsilon_3, f_3, \varphi_3$. Applying the signature rule,
$$\begin{ytableau}
    \none[-] & \none & \none[-] & \none & \none[+] & \none & \none & \none & \none & \none & \none & \none & \none\\
    3 & \none[\otimes] & 3 & \none[\otimes] & 4 & \none[\otimes] & 5 & \none[\otimes] & 1 & \none[\otimes] & 4 & \none[\otimes] & 3\\
\end{ytableau}$$
and pulling the result back to $T$ we see 
$$\varepsilon_3(T) = 1, e_3(T) = \begin{ytableau}
        \none & \none & \none & \none & \none[1]\\
        \none & \none & \none & \none & \none[2]\\
        \none & \none & \none & 3 & \none[3]\\
        \none & \none & 1 & 4 & \none[4]\\
        3 & 3 & *(yellow) 3 & 5 & \none[5]\\
    \end{ytableau}, \varphi_3(T) = 2, \text{ and } f_3(T) = \begin{ytableau}
        \none & \none & \none & \none & \none[1]\\
        \none & \none & \none & \none & \none[2]\\
        \none & \none & \none & 3 & \none[3]\\
        \none & \none & 1 & 4 & \none[4]\\
        3 & *(yellow) 4 & 4 & 5 & \none[5]\\
    \end{ytableau}$$

\end{example}

Next, we establish notation for the highest- and lowest-weight elements in $\B_n^\lambda$ and $\B_n^{\lambda^\pi}$. Let $H_n^\lambda$ denote the tableau on $\lambda$ in which a box $b$ is filled by the number of boxes in the same column as $b$ that are weakly above $b$, i.e.
\[
\#\{\text{boxes in the same column as } b \text{ that are above or equal to } b\}.
\]
Define $H_n^{\lambda^\pi} \in \B_n^{\lambda^\pi}$ in the same way. These are the \textit{highest-weight tableaux}.

Similarly, let $L_n^\lambda$ and $L_n^{\lambda^\pi}$ denote the \textit{lowest-weight tableaux}. A box $b$ in a lowest-weight tableau is filled by
\[
n-\#\{\text{boxes in the same column as } b \text{ that are strictly below } b\}.
\]

\begin{example}
    For $\lambda = (4,2,1)$, we have
    $$H^\lambda_5 = \begin{ytableau}
        \none[1] & 1 & 1 & 1 & 1\\
        \none[2] & 2 & 2 \\
        \none[3] & 3\\
        \none[4]\\
        \none[5]\\
    \end{ytableau}, \ L^\lambda_5 = \begin{ytableau}
        \none[1] & 3 & 4 & 5 & 5\\
        \none[2] & 4 & 5 \\
        \none[3] & 5\\
        \none[4]\\
        \none[5]\\
    \end{ytableau}$$
    and 
    $$H^{\lambda^\pi}_5 = \begin{ytableau}
        \none & \none & \none & \none & \none[1]\\
        \none & \none & \none & \none & \none[2]\\
        \none & \none & \none & 1 & \none[3]\\
        \none & \none & 1 & 2 & \none[4]\\
        1 & 1 & 2 & 3 & \none[5]\\
    \end{ytableau}, \ L^{\lambda^\pi}_5 = \begin{ytableau}
        \none & \none & \none & \none & \none[1]\\
        \none & \none & \none & \none & \none[2]\\
        \none & \none & \none & 3 & \none[3]\\
        \none & \none & 4 & 4 & \none[4]\\
        5 & 5 & 5 & 5 & \none[5]\\
    \end{ytableau}$$
\end{example}

We see that the highest-weight element is easy to write in $\B^\lambda_n$ since $\lambda$ is exactly the highest weight. Similarly, lowest-weight element is easy to write in $\B^{\lambda^\pi}_n$. Next, we justify our definitions in the case of $\B^\lambda_n$. The analogous fact for $\B^{\lambda^\pi}_n$ is in Theorem \ref{thm:crystal_rotation_isomorphisms} below.

\begin{thm}
$\B_n^\lambda$ is a connected crystal with unique highest-weight element $H_n^\lambda$ and lowest-weight element $L_n^\lambda$.
\end{thm}

\begin{proof}
    The connected assertion and highest weight assertion are contained in \cite[Theorem 3.2]{BumpSchilling2017}. The definition of $L^\lambda_n$ ensures $\RR(L^\lambda_n)$ is an anti-Yamanouchi word: if box $b$ in $L^\lambda_n$ is filled by $n -i$ with $i > 0$, then there are $i$ boxes below $b$. In particular, the box directly below $b$ is filled with $n-i+1$. This comes before $b$ in the row reading. It follows that $\# \{ n-i +1 \} \geq \# \{n-i\}$ for all $i \in [n-1]$ and for all initial segments of $\RR(L^\lambda_n)$.
\end{proof}

\subsection{Plactic equivalence} \label{subsec:plactic equiv}
There are two equivalence relations that are important to us: Knuth equivalence $\equiv_K$ discussed above and plactic equivalence denoted $\equiv$.  Plactic equivalence comes from isomorphisms in the monoidal category of normal $\GL_n$-crystals. 

\begin{defn}[{\cite[Definition 8.1]{BumpSchilling2017}}]
    Let $\B_1$ and $\B_2$ be normal $\GL_n$-crystals and $x_i \in \B_i$. Let $\B_i'$ be the connected component of $\B_i$ that contains $x_i$. If $\B_1'$ is isomorphic to $\B_2'$ and if the unique isomorphism $\B_1' \to \B_2'$ takes $x_1$ to $x_2$, then we say $x_1$ and $x_2$ are \textit{plactically equivalent}, and write $x_1 \equiv x_2$.
\end{defn}

\begin{example}\label{example:T=RR(T)}
    For any $T \in \B_n^{\lambda/\mu}$, $T \equiv \RR(T)$ by construction of the crystal structure on $\B^{\lambda/\mu}_n$.
\end{example}

\begin{example}
    The tableau $T$ from Example \ref{ex:sig rule 2} is plactically equivalent to the tableau $T$ from Example \ref{ex:sig rule 3} via the isomorphism of Theorem \ref{thm:crystal_rotation_isomorphisms} (2) below.
\end{example}

\begin{thm} \label{thm:plactic=knuth}
    $x, y \in \mathbb{B}_n^{\otimes k}$ are plactically equivalent if and only if they are Knuth equivalent.
\end{thm}

\begin{proof}
     By \cite[Theorem 8.6(i)]{BumpSchilling2017}, we have \(x\equiv P(x)\) and \(y\equiv P(y)\). Hence
     \[ x\equiv y \iff P(x)\equiv P(y).\]
     We claim that \[P(x)\equiv P(y) \iff P(x)=P(y).\]
     Indeed, say \(P(x)\in \B_n^\mu\) and \(P(y)\in \B_n^\nu\). If \(P(x)\equiv P(y)\), then by definition of plactic equivalence the connected crystals \(\B_n^\mu\) and \(\B_n^\nu\) are isomorphic, so \(\mu=\nu\). Moreover, the only automorphism of the connected normal crystal \(\B_n^\mu\) is the identity, and hence \(P(x)=P(y)\). The converse is immediate.
     
     Finally, by Theorem \ref{thm:knuth_equiv_tableau}, \(x\) is Knuth equivalent to \(P(x)\), and \(y\) is Knuth equivalent to \(P(y)\). Thus 
     \[P(x)=P(y) \quad \Longleftrightarrow \quad x\equiv_K y.\]
     Combining these equivalences gives
 \[x \equiv y \iff P(x) \equiv P(y) \iff P(x) = P(y) \iff x \equiv_K y.\]
\end{proof}

Hence, plactic equivalence and Knuth equivalence both yield the same monoid on the set of all words in the alphabet $[n]$. 

\begin{prop}\label{prop:tensor_plactic}
    Let $T_1 \in \B_n^{\lambda_1/\mu_1}$ and $T_2 \in \B_n^{\lambda_2/\mu_2}$. The following are plactically equivalent.
    \begin{enumerate}
        \item $T_1 \otimes T_2$ in $\B_n^{\lambda_1/\mu_1} \otimes \B_n^{\lambda_2/\mu_2}$,
        \item $\RR(T_1) \otimes \RR(T_2)$ in $\mathbb{B}_n^{\otimes (|\lambda_1/\mu_1| + |\lambda_2 / \mu_2|)}$,
        \item $T_1 \cdot T_2$ in $\B_n^\nu$, where $\nu$ is the shape of $T_1 \cdot T_2$.
    \end{enumerate}
\end{prop}

\begin{proof}
    The equivalence between (1) and (2) follows from Example~\ref{example:T=RR(T)} and the fact that plactic equivalence is compatible with tensor products, see \cite[pg. 113]{BumpSchilling2017}. Since $T_i \equiv \RR(T_i)$ for $i=1,2$, we have
    \[T_1\otimes T_2 \equiv \RR(T_1)\otimes \RR(T_2).\]

    Next we show (2) and (3) are equivalent. By Lemma \ref{lemma:product_of_skew}, $\RR(T_1 \cdot T_2) \equiv_K \RR(T_1) \cdot \RR(T_2)$ so by Theorem \ref{thm:plactic=knuth}, $\RR(T_1 \cdot T_2) \equiv \RR(T_1) \otimes \RR(T_2)$. Since $\RR(T_1 \cdot T_2) \equiv T_1 \cdot T_2$, the result follows.
\end{proof}

\begin{prop} \label{prop:rect_plactic_equiv}
    $\Rect(T) \equiv T$ for any skew tableau $T$.
\end{prop}

\begin{proof}
    $\RR(\Rect(T)) \equiv_K \RR(T)$ by Proposition \ref{prop:knuth_equiv_of_simple_ops} so by Theorem \ref{thm:plactic=knuth}, $\RR(\Rect(T)) \equiv \RR(T)$. Hence by Example \ref{example:T=RR(T)}, $\Rect(T) \equiv T$.
\end{proof}

\subsection{\texorpdfstring{$\GL_n$}{GLn}-crystal anti-isomorphisms} \label{subsec:gln anti-iso}
There are two $\GL_n$-crystal anti-isomorphisms we will use. First, Berenstein and Zelevinsky showed that the unique crystal involution on a connected $\GL_n$-crystal of tableaux is given by \textit{Sch\"utzenberger evacuation}, see \cite[Section 8]{BerensteinZelevinsky1996}. We recall from \cite[6-7]{Lenart2007} that this map can be realized by the following procedure:

\begin{enumerate}
    \item Rotate the tableau $180^\circ$.
    \item Complement the entries via the permutation $i \mapsto w_0(i)= n+1 - i$.
    \item Rectify the skew tableau.
\end{enumerate}

Let $S: \B_n^\lambda \to \B_n^\lambda$ be the map defined by applying the Sch\"utzenberger evacuation. We extend $S$ to tensor products as in Definition \ref{def:S}. This is the first crystal anti-isomorphism that we will use. 

The second comes from just applying the first two steps above. Denote this map $S^\pi$. Notice the map $S^\pi$ makes sense for any crystal of skew tableaux. We define it in this generality, so $S^\pi : \B^{\lambda/\mu}_n \to \B^{\lambda^\pi/\mu^\pi}_n$. Here $\lambda^\pi/ \mu^\pi$ is constructed identically to $\lambda^\pi$ in Section \ref{section:types_of_tableaux}, namely rotating $\lambda/\mu$ by $180^\circ$ about the origin and translating by $n+1$. By construction $S: \B^\lambda_n \to \B^\lambda_n$ is just $\Rect \circ S^\pi$.

To show that $S$ and $S^\pi$ are crystal anti-isomorphisms, we start with the unique involution $I : \mathbb{B}_n \to \mathbb{B}_n$ which takes $i \mapsto w_0(i) = n+1 -i $. We denote $w_0(i)$ by $\widetilde{\iota}$. By the construction in Definition \ref{def:S}, this leads to a crystal involution on $\mathbb{B}^{\otimes k}_n$, namely 
$$\begin{ytableau}
    i_1 & \none[\otimes] & i_2 & \none[\otimes] & \none[\cdots] & \none[\otimes] & i_k & \none[\mapsto] & \widetilde{\iota}_k & \none[\otimes] & \none[\cdots] & \none[\otimes] & \widetilde{\iota}_2 & \none[\otimes] & \widetilde{\iota}_1\\
\end{ytableau}.$$
Denote this involution $I^{\otimes k} : \mathbb{B}_n^{\otimes k} \to \mathbb{B}_n^{\otimes k}$.

\begin{lemma}\label{lemma:S_pi_antiiso}
$S^\pi : \B^{\lambda/\mu}_n \to \B^{\lambda^\pi/\mu^\pi}_n$ is a crystal anti-isomorphism.
\end{lemma}

\begin{proof}
    Since the crystal structures on $\B_n^{\lambda/\mu}$ and $\B_n^{\lambda^\pi/\mu^\pi}$ are inherited from the subcrystals
\[
\RR(\B_n^{\lambda/\mu}) \subset \mathbb{B}_n^{\otimes k}
\qquad\text{and}\qquad
\RR(\B_n^{\lambda^\pi/\mu^\pi}) \subset \mathbb{B}_n^{\otimes k},
\]
it suffices to work on the level of row words. Under row reading, step (1) in the definition of $S^\pi$ reverses the tensor factors, and step (2) replaces each entry $i_j$ by $\widetilde{i_j}=w_0(i_j)$. Hence, on row words, $S^\pi$ is exactly the restriction of $I^{\otimes k}$. Since $I^{\otimes k}$ is a crystal anti-isomorphism and $S^\pi$ is bijective, it follows that $S^\pi$ is a crystal anti-isomorphism.
\end{proof}

We extend $S^\pi$ to tensor products as in Definition \ref{def:S}. 

\begin{thm}\label{thm:crystal_rotation_isomorphisms}~
    \begin{enumerate}
        \item $\B_n^{\lambda^\pi}$ is a connected crystal with unique highest-weight element $H_n^{\lambda^\pi}$ and lowest-weight element $L_n^{\lambda^\pi}$.
        \item $\Rect: \B_n^{\lambda^\pi} \to \B_n^\lambda$ taking $T$ to $\Rect(T)$ is a crystal isomorphism. 
        \item  $\B_n^{\lambda/\mu}$ is isomorphic to $\B_n^{\lambda^\pi/\mu^\pi}$.
        \item $S: \B^\lambda_n \to \B^\lambda_n$ is a crystal involution.
    \end{enumerate}
\end{thm}

\begin{proof}
(1) By Lemma \ref{lemma:S_pi_antiiso}, the map
\(S^\pi:\B_n^\lambda \to \B_n^{\lambda^\pi}\) is a crystal anti-isomorphism.
Since \(\B_n^\lambda\) is connected with unique highest-weight element
\(H_n^\lambda\) and unique lowest-weight element \(L_n^\lambda\), the same is
true of \(\B_n^{\lambda^\pi}\), with highest and lowest weights exchanged by
\(S^\pi\). By the definitions of \(H_n^{\lambda^\pi}\) and
\(L_n^{\lambda^\pi}\), we have
\[
S^\pi(H_n^\lambda)=L_n^{\lambda^\pi},
\qquad
S^\pi(L_n^\lambda)=H_n^{\lambda^\pi}.
\]
Thus \(\B_n^{\lambda^\pi}\) is connected with unique highest-weight element
\(H_n^{\lambda^\pi}\) and unique lowest-weight element \(L_n^{\lambda^\pi}\).

(2) Since $\B^\lambda_n$ and $\B^{\lambda^\pi}_n$ are connected with the same highest weight by (1), there is a unique crystal isomorphism $\iota : \B^{\lambda^\pi}_n \to \B^{\lambda}_n$. We show $\iota = \Rect$. By Proposition \ref{prop:rect_plactic_equiv}, $\Rect(T) \equiv T$. Hence, it suffices to show $\Rect$ is bijective from $\B^{\lambda^\pi}_n \to \B^\lambda_n$. Since $|\B^{\lambda^\pi}_n| = |\B^\lambda_n|$ it suffices to show $\Rect$ is injective. Let $T_1, T_2 \in \B^{\lambda^\pi}_n$ with $T_1 \neq T_2$. Then $T_1 \not\equiv T_2$ since these are distinct elements in a connected normal crystal and the only crystal automorphism is the identity. Hence, by Proposition \ref{prop:rect_plactic_equiv}, $\Rect(T_1) \equiv T_1 \not \equiv T_2 \equiv \Rect(T_2)$, so $\Rect(T_1) \neq \Rect(T_2)$.

(3) $S^\pi$ gives a bijection between highest weights of $\B^{\lambda/\mu}_n$ and lowest weights of $\B^{\lambda^\pi/ \mu^\pi}_n$, the lowest weights of $\B^{\lambda^\pi/\mu^\pi}_n$ being the same as the lowest weights of $\B^{\lambda/\mu}_n$. Since normal crystals are uniquely determined by lowest weights, $\B^{\lambda/\mu}_n \cong \B^{\lambda^\pi/\mu^\pi}_n$.

(4) $S = \Rect \circ S^\pi$. By Lemma \ref{lemma:S_pi_antiiso}, $S^\pi$ is a crystal anti-isomorphism. $\Rect$ is a crystal isomorphism by (2). The composition of a crystal isomorphism  and crystal anti-isomorphism is clearly a crystal anti-isomorphism. In this case, an involution.
\end{proof}

\begin{example}
    For $\lambda = (4,2,1)$, we have
    $H^\lambda_5 = \begin{ytableau}
        \none[1] & 1 & 1 & 1 & 1\\
        \none[2] & 2 & 2 \\
        \none[3] & 3\\
        \none[4]\\
        \none[5]\\
    \end{ytableau}$ applying $S^\pi$ we have $\ L^{\lambda^\pi}_5 = \begin{ytableau}
        \none & \none & \none & \none & \none[1]\\
        \none & \none & \none & \none & \none[2]\\
        \none & \none & \none & 3 & \none[3]\\
        \none & \none & 4 & 4 & \none[4]\\
        5 & 5 & 5 & 5 & \none[5]\\
    \end{ytableau}$. Rectifying, we get $\ L^\lambda_5 = \begin{ytableau}
        \none[1] & 3 & 4 & 5 & 5\\
        \none[2] & 4 & 5 \\
        \none[3] & 5\\
        \none[4]\\
        \none[5]\\
    \end{ytableau}$. Applying $S^\pi$, we have $H^{\lambda^\pi}_5 = \begin{ytableau}
        \none & \none & \none & \none & \none[1]\\
        \none & \none & \none & \none & \none[2]\\
        \none & \none & \none & 1 & \none[3]\\
        \none & \none & 1 & 2 & \none[4]\\
        1 & 1 & 2 & 3 & \none[5]\\
    \end{ytableau}$ and rectifying again we arrive back at $H^\lambda_5$.
\end{example}

\subsection{The Littlewood-Richardson rule} \label{subsec:lR rule}
The key to our approach to the Littlewood-Richardson rule is the following seesaw pair, see \cite[Chapter 9]{BumpSchilling2017},

\begin{center}
    \begin{tikzcd}
    \GL_n^2 \arrow[r,dash] & \GL_n^2   \\
    \GL_n  \arrow[r,dash]\arrow[u, phantom, sloped, "\subset"] & \GL_{2n}, \arrow[u, phantom, sloped, "\supset"]
\end{tikzcd}
\end{center}

\noindent where on the left $\GL_n$ is included diagonally in $\GL_n^2$ and on the right $\GL_n^2$ is included as the Levi subgroup of a parabolic in $\GL_{2n}$. All groups act on the symmetric algebra on $2n \times n $ matrices, $S(M_{2n, n}(\C))$, and each row of the diagram above is a dual pair, see \cite[Appendix 2]{BumpSchilling2017} for details. The usual seesaw argument allows us to conclude,

\begin{prop}[{\cite[Theorem 9.4]{BumpSchilling2017}}]
    $$\mult(\B^\lambda_n, \B^\mu_n\boxtimes \B^\nu_n) = \mult(\B^\mu_n \boxtimes \B^\nu_n, \B^\lambda_{2n}).$$
\end{prop}

We can find the non-Levi branching multiplicity of $\B^\lambda_{n}$ in $\Res^{\GL^2_n}_{\GL_n}(\B^\mu_{n} \boxtimes \B^\nu_{n}) = \B^\mu_n \otimes \B^\nu_n$ coming from the left side by counting highest-weight elements in $\B^\mu_{n} \otimes \B^\nu_{n}$ with weight $\lambda$. Alternatively, we could count lowest-weight elements. In either case, the underlying tool is the tensor product rule. Notice by Proposition \ref{prop-highest/lowest-in-tensor-product} the highest-weight elements in $\B^\mu_{n} \otimes \B^\nu_{n}$ are of the form $T \otimes H_n^\nu$ for some $T \in \B^\mu_n$. Hence we can parametrize highest weights in terms of elements of $\B^\mu_n$. Similarly, lowest-weight elements are of the form $L_n^\mu \otimes T$ for some $T \in \B^\nu_n$ and can be parameterized in terms of elements of $\B^\nu_n$. We make the following definitions:
\begin{align*}
    \HC_{\mu \nu}^\lambda(\GL_n) &\coloneqq \{T \in \B^\mu_n : T \otimes H_n^\nu \equiv H_n^\lambda \},\\
    \LC_{\mu \nu}^\lambda(\GL_n) & \coloneqq \{T \in \B^\nu_n : L_n^\mu \otimes T \equiv L_n^\lambda \}.
\end{align*}

$\HC_{\mu \nu}^\lambda(\GL_n)$ are called \textit{(highest-weight) companion tableaux} and $\LC_{\mu \nu}^\lambda(\GL_n)$ are called \textit{lowest-weight companion tableaux}. By Theorem \ref{thm:crystal_rotation_isomorphisms}, $\B_n^\delta \cong \B_n^{\delta^\pi}$ so we can replace any Young diagram $\delta$ in this definition with $\delta^\pi$. For example $\HC^{\lambda^\pi}_{\mu \nu^\pi}(\GL_n)$ has the same cardinality as  $\HC_{\mu \nu}^\lambda(\GL_n)$.

The Levi branching coming from the right side of the seesaw leads to the standard Littlewood-Richardson tableaux. Again, we can determine these multiplicities by counting highest or lowest-weight elements of $\B_n^\mu \boxtimes \B_n^\nu$ in the $\GL_n \times \GL_n$-crystal obtained by branching $\B_{2n}^\lambda$. 

\begin{thm}[{\cite[Theorem 8.14]{BumpSchilling2017}}]
    Let $\lambda$ be a partition of length $\leq 2n$.
    $$\Res^{\GL_{2n}}_{\GL_n^2}(\B^\lambda_{2n}) = \bigoplus_{\mu \subseteq \lambda} (\B_n^\mu \boxtimes \B_n^{\lambda/\mu}).$$
\end{thm}

By Proposition \ref{prop:exterior_tensor}, to find the multiplicity of $\B_n^\mu \boxtimes \B_n^\nu$ in $\Res^{\GL_{2n}}_{\GL_n^2}(\B^\lambda_{2n})$ we just need to count the highest or lowest weights of $\B_n^\nu$ in $\B_n^{\lambda/\mu}$. We make the following definitions:
\begin{align*}
    \HL_{\mu \nu}^\lambda(\GL_n) &\coloneqq \{T \in \B_n^{\lambda/\mu} : T  \equiv H_n^\nu \},\\
    \LL_{\mu \nu}^\lambda(\GL_n) & \coloneqq \{T \in \B_n^{\lambda/\mu} : T  \equiv L_n^\nu \}.
\end{align*}

$\HL_{\mu \nu}^\lambda(\GL_n)$ are called \textit{(highest-weight) Littlewood-Richardson tableaux} and $\LL_{\mu \nu}^\lambda(\GL_n)$ are called \textit{lowest-weight Littlewood-Richardson tableaux}. Again, there are many variations of this definition. By Theorem \ref{thm:crystal_rotation_isomorphisms}, $\B^{\lambda/\mu}_n \cong \B^{\lambda^\pi/\mu^\pi}_n$ and $\B^\nu_n \cong \B^{\nu^\pi}_n$ so we can apply $\pi$ to $\lambda/\mu$ or $\nu$ or both. An example that will be important later is $\HL^{\lambda^\pi}_{\mu^\pi \nu}(\GL_n)$. From the above discussion, we conclude

\begin{thm} \label{thm:HC=LC=HL=LL}
For $\lambda, \mu,\nu$ all of ambient length $n$,
    $$|\HC_{\mu \nu}^\lambda(\GL_n)| = |\LC_{\mu \nu}^\lambda(\GL_n)| = |\HL_{\mu \nu}^\lambda(\GL_n)| = |\LL_{\mu \nu}^\lambda(\GL_n)|.$$
\end{thm}

$\pi$ can be placed on the $\lambda, \mu, \nu$ in this theorem as described above. The common cardinality of these sets is the \textit{Littlewood-Richardson coefficient} denoted $c^\lambda_{\mu \nu}$ and as we recall below, is independent of $n$.

 Next, we provide characterizations of these highest  and lowest-weight elements in the two cases we will use in applications. Similar statements can be made for the other cases. Recall that $\mu^*$ denotes the highest weight of the dual representation and that $\mu_0$ denotes the projection of $\mu$ onto the semisimple weight lattice of $\sll_n$, see Theorem \ref{thm:GLn_highest_lowest} and the surrounding discussion. In (2) below, $\mu_0, \varphi(T)^* \in P_{++}(\sll_n)$ are elements of the semisimple weight lattice and the partial order $\ge$ is the one defined on $P(\sll_n) = P(\g')$ in Section \ref{sec:reductive Lie algebras}.

\begin{prop} \label{prop:LC_equivalent_sets}
    The following sets are equal.
    \begin{enumerate}
        \item $\LC^\lambda_{\mu \nu}(\GL_n)$,
        \item $\{T \in \B_n^\nu : \mu_0 \geq \varphi(T)^* \text{ and }  \wt(L^\mu_n) + \wt(T) = \wt(L^\lambda_n) \}$,
        \item $\{T \in \B_n^\nu : \RR(L_n^\mu) \otimes \RR(T) \text{ is an anti-Yamanouchi word and } \wt(L^\mu_n) + \wt(T) = \wt(L^\lambda_n) \}$,
        \item $\{T \in \B_n^\nu : L_n^\mu  \cdot T \equiv L_n^\lambda \}$.
    \end{enumerate}
    As above, we can place $\pi$ on any combination of $\lambda, \mu, \nu$ without impacting the cardinality of the set. 
\end{prop}

\begin{proof}
    We show (2), (3), and (4) are equivalent to (1).
    
    (2) By Proposition \ref{prop-highest/lowest-in-tensor-product}, $L^\mu_n \otimes T \equiv L_n^\lambda$ if and only if $-\wt_{ss}(L^\mu_n) \geq \varphi(T)$ and $\wt(L^\mu_n) + \wt(T) = \wt(L^\lambda_n)$. By the definition of \(L_n^\mu\) and Lemma \ref{lemma:hw-dual-rep}, $-\wt_{ss}(L_n^\mu)=-ss(w_0(\mu))=\mu_0^*$. Thus the inequality \(-\wt_{ss}(L_n^\mu)\geq \varphi(T)\) is equivalent to \(\mu_0^*\geq \varphi(T)\). Applying \(-w_0\) to both sides gives \(\mu_0\geq \varphi(T)^*\).
    
    (3) By Proposition \ref{prop:yamanouchi}, $\RR(L^\mu_n) \otimes \RR(T)$ is an anti-Yamanouchi word if and only if $\RR(L^\mu_n) \otimes \RR(T)$ is a lowest-weight element. By Example \ref{example:T=RR(T)}, $\RR(L^\mu_n) \otimes \RR(T) \equiv L^\mu_n \otimes T$. Since $\wt(L^\mu_n) + \wt(T) = \wt(L^\lambda_n)$ it follows that $\RR(L^\mu_n) \otimes \RR(T)$ is an anti-Yamanouchi word if and only if $L^\mu_n \otimes T \equiv L_n^\lambda$. 
    
    (4) $L_n^\lambda \equiv L_n^\mu \otimes T \equiv L_n^\mu \cdot T$. The first equivalence is the definition of $\LC^\lambda_{\mu \nu}(\GL_n)$. The second equivalence is Proposition \ref{prop:tensor_plactic}.
\end{proof}

Now for $\HL^{\lambda}_{\mu \nu}(\GL_n)$,

\begin{prop} \label{prop:HL equivalent sets}
    The following sets are equal.
    \begin{enumerate}
        \item $\HL^{\lambda}_{\mu \nu}(\GL_n)$,
        \item $\{T \in \B_n^{\lambda / \mu} : \varepsilon(T)=0 \text{ and } \wt(T) = \nu\}$,
        \item $\{T \in \B_n^{\lambda / \mu} : \RR(T) \text{ is a Yamanouchi word and } \wt(T) = \nu \}$,
        \item $\{T \in \B_n^{\lambda / \mu} : \Rect (T) = H_n^\nu \}$.
    \end{enumerate}
    As above, we can place $\pi$ on any combination of $\lambda/\mu$ and $\nu$ without impacting the cardinality of the set.
\end{prop}

\begin{proof}
    We show (2), (3), and (4) are equivalent to (1).
    
    (2) follows immediately from the definition of highest-weight element and the fact that $\B_n^{\lambda / \mu}$ is seminormal. 
    
    (3) By Proposition \ref{prop:yamanouchi}, $\RR(T)$ is a Yamanouchi word if and only if $\RR(T)$ is a highest-weight element. Since $\wt(T) = \nu$, this happens if and only if  $\RR(T) \equiv H^\nu_n$. By Example \ref{example:T=RR(T)}, $T \equiv \RR(T) \equiv H^\nu_n$.
    
    (4) By Proposition \ref{prop:knuth_equiv_of_simple_ops}, rectification preserves Knuth equivalence so $\Rect(T) = H^\nu_n$ if and only if $\RR(T) \equiv_K \RR(H^\nu_n)$. By Theorem \ref{thm:plactic=knuth}, $\RR(T) \equiv \RR(H^\nu_n)$ and by Example \ref{example:T=RR(T)}, $T \equiv H^\nu_n$.
\end{proof}

\begin{cor}
    For $n,m \geq \ell(\lambda)$, $|\HL^\lambda_{\mu \nu}(\GL_m)| = |\HL^\lambda_{\mu \nu}(\GL_n)| = c^\lambda_{\mu \nu}$, where $c^\lambda_{\mu \nu}$ is the standard Littlewood-Richardson coefficient. 
\end{cor}

\begin{proof}
    A skew tableau $T$ on shape $\lambda/\mu$ such that $\RR(T)$ is a Yamanouchi word can only involve numbers up to $\ell(\lambda)$, so the set in Proposition \ref{prop:HL equivalent sets} (3) is the same for any $n,m \geq \ell(\lambda)$.
\end{proof}

Next, we define bijections between these sets in a few instances that will be important later. Similar arguments can provide many more bijections. 

\begin{defn}
    Let $\phi$ be the map that takes $T \in \LC^\lambda_{\mu \nu}(\GL_n)$ to the skew tableau $\phi(T)$ of shape $\lambda^\pi/\mu^\pi$ defined as follows. The $i$th row of $\phi(T)$ is filled by the row numbers of the entries equal to $i$ in $T$, recorded with multiplicity and placed in weakly increasing order. Equivalently, an entry $i$ in the $j$th row of $T$ corresponds to an entry $j$ in the $i$th row of $\phi(T)$.
\end{defn}

See \cite{vanLeeuwen2001}. We show in Theorem \ref{thm:phi bijection LC to HL} that $\phi$ is well-defined.

\begin{thm} \label{thm:phi bijection LC to HL}
    $\phi$ is a bijection from  $\LC^\lambda_{\mu \nu}(\GL_n)$ to $\HL^{\lambda^\pi}_{\mu^\pi \nu}(\GL_n)$. $\phi^{-1}$ is given by the same rule as $\phi$.   
\end{thm}

\begin{proof}
    $T \in \LC^\lambda_{\mu \nu}(\GL_n)$, so $T$ has shape $\nu$ and $\wt(T) = \wt(L^\lambda_n) - \wt(L^\mu_n) = \lambda^\pi/\mu^\pi$. Hence $\phi$ is a well-defined map to $\lambda^\pi/\mu^\pi$ with filling $\nu$. 
    
    Next, notice that the strictly increasing columns condition on $T$ ensures that $\RR(\phi(T))$ is a Yamanouchi word. Assume $i>1$ and let $\{ x_{i1}, \dots , x_{il} \}$ be the filling of row $i$ in $T$. Then we record an $i$ in rows $x_{i1}, \dots , x_{il}$ of $\phi(T)$. Since the columns of $T$ are strictly increasing, there exists an $x_{i-1j} <  x_{ij} $ in row $i-1$ above $x_{ij}$ for $j = 1, \dots l$. Thus, for any $i$ in row $x_{ij}$ of $\phi(T)$, we record an $i-1$ in row $x_{i-1j}$. Since $x_{i-1j} < x_{ij}$, there are at least as many $i-1$ as there are $i$ in any final segment of $\RR(\phi(T))$. Hence $\RR(\phi(T))$ is a Yamanouchi word. 

    To show $\phi(T)$ is semistandard, it suffices to show the columns are strictly increasing. By our characterization of $\LC^\lambda_{\mu \nu}(\GL_n)$ in Proposition \ref{prop:LC_equivalent_sets}, we see $\RR(L^\mu_n) \otimes \RR(T)$ is an anti-Yamanouchi word. Write this word as $w =m_1 \cdots m_j x_1 x_2 \cdots x_k$, where $m_1 \otimes \cdots \otimes m_j = \RR(L^\mu_n)$ and $x_1 \otimes \cdots \otimes x_k = \RR(T)$. The fact that  $\RR(L^\mu_n) \otimes \RR(T)$ is anti-Yamanouchi implies that every initial segment $m_1 \cdots m_j x_1 \cdots x_i$ containing the first $j+i$ terms  of $w$ is plactically equivalent to a lowest-weight tableau $L^{\lambda_i^\pi}_n$. So we get an increasing chain of lowest-weight tableaux 
    $$L^{\mu^\pi}_n = L_n^{\lambda^\pi_0} \subset L^{\lambda^\pi_1}_n \subset L^{\lambda^\pi_2}_n \subset \cdots \subset L^{\lambda^\pi_k}_n = L^{\lambda^\pi}_n,$$
    \noindent where $\lambda_{i+1}^\pi/\lambda_{i}^\pi$ contains one box $b$ in row $x_{i+1}$. $\phi(T)$ is constructed by putting the $T$ row number of $x_{i+1}$, denoted $\text{row}(x_{i+1})$ in $b$ for $x_1, \dots , x_k$. There is no box above $b$, otherwise $\lambda^\pi_i$ would not be a rotated Young diagram and if $x_{i+1} < n$, there is a box $c$ below $b$. If $c \not\in \mu^\pi$, let $r$ be the filling of $c$. $r$ is the row number of some $x_l \in \{x_1, \dots , x_i \}$. Since $c$ is one row lower than $b$ in $\lambda_{i+1}^\pi$, $x_l= x_{i+1}+1$ and since $c \in \lambda_i^\pi$ and $b \not \in \lambda_i^\pi$, $l < i+1$. Because the rows of $T$ are weakly increasing, the fact that $l < i+1$ and $x_l > x_{i+1}$, implies that $x_l$ must be from a different row of $T$ than $x_{i+1}$. Since rows of $T$ are read bottom to top, $\text{row}(x_l) > \text{row}(x_{i+1})$. Recalling that the filling of $b$ is $\text{row}(x_{i+1})$ and the filling of $c$ is $\text{row}(x_l)$, we see the columns of $\phi(T)$ are strictly increasing. So $\phi$ maps to $\HL^{\lambda^\pi}_{\mu^\pi \nu}(\GL_n)$.

    It is easy to see $\phi$ is injective. Since $|\LC^\lambda_{\mu \nu}(\GL_n) | = | \HL^{\lambda^\pi}_{\mu^\pi \nu}(\GL_n) |$ by Theorem \ref{thm:HC=LC=HL=LL}, we conclude $\phi$ is a bijection. It is then clear that $\phi^{-1}$ is given by the same rule. 
\end{proof}

\begin{example}
    Let $\lambda = (4, 4,2,1)$, $\mu = (3,1)$, and $\nu = (4,2,1)$. The tableau \begin{ytableau}
        2 & 3 & 4& 4 \\
        3 & 5 \\
        4\\
    \end{ytableau} in $\LC^\lambda_{\mu \nu}(\GL_5)$ corresponds under $\phi$ to \begin{ytableau}
        \none & \none & \none & \none & \none[1] \\
        \none & \none & \none & 1 & \none[2] \\
        \none & \none & 1 & 2 & \none[3] \\
        1 & 1 & 3 & *(black) & \none[4] \\
        2 & *(black) & *(black) & *(black) & \none[5] \\
    \end{ytableau}.

\end{example}

\begin{thm} \label{thm:S-bijection}
    $S:\B^\nu_n \to \B^\nu_n$ restricts to a bijection between $\HC_{\nu \mu}^\lambda(\GL_n)$ and $\LC_{\mu \nu}^\lambda(\GL_n)$ and $S^\pi: \B^{\nu^\pi}_n \to \B^\nu_n$ restricts to a bijection between $\HC_{\nu^\pi \mu}^\lambda(\GL_n)$ and $\LC_{\mu \nu}^\lambda(\GL_n)$.
\end{thm}

\begin{proof}
    Let $T \in \HC^\lambda_{\nu \mu}(\GL_n)$ so $T \otimes H_n^\mu \equiv H_n^\lambda$. Applying $S$, we see $S(T \otimes H_n^\mu) \equiv S(H_n^\lambda) = L_n^\lambda$ by Theorem \ref{thm:crystal_rotation_isomorphisms} and the definition of crystal anti-isomorphism. Using the extension of $S$ to tensor products, Definition \ref{def:S}, $S(T \otimes H_n^\mu) = S(H_n^\mu) \otimes S(T) = L_n^\mu \otimes S(T)$. Hence $L_n^\lambda \equiv L_n^\mu \otimes S(T) $ and $S(T) \in \LC_{\mu \nu}^\lambda(\GL_n)$. Similarly, if $T \in \LC_{\mu \nu}^\lambda(\GL_n)$, then $L_n^\mu \otimes T \equiv L_n^\lambda$. Applying $S$, we find $S(T) \otimes H_n^\mu \equiv H_n^\lambda$, so $S(T) \in \HC^\lambda_{\nu \mu}(\GL_n)$.

    The argument for $\HC^\lambda_{\nu^\pi \mu}(\GL_n)$ and $\LC^\lambda_{\mu \nu}(\GL_n)$ is the same, using $S^\pi$ in place of $S$.
\end{proof}

\subsection{Rational \texorpdfstring{$\GL_n$}{GLn}-crystals} \label{subsec:rational gln crystals}
Let $\mathbb{B}_n^*$ denote the dual crystal of the defining representation of $\GL_n$. 

\begin{example}
The crystal graph of $\mathbb{B}_5^*$ is,

$$\tikzstyle{arrow} = [thick, <-, >=stealth]
\begin{tikzpicture}[squarednode/.style={rectangle, draw=black!100, minimum size=1.5em}]
    \node[squarednode] (1) {$\overline{1}$};
    \node[squarednode] (2) [right=of 1] {$\overline{2}$};
    \node[squarednode] (3) [right=of 2] {$\overline{3}$};
    \node[squarednode] (4) [right=of 3] {$\overline{4}$};
    \node[squarednode] (5) [right=of 4] {$\overline{5}$};

    \draw[arrow] (1.east) -- (2.west) node[midway,above] {1};
    \draw[arrow] (2.east) -- (3.west) node[midway,above] {2};
    \draw[arrow] (3.east) -- (4.west) node[midway,above] {3};
    \draw[arrow] (4.east) -- (5.west) node[midway,above] {4};
\end{tikzpicture}.$$

$\wt(\begin{ytableau}
    \overline{i}\\
\end{ytableau}) = -\epsilon_i$ and $e_i, f_i, \varepsilon_i, \varphi_i$ can be read from the graph.
\end{example}

Let $\lambda = (\lambda^+, \lambda^-) \in \widehat{\GL}_n$. Similar to the case of polynomial $\GL_n$-crystals, we put a crystal structure on the rational $\GL_n$-tableaux $\mathcal{T}^\lambda_{\GL_n}$ of Definition \ref{def:K-tableaux} by embedding into $(\mathbb{B}_n^*)^{\otimes k^-} \otimes \mathbb{B}_n^{\otimes k^+}$ where $k^- = |\lambda^-|$ and $k^+ = |\lambda^+|$. Let $p=\ell(\lambda^+)$ and $q=\ell(\lambda^-)$. As discussed in
\cite{STEMBRIDGE1987}, we can think of the ordered pair
$\lambda=(\lambda^+,\lambda^-)$ as the staircase
\[
(\lambda_1^+,\ldots,\lambda_p^+,0,\ldots,0,
-\lambda_q^-,\ldots,-\lambda_1^-)\in \Z^n,
\]
which corresponds to the highest weight
\(
\lambda_1^+\epsilon_1+\cdots+\lambda_p^+\epsilon_p
-\lambda_q^-\epsilon_{n-q+1}-\cdots-\lambda_1^-\epsilon_n.
\) Similarly, we can think of a rational tableau as a staircase tableau by rotating the tableau on $\lambda^-$ by $180^\circ$ and negating the entries. We embed the staircase tableaux into  $(\mathbb{B}_n^*)^{\otimes k^-} \otimes \mathbb{B}_n^{\otimes k^+}$ by row reading. 

\begin{example}
     Let $\lambda = (\lambda^+, \lambda^-) = ((3,2), (3,1))$ corresponding to the irreducible rational representation of $\GL_5$ with highest weight $3 \epsilon_1 + 2\epsilon_2 - \epsilon_4 - 3 \epsilon_5$. Here is an example of the correspondence between rational tableaux and staircase tableaux.
    $$\Bigg ( \begin{ytableau}
        1 & 4 & 5\\
        2 & 5\\
    \end{ytableau}, \begin{ytableau}
        3 & 4 & 5\\
        5 \\
    \end{ytableau} \Bigg ) \longleftrightarrow \begin{ytableau}
        \none & \none & \none & \none[1] & 1 & 4 & 5\\
        \none & \none & \none & \none[2] & 2 & 5\\
        \none & \none & \none & \none[3]\\ 
        \none & \none & \overline{5} & \none[4]\\
        \overline{5} & \overline{4} & \overline{3} & \none[5]\\
    \end{ytableau} = T.$$

    The weight of this tableau is $\wt(T) = \epsilon_1 + \epsilon_2 - \epsilon_3$. It is embedded into $(\mathbb{B}_5^*)^{\otimes 4} \otimes \mathbb{B}_5^{\otimes 5}$ by row reading as
    $$\RR(T) = \begin{ytableau}
        \overline{5} & \none[\otimes] & \overline{4} & \none[\otimes] & \overline{3} & \none[\otimes] & \overline{5} & \none[\otimes] & 2 & \none[\otimes] & 5 & \none[\otimes] & 1 & \none[\otimes] & 4 & \none[\otimes] & 5\\
    \end{ytableau}.$$
\end{example}

As in the case of polynomial $\GL_n$-crystals, $\RR(\mathcal{T}^\lambda_{\GL_n})$ is a subcrystal of $(\mathbb{B}_n^*)^{\otimes k^-} \otimes \mathbb{B}_n^{\otimes k^+}$. We give $\mathcal{T}^\lambda_{\GL_n}$ the crystal structure it inherits via this embedding and denote the resulting crystal $\B^\lambda_n$, see \cite[Section 3.4]{Kwon2009}.

We extend the signature rule to $(\mathbb{B}_n^*)^{\otimes k^-} \otimes \mathbb{B}_n^{\otimes k^+}$ as follows. To compute the action of $e_i$ and $f_i$, write a $-$ above every unbarred $i$ and a $+$ above every unbarred $i+1$, as in the $\mathbb{B}_n$ case. In addition, write a $-$ above every $\overline{i+1}$ and a $+$ above every $\overline{i}$. Then proceed as before: successively cancel adjacent $(+,-)$ pairs, ignoring blanks, until no such pair remains. Suppose there are $a$ remaining $+$'s and $b$ remaining $-$'s. Then $e_i$ acts on the tensor factor corresponding to the leftmost remaining $+$, and $\varepsilon_i(x_1 \otimes \cdots \otimes x_k)=a.$
Similarly, $f_i$ acts on the tensor factor corresponding to the rightmost remaining $-$, and
\(
\varphi_i(x_1 \otimes \cdots \otimes x_k)=b.
\)

\begin{example} \label{ex:sig rule 4}
    Suppose we want to compute $e_4, \varepsilon_4, f_4, \varphi_4$ for the example above. We start by putting $-$ over every $4$ and every $\overline{5}$ and $+$ over every $5$ and every $\overline{4}$.
    $$\begin{ytableau}
         \none[-] & \none & \none[+] & \none & \none & \none & \none[-] & \none & \none & \none & \none[+] & \none & \none & \none & \none[-] & \none & \none[+]\\
         \overline{5} & \none[\otimes] & \overline{4} & \none[\otimes] & \overline{3} & \none[\otimes] & \overline{5} & \none[\otimes] & 2 & \none[\otimes] & 5 & \none[\otimes] & 1 & \none[\otimes] & 4 & \none[\otimes] & 5\\
    \end{ytableau}$$
    Then cancel $(+ -)$ pairs,
    $$\begin{ytableau}
         \none[-] & \none & \none & \none & \none & \none & \none & \none & \none & \none & \none & \none & \none & \none & \none & \none & \none[+]\\
         \overline{5} & \none[\otimes] & \overline{4} & \none[\otimes] & \overline{3} & \none[\otimes] & \overline{5} & \none[\otimes] & 2 & \none[\otimes] & 5 & \none[\otimes] & 1 & \none[\otimes] & 4 & \none[\otimes] & 5\\
    \end{ytableau}$$

    From this we see $\varepsilon_4(\RR(T)) = \varphi_4(\RR(T)) = 1$, $e_4$ changes the $5$ below the $+$ to a $4$ and $f_4$ changes the $\overline{5}$ below the $-$ to a $\overline{4}$. This pulls back to the staircase tableau. For example,
    $$e_4(T) = \begin{ytableau}
        \none & \none & \none & \none[1] & 1 & 4 & *(yellow) 4\\
        \none & \none & \none & \none[2] & 2 & 5\\
        \none & \none & \none & \none[3]\\ 
        \none & \none & \overline{5} & \none[4]\\
        \overline{5} & \overline{4} & \overline{3} & \none[5]\\
    \end{ytableau}$$
\end{example}

We extend the definition of $H_n^\lambda$ and $L_n^\lambda$ to $\lambda=(\lambda^+,\lambda^-)$. Continue to view $\lambda$ as a staircase. Let $H_n^\lambda$ denote the tableau on $\lambda$ defined as follows. A box $b^+ \in \lambda^+$ is filled by 
\[
\#\{\text{boxes weakly above } b^+ \text{ in its column}\},
\]
while a box $b^- \in \lambda^-$ is filled by
\[
\overline{\,n-\#\{\text{boxes strictly above } b^- \text{ in its column}\}\,}.
\]
Similarly, let $L_n^\lambda$ denote the tableau on $\lambda$ defined as follows. A box $b^+ \in \lambda^+$ is filled by
\[
n-\#\{\text{boxes strictly below } b^+ \text{ in its column}\},
\]
while a box $b^- \in \lambda^-$ is filled by 
\[
\overline{\,\#\{\text{boxes weakly below } b^- \text{ in its column}\}\,}.
\]

\begin{example}
    Let $\lambda = (\lambda^+, \lambda^-)$ with $\lambda^+ = (3,2)$ and $\lambda^-=(3,1)$. The associated staircase is $(3,2,0,-1,-3)$ and the corresponding highest weight is $3 \epsilon_1 + 2 \epsilon_2 - \epsilon_4 - 3\epsilon_5$. The longest Weyl group element $w_0$ takes the highest weight $\lambda$ to the lowest weight $w_0(\lambda) = 3 \epsilon_5 + 2 \epsilon_4 - \epsilon_2 - 3\epsilon_1$. We can see this reflected in $H^\lambda_5$ and $L^\lambda_5$ below.

    $$H^\lambda_5 = \begin{ytableau}
        \none & \none & \none & \none[1] & 1 & 1 & 1\\
        \none & \none & \none & \none[2] & 2 & 2\\
        \none & \none & \none & \none[3]\\ 
        \none & \none & \overline{5} & \none[4]\\
        \overline{5} & \overline{5} & \overline{4} & \none[5]\\
    \end{ytableau} \text{ and } L^\lambda_5 = \begin{ytableau}
        \none & \none & \none & \none[1] & 4 & 4 & 5\\
        \none & \none & \none & \none[2] & 5 & 5\\
        \none & \none & \none & \none[3]\\ 
        \none & \none & \overline{2} & \none[4]\\
        \overline{1} & \overline{1} & \overline{1} & \none[5]\\
    \end{ytableau}$$
\end{example}

\begin{thm}
For $\lambda = (\lambda^+, \lambda^-)\in \widehat{\GL}_n$, $\B_n^\lambda$ is a connected normal crystal with unique highest-weight element $H_n^\lambda$ and lowest-weight element $L_n^\lambda$.
\end{thm}

\begin{proof}
    First observe $(\mathbb{B}_n^*)^{\otimes k^-} \otimes \mathbb{B}_n^{\otimes k^+}$ is a normal crystal since $\mathbb{B}_n$ and $\mathbb{B}_n^*$ are normal, and normality is preserved under tensor product, \cite[Theorem 5.20]{BumpSchilling2017}. Hence, $\B^\lambda_n$ is a normal $\GL_n$-crystal since we constructed it as a subcrystal of $(\mathbb{B}_n^*)^{\otimes k^-} \otimes \mathbb{B}_n^{\otimes k^+}$. Further, by \cite[Proposition 2.4]{STEMBRIDGE1987}, the character of $\B^\lambda_n$ equals the character of the irreducible $\GL_n$-representation with highest weight $\lambda$. Since the character theory of normal crystals and finite-dimensional representations of $\GL_n$ is the same, $\B^\lambda_n$ must be the unique connected normal crystal with highest weight $\lambda = \wt(H^\lambda_n)$ and lowest weight $w_0(\lambda) = \wt(L^\lambda_n)$. 
\end{proof}

The definition of plactic equivalence was given for normal $\GL_n$-crystals, so applies to rational $\GL_n$-crystals. Let $\lambda, \mu, \nu \in \widehat{\GL}_n$, i.e. rational $\GL_n$-representations are allowed, and set
$$\LC^{\lambda}_{\mu \nu}(\GL_n) \coloneqq \{T \in \B^{\nu}_n : L^{\mu}_n \otimes T \equiv L^{\lambda}_n \}.$$
\noindent We extend the notation for Littlewood-Richardson coefficients to rational $\GL_n$-representations, setting $c^\lambda_{\mu \nu} \coloneqq |\LC^\lambda_{\mu \nu}(\GL_n)|$.

\begin{prop} \label{prop:LC rational equiv set}
    $$\LC^{\lambda}_{\mu \nu}(\GL_n) = \{ T \in \B^{\nu}_n : \mu_0 \geq \varphi(T)^* \text{ and } \wt(L^{\mu}_n) + \wt(T) = \wt(L^{\lambda}_n) \}.$$
\end{prop}

\begin{proof}
    This follows identically to the proof of Proposition \ref{prop:LC_equivalent_sets} (2).
\end{proof}

\begin{prop}\label{prop:GLn branching rule}
    For $[\lambda_1, \lambda_2] \in \widehat{\GL}_n^2$ with $\lambda_i = (\lambda_i^+, \lambda_i^-)$, and $\nu = (\nu^+, \nu^-) \in \widehat{\GL}_n$, 
    $$b^{[\lambda_1, \lambda_2]}_\nu = |\LC^\nu_{\lambda_1 \lambda_2}(\GL_n)| = c^\nu_{\lambda_1 \lambda_2}.$$
\end{prop}

\begin{proof}
    $\Res^{\GL_n^2}_{\GL_n} \pi^{[\lambda_1, \lambda_2]}_{\GL_n^2} = \pi^{\lambda_1}_{\GL_n} \otimes \pi^{\lambda_2}_{\GL_n}$ and the result follows from the correspondence between normal $\GL_n$-crystals and finite-dimensional $\GL_n$-representations.
\end{proof}

\subsection{Some symmetries}\label{subsec:LR-symmetries}
We continue to work with rational $\GL_n$-representations in this section. By the correspondence between normal $\GL_n$-crystals and finite-dimensional $\GL_n$-representations,
$$ c^\lambda_{\mu \nu} = \dim \Hom_{\GL_n}(\pi^\lambda_{\GL_n}, \pi^\mu_{\GL_n} \otimes \pi^\nu_{\GL_n}).$$

\noindent Notice $\dim \Hom_{\GL_n}(\pi^\lambda_{\GL_n}, \pi^\mu_{\GL_n} \otimes \pi^\nu_{\GL_n})$ is equal to the dimension of the $\GL_n$ invariants in a triple product,
$$\dim \bigg ((\pi^\lambda_{\GL_n})^* \otimes \pi^\mu_{\GL_n} \otimes \pi^\nu_{\GL_n} \bigg)^{\GL_n}.$$

This reveals many symmetries in the coefficients themselves, see \cite{BerensteinZelevinsky1992}. In particular, $S_3$ acts by permuting the three irreducible representations, yielding:
$$c^\lambda_{\mu \nu} = c^\lambda_{\nu \mu} = c^{\mu^*}_{\lambda^* \nu} = c^{\nu^*}_{\mu \lambda^*} = c^{\mu^*}_{\nu \lambda^*} = c^{\nu^*}_{\lambda^* \mu}.$$

Taking the dual of the triple product preserves the multiplicity of invariants. Hence, the following coefficients are also equal to $c^\lambda_{\mu \nu}$:
$$c^{\lambda^*}_{\mu^* \nu^*} = c^{\lambda^*}_{\nu^* \mu^*} = c^{\mu}_{\lambda \nu^*} = c^{\nu}_{\mu^* \lambda} = c^{\mu}_{\nu^* \lambda} = c^{\nu}_{\lambda \mu^*}.$$

\noindent We will make use of these symmetries to establish the graded multiplicity formula for $(\GL_n^2 , \GL_n)$.

\section{Branching rules}\label{sec:branching rules}
In this section we develop generalizations of the Littlewood restriction rules: combinatorial rules for computing 
\[
\mult(\pi^\nu_{\Or_n}, \pi^\lambda_{\GL_n})
\qquad \text{and} \qquad
\mult(\pi^\nu_{\Sp_{2n}}, \pi^\lambda_{\GL_{2n}}),
\]
where the $\GL_n$- and $\GL_{2n}$-representations are polynomial. We focus on lowest-weight companion tableaux, as defined in Section \ref{subsec:lR rule}, since this makes the connection with the familiar $K$-tableaux from Definition \ref{def:K-tableaux} most natural. These results can also be translated into various equivalent statements via the maps $S$, $S^\pi$, and $\phi$ from Section \ref{subsec:lR rule}. 

In Section \ref{subsec:branch-examples} we state the rules, along with examples and corollaries. The examples illustrate that these rules are computationally no more difficult than the classical Littlewood restriction rules. We also show that these rules stabilize to the Littlewood restriction rules. In our formulation, the stable range becomes visually transparent: the skew shape $\lambda^\pi/\mu^\pi$ has enough empty rows above it that the $\Or_n$ and $\Sp_{2n}$ highlighting conditions described below are automatically satisfied by any Littlewood-Richardson tableau. Proofs, which rely heavily on \cite{Kwon2018, KwonJang2021, LecouveyLenart2020}, are given in Section \ref{subsec:proof-symp} and Section \ref{subsec:proof-orth}.

\subsection{Rules, examples, and Littlewood restriction} \label{subsec:branch-examples}
To state the rules, we define two variants of the Littlewood-Richardson coefficient. 

\begin{defn}~
    \begin{enumerate}
        \item Let $\lambda,\mu,\nu \in \Par_{2n}$ and define $\LC^\lambda_{\mu \nu}(\Sp_{2n}) \coloneqq \LC^\lambda_{\mu \nu}(\GL_{2n}) \cap \mathcal{T}_{\Sp_{2n}}^{\nu, H}$. Call elements of this set \textit{symplectic companion tableaux}. Define $c^\lambda_{\mu \nu}(\Sp_{2n}) \coloneqq |\LC^\lambda_{\mu \nu}(\Sp_{2n})|$ and call this a \textit{symplectic Littlewood-Richardson coefficient}.
        \item Let $\lambda,\mu,\nu \in \Par_{n}$ and define $\LC^\lambda_{\mu \nu}(\Or_n) \coloneqq \LC^\lambda_{\mu \nu}(\GL_n) \cap \mathcal{T}_{\Or_n}^\nu$. Call elements of this set \textit{orthogonal companion tableaux}. Define $c^\lambda_{\mu \nu}(\Or_{n}) \coloneqq |\LC^\lambda_{\mu \nu}(\Or_n)|$ and call this an \textit{orthogonal Littlewood-Richardson coefficient}.
    \end{enumerate}
\end{defn}
$\LC^\lambda_{\mu \nu}(\Sp_{2n})$ is just the set of $\Sp_{2n}$-tableaux $T$ on $\nu$ such that $L_{2n}^{\mu} \otimes T \equiv L_{2n}^\lambda$. Similarly, $\LC^\lambda_{\mu \nu}(\Or_n)$ is the set of $\Or_{n}$-tableaux $T$ on $\nu$ such that $L_n^{\mu} \otimes T \equiv L_n^\lambda$.

We can also characterize these modified Littlewood-Richardson coefficients in terms of Littlewood-Richardson tableaux. For a semistandard skew tableau $T$, we say a box $b$ in $T$ is \textit{$\Or_n$-highlighted} if the entry $i$ in $b$ is the first or second occurrence of $i$ when $w(T)$ is read right to left. Similarly, a box is \textit{$\Sp_{2n}$-highlighted} if the  entry $i$ in $b$ is the first occurrence of $i$ when $w(T)$ is read right to left.

\begin{defn}~
    \begin{enumerate}
        \item $\HL^{\lambda^\pi}_{\mu^\pi \nu}(\Sp_{2n}) \coloneqq \{T \in \HL^{\lambda^\pi}_{\mu^\pi \nu}(\GL_{2n}) : \# \{\text{$\Sp_{2n}$-highlighted boxes in the top $2i$ rows} \} \leq i 
        \text{ for all } i \in [n]\}$ and call the elements of this set \textit{symplectic Littlewood-Richardson tableaux}. 
        \item $\HL^{\lambda^\pi}_{\mu^\pi \nu}(\Or_n) \coloneqq \{T \in \HL^{\lambda^\pi}_{\mu^\pi \nu}(\GL_n) : \# \{\text{$\Or_n$-highlighted boxes in the top $i$ rows} \} \leq i 
        \text{ for all } i \in [n]\}$ and call the elements of this set \textit{orthogonal Littlewood-Richardson tableaux}.  
    \end{enumerate}
\end{defn}

The following lemma says that $c^\lambda_{\mu \nu}(\Sp_{2n})$ and $c^\lambda_{\mu \nu}(\Or_{n})$ can be computed via companion tableaux or LR tableaux, just like $c^\lambda_{\mu \nu}$.
\begin{lemma}
    $$|\LC_{\mu \nu}^\lambda(\Sp_{2n})| = |\HL^{\lambda^\pi}_{\mu^\pi \nu}(\Sp_{2n})|,$$
    $$|\LC_{\mu \nu}^\lambda(\Or_{n})| = |\HL^{\lambda^\pi}_{\mu^\pi \nu}(\Or_{n})|.$$
\end{lemma}

\begin{proof}
By Theorem \ref{thm:phi bijection LC to HL}, the map
\[\phi:\LC^\lambda_{\mu\nu}(\GL_n)\to \HL^{\lambda^\pi}_{\mu^\pi\nu}(\GL_n)\]
is a bijection. We show that this bijection restricts to a bijection between $\LC_{\mu \nu}^\lambda(\Or_{n})$ and $\HL^{\lambda^\pi}_{\mu^\pi \nu}(\Or_{n})$. Let $T\in \LC^\lambda_{\mu\nu}(\GL_n)$. Recall that $T$ is an $\Or_n$-tableau if and only if 
\[\#\{ \text{boxes in the first two columns of $T$ with entry $\leq i$}\} \leq i\]
for all $i \in [n]$. We compare this condition with the highlighting condition on $\phi(T)$.

Fix a row $j$ of $T$ and write its entries from left to right as
\[x_{j1}\leq x_{j2}\leq \cdots \leq x_{jm}.\]
Under the map $\phi$, the entry $x_{jk}$ contributes an entry $j$ in row $x_{jk}$ of $\phi(T)$. Since the row of $T$ is weakly increasing, the entries coming from the leftmost boxes of row $j$ are placed in weakly higher rows of $\phi(T)$. Therefore, when $w(\phi(T))$ is read right to left, the earliest occurrences of $j$ come from the leftmost boxes of row $j$ of $T$. If several of the $x_{jk}$ are equal, the corresponding entries $j$ lie in the same row of $\phi(T)$ and are indistinguishable for the present counting. In particular, the $\Or_n$-highlighted boxes of $\phi(T)$ filled with $j$ correspond precisely to the first two boxes in row $j$ of $T$, when they exist.

Moreover, entries of $T$ with values $\leq i$ are mapped to rows $\leq i$ in $\phi(T)$. It follows that, for each $i$, the number of $\Or_n$-highlighted boxes of $\phi(T)$ in the top $i$ rows is equal to the number of boxes in the first two columns of $T$ with entry $\leq i$. Thus the $\Or_n$-tableau condition on $T$ is equivalent to the defining condition for $\phi(T)$ to lie in $\HL_{\mu^\pi\nu}^{\lambda^\pi}(\Or_n)$. Therefore $\phi$ restricts to a bijection
\[\LC_{\mu\nu}^\lambda(\Or_n)\to \HL_{\mu^\pi\nu}^{\lambda^\pi}(\Or_n).\]

Similarly, recall that $T$ is an $\Sp_{2n}$-tableau if and only if 
$$\#\{ \text{boxes in first column of $T$ with entry $\leq 2i$}\} \leq i$$
for all $i \in [n]$. The $\Sp_{2n}$-highlighted boxes of $\phi(T)$ correspond to the boxes coming from the first column of $T$. The entries of $T$ with values $\leq 2i$ are mapped to rows $\leq 2i$ in $\phi(T)$. So the number of boxes in the first column of $T$ with entry $\leq 2i$ is equal to the number of highlighted boxes in the first $2i$ rows of $\phi(T)$. The lemma follows.
\end{proof}

Next, we give a few examples illustrating the computation of $c^\lambda_{\mu \nu}(\Or_n)$ and $c^\lambda_{\mu \nu}(\Sp_{2n})$. We see through the examples that for $n$ sufficiently large, $c^\lambda_{\mu \nu}( \Or_n) = c^\lambda_{\mu \nu}(\Sp_{2n}) = c^\lambda_{\mu \nu}$.

\begin{example}
Consider $(\GL_4,\Or_4)$. Suppose we would like to compute $c^\lambda_{\mu \nu}(\Or_4)$ for $\lambda = (7, 3, 3, 3), \mu = (2, 2, 2)$, and $\nu = (5, 3, 1, 1)$. We start by determining  $\HL^{\lambda^\pi}_{\mu^\pi \nu}(\GL_4)$, the set of LR tableaux on $\lambda^\pi / \mu^\pi$ with filling $\nu$. There is only one such tableau, displayed on the left below. Its companion tableau is displayed next to it on the right. The companion tableau is in $\LC^\lambda_{\mu \nu}(\GL_4)$.
\[ \begin{ytableau}
    \none & \none & \none & \none & 1 & *(green) 1 & *(green) 1 & \none[1]\\
    \none & \none & \none & \none & *(green) 2 & *(black) & *(black) & \none[2]\\
    \none & \none &\none & \none & *(green) 3 & *(black) & *(black) & \none[3]\\
    1 & 1 & 2 & *(green) 2 & *(green) 4 & *(black) & *(black) & \none[4]\\
\end{ytableau}
\begin{ytableau}
            *(green) 1 & *(green) 1 & 1 & 4 & 4\\
            *(green) 2 & *(green) 4 & 4\\
            *(green) 3\\
            *(green) 4\\
        \end{ytableau}
        \]

Next, we determine which members of $\HL^{\lambda^\pi}_{\mu^\pi \nu}(\GL_4)$ are in $\HL^{\lambda^\pi}_{\mu^\pi \nu}(\Or_4)$. Equivalently, we could determine which members of $\LC^\lambda_{\mu \nu}(\GL_4)$ are in $\LC^\lambda_{\mu \nu}(\Or_4)$. On the LR tableau side, this amounts to checking the condition on highlighted boxes. The $\Or_4$-highlighted boxes are colored in the example above. There are two highlighted boxes in the first row, violating the condition for $i=1$. Thus, this LR tableau is excluded. The highlighted boxes in the LR tableau correspond to the boxes in the first two columns of the companion tableau.  On the companion tableau side, we see the companion tableau is not an $\Or_4$-tableau, so it is excluded. 

To get a sense of how these modified LR coefficients lead to a result that stabilizes to the classical Littlewood restriction rule, consider $(\GL_5, \Or_5)$. We present the relevant tableaux below. We keep the same partitions and fillings, now viewed with ambient length $5$ by appending trailing zeros.
\[ \begin{ytableau}
    \none & \none & \none & \none & \none & \none & \none & \none[1]\\
    \none & \none & \none & \none & 1 & *(green) 1 & *(green) 1 & \none[2]\\
    \none & \none & \none & \none & *(green) 2 & *(black) & *(black) & \none[3]\\
    \none & \none &\none & \none & *(green) 3 & *(black) & *(black) & \none[4]\\
    1 & 1 & 2 & *(green) 2 & *(green) 4 & *(black) & *(black) & \none[5]\\
\end{ytableau}
\begin{ytableau}
            *(green) 2 & *(green) 2 & 2 & 5 & 5\\
            *(green) 3 & *(green) 5 & 5\\
            *(green) 4\\
            *(green) 5\\
        \end{ytableau}
        \]
The tableau in $\HL^{\lambda^\pi}_{\mu^\pi \nu}(\GL_5)$ is the same as the $n=4$ case, except the bottom row is now 5. On the companion side, this has the effect of adding one to every entry: recall, the entries in row $i$ of the companion tableau come from the rows containing entry $i$ in the LR tableau. For example, the  1's in the LR tableau occur (written with multiplicity) in rows $\{2, 2, 2, 5, 5 \}$ when $n=5$ versus rows $\{1, 1, 1, 4, 4 \}$ when $n=4$. Hence, the entries in row 1 of the companion tableau are $22255$ when $n=5$ versus $11144$ when $n=4$.

This tableau still does not contribute to $c^\lambda_{\mu \nu}(\Or_5)$. On the LR side, the problem does not occur until row 5: we see six highlighted boxes in the first five rows. Similarly, on the companion side, we are okay until $i=5$. The shape supported on $[5]$ is not an $\Or_5$-tableau. For $i<5$, however, the shape supported on $[i]$ is an $\Or_i$-tableau.

Next consider $(\GL_6, \Or_6)$. The tableaux are below. 
\[ \begin{ytableau}
    \none & \none & \none & \none & \none & \none & \none & \none[1]\\
    \none & \none & \none & \none & \none & \none & \none & \none[2]\\
    \none & \none & \none & \none & 1 & *(green) 1 & *(green) 1 & \none[3]\\
    \none & \none & \none & \none & *(green) 2 & *(black) & *(black) & \none[4]\\
    \none & \none &\none & \none & *(green) 3 & *(black) & *(black) & \none[5]\\
    1 & 1 & 2 & *(green) 2 & *(green) 4 & *(black) & *(black) & \none[6]\\
\end{ytableau}
\begin{ytableau}
            *(green) 3 & *(green) 3 & 3 & 6 & 6\\
            *(green) 4 & *(green) 6 & 6\\
            *(green) 5\\
            *(green) 6\\
        \end{ytableau}
        \]

Finally, we reach the stable situation where this tableau is counted. On the left, the number of highlighted boxes in the first $i$ rows never exceeds $i$. On the right, the companion tableau is an $\Or_6$-tableau. For $n > 6$, it is easy to see this tableau continues to be counted when computing $c^\lambda_{\mu \nu}(\Or_n)$. From the above, we see
\[
c^\lambda_{\mu \nu}(\Or_n) =
\begin{cases}
    0 & \text{if } n < 6,\\
    1 & \text{if } n \geq 6.
\end{cases}
\]
Thus $c^\lambda_{\mu \nu}(\Or_n)=c^\lambda_{\mu \nu}$ when $n\geq 6$.
\end{example}

\begin{example}
    Similarly, start with $(\GL_4,\Sp_4)$ and consider
$\lambda = (5,4,3,2)$, $\mu = (2,2,1,1)$, and $\nu = (4,3,1)$.
One LR tableau in $\HL^{\lambda^\pi}_{\mu^\pi \nu}(\GL_4)$ is displayed on the left below, and its corresponding companion tableau in $\LC^\lambda_{\mu \nu}(\GL_4)$ is displayed on the right.
    \[ \begin{ytableau}
    \none & \none & \none & *(green) 1 & *(black) & \none[1]\\
    \none & \none & 1 & *(green) 2 & *(black) & \none[2]\\
    \none & 1 & 2 & *(black) & *(black) & \none[3]\\
    1 & 2 & *(green) 3 & *(black) & *(black) & \none[4]\\
\end{ytableau}
\begin{ytableau}
            *(green) 1 & 2 & 3 & 4\\
            *(green) 2 & 3 & 4\\
            *(green) 4\\
        \end{ytableau}
        \]

    The $\Sp_4$-highlighted boxes are colored. On the left, there are two highlighted boxes in the first two rows, violating the condition on highlighted boxes. On the right, the companion tableau is not an $\Sp_4$-tableau. So this tableau is not counted when computing $c^\lambda_{\mu \nu}(\Sp_4)$.

    Next, consider $(\GL_6,\Sp_6)$, viewing the same partitions with ambient length $6$ by appending two trailing zeros. The relevant tableaux are below.

\[ \begin{ytableau}
    \none & \none & \none & \none  & \none  & \none[1]\\
    \none & \none & \none & \none  & \none  & \none[2]\\
    \none & \none & \none & *(green) 1 & *(black) & \none[3]\\
    \none & \none & 1 & *(green) 2 & *(black) & \none[4]\\
    \none & 1 & 2 & *(black) & *(black) & \none[5]\\
    1 & 2 & *(green) 3 & *(black) & *(black) & \none[6]\\
\end{ytableau}
\begin{ytableau}
            *(green) 3 & 4 & 5 & 6\\
            *(green) 4 & 5 & 6\\
            *(green) 6\\
        \end{ytableau}
        \]
Notice we add 2 to all entries in the companion tableau (coming from the 2 additional empty rows added above the $n=4$ LR tableau). The condition on $\Sp_{6}$-highlighted boxes is satisfied, so this tableau is counted in $c^\lambda_{\mu \nu}(\Sp_6)$ and it is easy to see it will be counted for all $\Sp_{2n}$ with $n \geq 3$.
\end{example}

For us, the importance of these modified Littlewood-Richardson coefficients is the following result.  

\begin{thm}\label{thm:main_branching}~
\begin{enumerate}
    \item For $\lambda \in \Par_{2n}$, $\nu \in \widehat{\Sp}_{2n}$,
    $$\mult(\pi^\nu_{\Sp_{2n}}, \pi^\lambda_{\GL_{2n}}) = \sum_{\mu \in \Par^{(1,1)}_{2n}} c^\lambda_{\mu \nu }(\Sp_{2n}).$$
    \item For $\lambda \in \Par_n$, $\nu \in \widehat{\Or}_n$,
    $$\mult(\pi^\nu_{\Or_{n}}, \pi^\lambda_{\GL_{n}}) =  \sum_{\mu \in \Par_n^{(2)}} c^\lambda_{\mu \nu }(\Or_n).$$
\end{enumerate}

\end{thm}

Notice that, although we follow the convention of displaying this as a sum over an infinite set, $\Par_n^{(2)}$ or $\Par_{2n}^{(1,1)}$, only finitely many terms actually contribute. Indeed, $c^\lambda_{\mu\nu}(K)=0$ unless $\mu$ is contained in $\lambda$ and $|\lambda|-|\mu|=|\nu|$. Next, we give a couple of examples using Theorem \ref{thm:main_branching}.

\begin{example}
We compute the decomposition of the $\GL_n$ irreducible representation $\lambda = (2,2)$ when restricted to $\Or_n$. This decomposition depends on $n$, but as we show, stabilizes for $n \geq 4$. Notice, for all but three $\mu \in \Par^{(2)}_n$, $c^\lambda_{\mu \nu}(\Or_n)$ is automatically zero regardless of $\nu$ because $\mu \nsubseteq \lambda$. We only need to check $c^\lambda_{\mu \nu}(\Or_n)$ for $\mu = (0), (2), (2,2)$. The easiest way to approach this is to find all LR tableaux on $\lambda^\pi/\mu^\pi$ for $\mu = (0), (2), (2,2)$ and then see which satisfy the condition on $\Or_n$-highlighted boxes. 

For $n =2$, the LR tableaux are:
\[\begin{ytableau}
    *(green) 1 & *(green) 1 & \none[1]\\
    *(green) 2 & *(green) 2 & \none[2]\\
\end{ytableau}, 
\begin{ytableau}
    *(green) 1 & *(green) 1 & \none[1]\\
    *(black) & *(black) & \none[2]\\
\end{ytableau}, 
\begin{ytableau}
    *(black) & *(black) & \none[1]\\
    *(black) & *(black) & \none[2]\\
\end{ytableau}.\]

The $\mu^\pi$ are filled in black. The corresponding $\nu$ is just the weight of the skew tableau. In this example $\nu = (2,2), (2), (0)$ in order from left to right. Observe that the first two skew tableaux violate the condition on highlighted boxes: there are two highlighted boxes in the first row. The final skew tableau vacuously satisfies the condition on highlighted boxes. Hence $\pi^{(2,2)}_{\GL_2}|_{\Or_2} = \pi_{\Or_2}^{(0)}$ which is expected as $\pi^{(2,2)}_{\GL_2} = \det^2$ and the square of the determinant is trivial on $\Or_2$.

For $n= 3$, the LR tableaux are:
\[\begin{ytableau}
    \none & \none & \none[1]\\
    *(green) 1 & *(green) 1 & \none[2]\\
    *(green) 2 & *(green) 2 & \none[3]\\
\end{ytableau}, 
\begin{ytableau}
    \none & \none & \none[1]\\
    *(green) 1 & *(green) 1 & \none[2]\\
    *(black) & *(black) & \none[3]\\
\end{ytableau}, 
\begin{ytableau}
    \none & \none & \none[1]\\
    *(black) & *(black) & \none[2]\\
    *(black) & *(black) & \none[3]\\
\end{ytableau}.\]

The second tableau is allowed, but the first still is not since there are 4 highlighted boxes in first 3 rows. Hence, $\pi^{(2,2)}_{\GL_3}|_{\Or_3} = \pi^{(2)}_{\Or_3} \oplus \pi_{\Or_3}^{(0)}$.

For $n= 4$, the LR tableaux are:
\[\begin{ytableau}
    \none & \none & \none[1]\\
    \none & \none & \none[2]\\
    *(green) 1 & *(green) 1 & \none[3]\\
    *(green) 2 & *(green) 2 & \none[4]\\
\end{ytableau}, 
\begin{ytableau}
    \none & \none & \none[1]\\
    \none & \none & \none[2]\\
    *(green) 1 & *(green) 1 & \none[3]\\
    *(black) & *(black) & \none[4]\\
\end{ytableau}, 
\begin{ytableau}
    \none & \none & \none[1]\\
    \none & \none & \none[2]\\
    *(black) & *(black) & \none[3]\\
    *(black) & *(black) & \none[4]\\
\end{ytableau}.\]

All tableaux are allowed. Hence, $\pi^{(2,2)}_{\GL_4}|_{\Or_4} = \pi^{(2,2)}_{\Or_4} \oplus \pi^{(2)}_{\Or_4} \oplus \pi_{\Or_4}^{(0)}$. It is also clear that at $n=4$ the decomposition stabilizes, that is for all $n \geq 4$ we have $\pi^{(2,2)}_{\GL_n}|_{\Or_n} = \pi^{(2,2)}_{\Or_n} \oplus \pi^{(2)}_{\Or_n} \oplus \pi_{\Or_n}^{(0)}$. As discussed in \cite[Proposition 10.3.6 and 10.3.7]{GoodmanWallach2009}, this is the decomposition of the space of curvature tensors on $\C^n$ relative to the natural action of $\Or_n$ into Weyl conformal curvature tensors, traceless Ricci curvature tensors, and scalar curvature tensors, respectively.
\end{example}

\begin{example}
    Next, we compute the decomposition of the $\GL_{2n}$ irreducible representation $\lambda = (2,2)$ when restricted to $\Sp_{2n}$. This is very similar to the orthogonal example so we move faster. 

    For $(\GL_2, \Sp_2)$, the relevant Littlewood-Richardson tableaux are:

    \[\begin{ytableau}
    1 &  *(green) 1 & \none[1]\\
    2 &  *(green) 2 & \none[2]\\
\end{ytableau}, 
\begin{ytableau}
    *(green) 1 & *(black) & \none[1]\\
    *(green) 2 & *(black) & \none[2]\\
\end{ytableau}, 
\begin{ytableau}
    *(black) & *(black) & \none[1]\\
    *(black) & *(black) & \none[2]\\
\end{ytableau}.\]

\noindent Only the LR tableau corresponding to $\mu = (2,2)$ is allowed. Hence $\pi^{(2,2)}_{\GL_2}|_{\Sp_2} = \pi^{(0)}_{\Sp_2}$ which is expected since $\pi^{(2,2)}_{\GL_2} = \det^2$ and the square of the determinant is trivial on $\Sp_2$.

For $(\GL_4, \Sp_4)$, we have: 

    \[\begin{ytableau}
    \none & \none & \none[1]\\
    \none & \none & \none[2]\\
    1 &  *(green) 1 & \none[3]\\
    2 &  *(green) 2 & \none[4]\\
\end{ytableau}, 
\begin{ytableau}
    \none & \none & \none[1]\\
    \none & \none & \none[2]\\
    *(green) 1 & *(black) & \none[3]\\
    *(green) 2 & *(black) & \none[4]\\
\end{ytableau}, 
\begin{ytableau}
    \none & \none & \none[1]\\
    \none & \none & \none[2]\\
    *(black) & *(black) & \none[3]\\
    *(black) & *(black) & \none[4]\\
\end{ytableau},\]

\noindent all of which are allowed, so $\pi^{(2,2)}_{\GL_4}|_{\Sp_4} = \pi^{(2,2)}_{\Sp_4} \oplus \pi^{(1,1)}_{\Sp_4} \oplus \pi_{\Sp_4}^{(0)}$. The decomposition stabilizes at $n = 2$, i.e. this decomposition holds for all $(\GL_{2n}, \Sp_{2n})$ with $n \geq 2$.
\end{example}

Next we formalize the stability condition which we observed in the examples above. The following corollary is the classical Littlewood restriction rule, see \cite{HoweTanWillenbring2005, Littlewood1944}. It is included to show how the branching rules above specialize to the classical rules in the stable range.

\begin{cor}[Littlewood Restriction Rules]~
    \begin{enumerate}
        \item For $\ell(\lambda) \leq n$,
        $$\mult(\pi^\nu_{\Sp_{2n}}, \pi^\lambda_{\GL_{2n}}) = \sum_{\mu \in \Par_{2n}^{(1,1)}} c^\lambda_{\mu \nu }.$$
        \item For $\ell(\lambda) \leq n/2$,
        $$\mult(\pi^\nu_{\Or_{n}}, \pi^\lambda_{\GL_{n}}) = \sum_{\mu \in \Par_n^{(2)}} c^\lambda_{\mu \nu }.$$
    \end{enumerate}
\end{cor}

\begin{proof}
    (1) For $\Sp_{2n}$, suppose $\ell(\lambda) \leq n$. Then any LR tableau 
    $T \in \HL^{\lambda^\pi}_{\mu^\pi, \nu}(\GL_{2n})$ has at least $n$ empty rows above it. Notice that there can be at most one $\Sp_{2n}$-highlighted box in each nonempty row thanks to the Yamanouchi condition: the first row in which $a$ occurs cannot contain $a+1$. Thus there are no highlighted boxes in rows $1$ through $n$, and at most one highlighted box in each row below row $n$. Hence,
    \[
    \#\{\text{$\Sp_{2n}$-highlighted boxes in the top } 2i\text{ rows} \} \leq \max \{2i - n, 0 \}.
    \]
    The $\Sp_{2n}$ condition is satisfied as long as $\max \{2i - n, 0 \} \leq i$ for all $i \in [n]$, which is automatic. Hence no elements of $\HL^{\lambda^\pi}_{\mu^\pi \nu}(\GL_{2n})$ are excluded, so
    \[
    \HL^{\lambda^\pi}_{\mu^\pi \nu}(\GL_{2n}) = \HL^{\lambda^\pi}_{\mu^\pi \nu}(\Sp_{2n}).
    \]
    The result follows from Theorem \ref{thm:main_branching}.
    
    (2) The $\Or_n$ situation is similar. Suppose $\ell(\lambda) \leq n/2$. Then any LR tableau 
    $T \in \HL^{\lambda^\pi}_{\mu^\pi, \nu}(\GL_n)$ has at least $n/2$ empty rows above it. We first observe that there can be at most two $\Or_n$-highlighted boxes in each nonempty row. Indeed, the $\Or_n$-highlighted boxes are the first and second occurrences of each entry when $w(T)$ is read right to left. By the Yamanouchi condition and the weakly increasing row condition, the first occurrence of $a+1$ must occur strictly below the first occurrence of $a$. Thus the first occurrences contribute at most one highlighted box to each row. The same argument applies to the second occurrences. Hence each nonempty row contains at most two $\Or_n$-highlighted boxes.
    
    Since the top $i$ rows contain at most $\max\{i-n/2,0\}$ nonempty rows and each such row contains $\leq 2$ highlighted boxes, we have
    \[
    \#\{\text{$\Or_n$-highlighted boxes in the top } i \text{ rows} \}
    \leq 2\max \{i-n/2,0 \}
    =
    \max \{2i- n, 0 \}.
    \]
    The $\Or_{n}$ condition is satisfied as long as $\max \{2i- n, 0 \} \leq i$ for all $i \in [n]$, which is automatic. Hence
    \[
    \HL^{\lambda^\pi}_{\mu^\pi \nu}(\GL_n) = \HL^{\lambda^\pi}_{\mu^\pi \nu}(\Or_{n}),
    \]
    and the result follows from Theorem \ref{thm:main_branching}.
\end{proof}

\subsection{Proof of the symplectic group rule} \label{subsec:proof-symp}
We prove
    $$\mult(\pi^\nu_{\Sp_{2n}}, \pi^\lambda_{\GL_{2n}}) = \sum_{\mu \in \Par^{(1,1)}_{2n}} |\LC_{\mu \nu}^\lambda(\Sp_{2n})|.$$

\begin{proof}
    In \cite{LecouveyLenart2020}, a subset $\overline{\HC}^{\lambda'}_{\nu' \mu'}(\GL_{2n})$ of $\HC^{\lambda'}_{\nu' \mu'}(\GL_{2n})$ is defined and $\overline{c}^\lambda_{\nu \mu} \coloneqq |\overline{\HC}^{\lambda'}_{\nu' \mu'}(\GL_{2n})|$, (in \cite{LecouveyLenart2020} these  sets are denoted $\overline{{\LR}}^{\lambda'}_{\nu', \mu'}$ and $\LR^{\lambda'}_{\nu', \mu'}$). \cite[Corollary 6.9]{LecouveyLenart2020} gives us

\begin{equation}\label{eq:Sp2nLL}
    \mult(\pi^\nu_{\Sp_{2n}}, \pi^\lambda_{\GL_{2n}}) = \sum_{\mu \in \Par^{(1,1)}_{2n}} \bar{c}^\lambda_{\nu \mu}.
\end{equation}

Now consider $b \in \HC^\lambda_{\nu \mu}(\GL_{2n})$ with $\mu \in \Par^{(1,1)}_{2n}$. Let $\sigma_\mu$ be the bijection $\HC^\lambda_{\nu \mu}(\GL_{2n}) \to \HC^{\lambda'}_{\nu' \mu'}(\GL_{2n})$ as in \cite{LecouveyLenart2020} and $S: \B^\nu_{2n} \to \B^\nu_{2n}$ the Sch\"utzenberger involution. By \cite[Lemma 6.11]{LecouveyLenart2020}, $\sigma_\mu(b) \in \overline{{\HC}}^{\lambda'}_{\nu' \mu'}(\GL_{2n})$ if and only if $S(b)$ has entries in row $i$ at least $2i-1$ for $i = 1, \dots , n$. As $S$ gives a bijection between $\HC^\lambda_{\nu \mu}(\GL_{2n})$ and $\LC^\lambda_{\mu \nu}(\GL_{2n})$ (Theorem \ref{thm:S-bijection} above), we see the set of $\LC^\lambda_{\mu \nu}(\GL_{2n})$ with entries in row $i$ at least $2i-1$ for $i = 1, \dots , n$ has cardinality $\bar{c}_{\nu, \mu}^\lambda$. By the discussion after Definition \ref{def:K-tableaux}, $\mathcal{T}^{\nu,H}_{\Sp_{2n}}$ can be characterized as the set of tableaux on $\nu$ with entries in row $i$ at least $2i-1$ for $i = 1, \dots , n$. Since $\LC_{\mu \nu}^\lambda(\Sp_{2n}) = \LC^\lambda_{\mu \nu}(\GL_{2n}) \cap \mathcal{T}^{\nu, H}_{\Sp_{2n}}$, $|\LC_{\mu \nu}^\lambda(\Sp_{2n})| = \bar{c}^\lambda_{\nu \mu}$ and the theorem follows from (\ref{eq:Sp2nLL}).

\end{proof}

\subsection{Proof of the orthogonal group rule} \label{subsec:proof-orth}

We prove
    $$\mult(\pi^\nu_{\Or_{n}}, \pi^\lambda_{\GL_{n}}) = \sum_{\mu \in \Par^{(2)}_n} |\LC_{\mu \nu}^\lambda(\Or_n)|.$$

Let $\nu \in \widehat{\Or}_n$, and write $\nu=(\nu_1,\dots,\nu_n)$ and $\nu'=(\nu_1',\dots,\nu_l')$ for its conjugate partition. In particular, $\nu_1'+\nu_2'\leq n$. We start with the following definition of Jang and Kwon,

\begin{defn}[{\cite[Definition 4.6]{KwonJang2021}}] \label{def_kwon}
     Let $T \in \B^{\nu^\pi}_n$ and if $n - 2 \nu_1' < 0$, set $r = n - \nu_1'$. Otherwise, $r=\nu_1'$. For $1 \leq i \leq \nu_1'$ and $1 \leq j \leq \nu_2'$, let 

    \begin{itemize}
        \item $\sigma_i = $ the $i$th entry in the rightmost column of $T$ from the bottom,
        \item $\tau_j = $ the $j$th entry in the second rightmost column of $T$ from the bottom,
        
        \item \begin{equation*}
          m_i =
        \begin{cases}
      \min\{n+1 - \sigma_i, 2i-1\} & \text{if $1 \leq i \leq r$,}\\
      \min\{n+1 - \sigma_i, n +i - \nu_1'\} & \text{if $r < i \leq \nu'_1$}
        \end{cases}       
        \end{equation*}
        
        \item $n_j = $ the $j$th smallest number in $\{j+1, \dots , n\} \setminus \{m_{j+1}, \dots , m_{\nu_1'} \}$.
    \end{itemize}

    \noindent Say $T$ is a \textit{Jang-Kwon orthogonal tableau} or a \textit{JK tableau} if 

    $$\tau_j \leq n+1 - n_j \quad \text{ for } \; 1 \leq j \leq \nu_2'.$$
    
    \noindent Let $\JK_n^{\nu^\pi}$ be the set of JK tableaux in $\B^{\nu^\pi}_n$ and $\JK^\lambda_{\nu^\pi \mu} = \JK_n^{\nu^\pi} \cap \HC^\lambda_{\nu^\pi \mu}(\GL_n)$.
    
\end{defn}

\noindent Jang and Kwon denote the set $\JK^\lambda_{\nu^\pi \mu}$  by $\underbar{\text{LR}}^\lambda_{\mu \nu}$.  They derive the following branching formula.

\begin{thm}[{\cite[Theorem 4.17]{KwonJang2021}}] \label{thm_KJ}
For $\lambda \in \Par_n$ and $\nu \in \widehat{\Or}_n$, 

$$\mult(\pi^\nu_{\Or_{n}}, \pi^\lambda_{\GL_{n}}) = \sum_{\mu \in \Par^{(2)}_n} |\JK^\lambda_{\nu^\pi \mu}|.$$
    
\end{thm}

In light of this result, it suffices to exhibit a bijection between 
$\JK_{\nu^\pi \mu}^\lambda$ and $\LC_{\mu \nu}^\lambda(\Or_n)$. Similar to the $\Sp_{2n}$ case, the map providing this bijection is a crystal anti-isomorphism, namely the rotation-complement map
\[
S^\pi:\B_n^{\nu^\pi}\to \B_n^\nu.
\]
By Theorem \ref{thm:S-bijection}, $S^\pi$ restricts to a bijection
\[
\HC_{\nu^\pi \mu}^\lambda(\GL_n)
\to
\LC_{\mu \nu}^\lambda(\GL_n).
\]
Thus it remains to show that $S^\pi$ restricts to a bijection
\[
\JK_n^{\nu^\pi}\to \mathcal{T}_{\Or_n}^\nu.
\]
Assuming this, we have
\begin{align*}
    S^\pi(\JK^\lambda_{\nu^\pi \mu}) 
    &= S^\pi(\JK_n^{\nu^\pi} \cap \HC^\lambda_{\nu^\pi \mu}(\GL_n)) \\
    &= S^\pi(\JK_n^{\nu^\pi}) \cap S^\pi(\HC^\lambda_{\nu^\pi \mu}(\GL_n)) \\
    &= \mathcal{T}^\nu_{\Or_n} \cap \LC^\lambda_{\mu \nu}(\GL_n) \\
    &= \LC^\lambda_{\mu \nu}(\Or_n).
\end{align*}

To this end, we begin by applying $S^\pi$ to $\JK_n^{\nu^\pi}$ and characterizing the image in $\B_n^\nu$. This is simply a matter of keeping track of how the conditions of Definition \ref{def_kwon} transform under $S^\pi$. Recall, $S^\pi$ rotates $\nu^\pi$ by $180^\circ$ so it becomes $\nu$ and sends the entry $i \mapsto n+1 - i$. It is easy to check that Definition \ref{def_kwon} becomes,

\begin{defn}
    Let $T \in \B^\nu_n$ and if $\nu_1' >n/2$, set $r = n - \nu_1'$. Otherwise, $r=\nu_1'$. For $1 \leq i \leq \nu_1'$ and $1 \leq j \leq \nu_2'$, let 

    \begin{itemize}
        \item $\sigma_i = $ the $i$th entry in the leftmost column of $T$ from the top,
        \item $\tau_j = $ the $j$th entry in the second leftmost column of $T$ from the top,
        
        \item \begin{equation*}
          m_i =
        \begin{cases}
      \min\{\sigma_i, 2i-1\} & \text{if $1 \leq i \leq r$,}\\
      \min\{\sigma_i, r +i \} & \text{if $r < i \leq \nu'_1$}
        \end{cases}       
        \end{equation*}
        
        \item $n_j = $ the $j$th smallest number in $\{j+1, \dots , n\} \setminus \{m_{j+1}, \dots , m_{\nu_1'} \}$.
    \end{itemize}

    \noindent Say $T$ is a \textit{preliminary rotated Jang-Kwon orthogonal tableau} or a \textit{pRJK tableau} if 

    $$\tau_j \geq  n_j \quad \text{ for } \; 1 \leq j \leq \nu_2'.$$
    
    \noindent Let $\pRJK_n^{\nu}$ be the set of pRJK tableaux in $\B^\nu_n$.
\end{defn}

By construction, $S^\pi$ restricts to a bijection $\JK_n^{\nu^\pi} \to \pRJK_n^\nu$. To finish the proof, we need to show $\pRJK_n^{\nu}= \mathcal{T}_{\Or_n}^\nu$. Before doing that, we replace the somewhat involved pRJK condition by a simpler condition involving only the complement of the entries in the first column.

\begin{defn} \label{def:rjk}
    Let $T \in \B^\nu_n$. For $1 \leq i \leq \nu_1'$ and $1 \leq j \leq \nu_2'$, let 

    \begin{itemize}
        \item $\sigma_i = $ the $i$th entry in the leftmost column of $T$ from the top,
        \item $\tau_j = $ the $j$th entry in the second leftmost column of $T$ from the top,
        
        \item $\beta_j = $ the $j$th smallest number in $[n] \setminus \{\sigma_{1}, \dots , \sigma_{\nu_1'} \}$.
    \end{itemize}

    \noindent Say $T$ is a \textit{rotated Jang-Kwon orthogonal tableau} or a \textit{RJK tableau} if 

    $$\tau_j \geq  \beta_j \quad \text{ for } \; 1 \leq j \leq \nu_2'.$$
    
    \noindent Let $\RJK_n^{\nu}$ be the set of RJK tableaux in $\B^\nu_n$.
\end{defn}

Thus, once the first column is fixed, the RJK condition says that the entries in the second column must occur no earlier than the available values in the complement
\[[n]\setminus\{\sigma_1,\dots,\sigma_{\nu_1'}\}.\]
The next lemma shows that this simpler condition is equivalent to the pRJK
condition.

\begin{lemma}
    $\pRJK_n^\nu=\RJK_n^\nu$.
\end{lemma}

\begin{proof}
    We begin by showing that the computation of $m_i$ in the definition of pRJK tableaux is unnecessary: we can just set $m_i = \sigma_i$. We do this in two steps. 
    
    First, we show that we can take $m_i = \min \{\sigma_i, 2i -1 \}$ for all $i \in [\nu_1']$. We do not need a separate case for $i> r$. To see this, for any $i > r$, write  $i = \nu_1' - j $ where $j$ is the number of boxes in the first column strictly below $\sigma_i$ and notice $r+ i = (n-\nu_1') + (\nu_1'-j) =n-j$. Since columns are strictly increasing, there are $j$ terms of $[n]$ that are greater than $\sigma_i$, hence $\sigma_i \leq n-j = r+ i$. Since $2i-1 \geq r+i$ for all $i > r$, there is no harm in replacing $\min\{\sigma_i, r+i \}$ with $\min \{\sigma_i, 2i -1 \}$. Thus we may replace the piecewise definition of $m_i$ by
    \[
    m_i=\min\{\sigma_i,2i-1\}
    \]
    for all $i\in[\nu_1']$.

   Second, we show that we may replace $m_i$ by $\sigma_i$ when $2i-1<\sigma_i$. Assume $2i-1<\sigma_i$, so that $m_i=2i-1$. We claim that for every $j$,
\[
n_j
=
\text{the $j$th smallest number in }
\{j+1,\dots,n\}\setminus
\{m_{j+1},\dots, m_{i-1}, \sigma_i, m_{i+1}, \dots,m_{\nu_1'}\}.
\]
    
  If $j \geq i$, then $m_i \notin \{m_{j+1},\dots,m_{\nu_1'}\}$, so $m_i$ does not affect the computation of $n_j$. Hence, there is no harm in replacing $m_i$ with $\sigma_i$. Now suppose $j<i$ and write $i=j+k$ with $k>0$. 

Since $j\leq \nu_2'$ and $\nu\in\widehat{\Or}_n$, we have $2j\leq 2\nu_2'\leq \nu_1'+\nu_2'\leq n$. Thus the $j$th smallest element of $\{j+1,\dots,n\}$ is $2j$. Observe that the sequence $m_1, m_2, \dots , m_{\nu_1'}$ is strictly increasing: $\sigma_i < \sigma_{i+1}$ since columns are strictly increasing and $2i-1 < 2i+1$. Hence, the only elements of
\(
\{m_{j+1},\dots,m_{\nu_1'} \}
\)
that are smaller than $m_i$ are among
\(
m_{j+1},\dots,m_{i-1}.
\)
Since $i-1 = j+k-1$, we see there are at most $k-1$ such elements. Therefore, after removing all elements smaller than $m_i$ from $\{j+1,\dots,n\}$, the $j$th smallest remaining element is at most
\[
2j+(k-1)=2j+k-1.
\]
But
\[
m_i=2i-1=2j+2k-1,
\]
so
\[
2j+k-1 < m_i.
\]
Thus the $j$th smallest remaining element is already strictly smaller than $m_i$, and removing $m_i$ cannot change $n_j$. Since $m_i=2i-1<\sigma_i$ by assumption, replacing $m_i$ by $\sigma_i$ does not change $n_j$. The elements $m_{i+1},m_{i+2},\dots$ are larger than $m_i$, so they also do not affect this conclusion. Therefore we obtain the same set of tableaux if we set $m_i=\sigma_i$ for all $i$.

    Notice, the only place where $r$ is used is in the computation of $m_i$, so we can drop $r$ from the characterization as well. At this point, we see that $n_j$ can be computed as the $j$th smallest element in 
    \[A_j=\{j+1, \dots, n\} \setminus \{\sigma_{j+1}, \dots , \sigma_{\nu_1'} \}.\]
    It remains to show that we can replace $n_j$ with $\beta_j$, the $j$th smallest element in 
    \[B=[n] \setminus \{\sigma_1, \dots , \sigma_{\nu_1'} \}.\]

    First observe both $A_j$ and $B$ have the same cardinality $n - \nu'_1$: since the elements $\sigma_1, \dots , \sigma_{\nu'_1}$ are strictly increasing in $[n]$ by the semistandard condition, $|B| = n - \nu_1'$ and $|A_j| = (n-j) - (\nu_1' - j) = n - \nu_1'$. Notice that
    \[A_j=    (B \setminus\{1,\dots,j\}) \cup (\{\sigma_1,\dots,\sigma_j\}\cap\{j+1,\dots,n\}).\]  
    Indeed, both sides consist of the elements of $\{j+1,\dots,n\}$ that are not among $\sigma_{j+1},\dots,\sigma_{\nu_1'}$. Passing from $B$ to $A_j$ removes some elements of $B$ that are at most $j$ and replaces them with some of the elements $\sigma_1,\dots,\sigma_j$ that are greater than $j$. Hence, $\beta_j \leq n_j$. These replacement elements are all at most $\sigma_j$. Therefore, $n_j \leq \max \{\beta_j, \sigma_j\}$. We conclude that
    \[\beta_j\leq n_j\leq \max\{\beta_j,\sigma_j\}.\]
    Since $T$ is semistandard, the entries in each row are weakly increasing. Hence $\tau_j\geq \sigma_j$.
    Now suppose $\tau_j\geq \beta_j$. Then
    \[\tau_j\geq \max\{\beta_j,\sigma_j\}\geq n_j.\]
    Conversely, if $\tau_j\geq n_j$, then $\tau_j\geq \beta_j$ because $\beta_j\leq n_j$. Thus
    \[\tau_j\geq \beta_j \quad\Longleftrightarrow\quad \tau_j\geq n_j,    \]
    as desired.
    
\end{proof}

The following lemma completes the proof that $S^\pi : \JK_n^{\nu^\pi} \to \mathcal{T}^\nu_{\Or_n}$ is a bijection and hence establishes the $\Or_n$ branching rule in Theorem \ref{thm:main_branching}.

\begin{lemma}
    $\RJK_n^\nu=\mathcal{T}_{\Or_n}^\nu$.
\end{lemma}

\begin{proof}
    Let $T\in \B_n^\nu$. Let
    \[
    [n]\setminus\{\sigma_1,\dots,\sigma_{\nu_1'}\}
    =
    \{\beta_1<\beta_2<\cdots\}.
    \]

    The $\Or_n$-tableau condition says that, for every $i\in[n]$, the number of entries $\leq i$ in the first two columns is at most $i$. Once the first column is fixed, the values $\sigma_1,\dots,\sigma_{\nu_1'}$ are already occupied. Therefore the second column can contribute at most as many entries $\leq i$ as there are elements $\leq i$ in the complement
    \[
    [n]\setminus\{\sigma_1,\dots,\sigma_{\nu_1'}\}.
    \]
    Equivalently, the $j$th entry in the second column must occur no earlier than the $j$th available value in this complement, i.e.
    \[
    \tau_j\geq \beta_j
    \qquad
    \text{for all }1\leq j\leq \nu_2'.
    \]

    Indeed, if $\tau_j<\beta_j$, then at $i=\tau_j$ the second column has at least $j$ entries $\leq i$, while the complement has fewer than $j$ entries $\leq i$, violating the $\Or_n$-tableau condition. Conversely, if $\tau_j\geq \beta_j$ for all $j$, then for every $i$, the number of second-column entries $\leq i$ is bounded by the number of complement entries $\leq i$, so the $\Or_n$-tableau condition holds.
\end{proof}

Applying $S^\pi$ to $\mathcal{T}^\nu_{\Or_n}$, we get a simpler characterization of the JK-tableaux,

\begin{cor}
    $T \in \B^{\nu^\pi}_n$ is a JK-tableau if and only if 
    $$\# \{ \text{boxes in rightmost two columns with entry} \geq n+1 -i \} \leq i$$
    for all $i \in [n]$.
\end{cor}

\section{Graded multiplicities of \texorpdfstring{$K$}{K}-types} \label{sec:graded-mults}
Recall from Definition \ref{def:M-weight delta}, $(\mathcal T^\nu_K)_0$ are the tableaux representing a basis for the zero $M$-weight space in $\pi^\nu_K$ and $m^{\nu,0}_{(G,K)}(q) = \mult_q(\pi^\nu_K, \Ha)$. The main result of this section is the following theorem,
\begin{thm}\label{thm:graded mults}
    For $(\GL_n, \Or_n)$, $(\GL_n^2, \GL_n)$, and $(\GL_{2n}, \Sp_{2n})$, 
    $$m^{\nu,0}_{(G,K)}(q)  = \sum_{T \in (\mathcal{T}_K^\nu)_0} q^{d(T)},$$
    \noindent where 
    \begin{center}
\begin{tabular}{c|l}
\multicolumn{1}{c}{$(G,K)$} & \multicolumn{1}{c}{$d(T)$}\\
    \hline
    $(\GL_n, \Or_n)$ & $\displaystyle \frac{1}{2} \bigg( |\nu| + \sum_{i = 1}^{n-1} (n-i) \cdot \lceil \varphi_{i}(T) \rceil_{\mathrm{ev}} \bigg )$,\\
    $(\GL_n^2, \GL_n) $ & $\displaystyle  \sum_{i = 1}^{n-1} (n-i) \cdot \varphi_{i}(T)$,\\
    $(\GL_{2n}, \Sp_{2n})$ & $\displaystyle \frac{1}{2} \bigg(|\nu|+\sum_{i = 1}^{n-1} (2n-2i) \cdot \varphi_{2i}(T) \bigg )$,\\
\end{tabular}
\end{center}
\noindent and where \(\lceil a\rceil_{\mathrm{ev}}\) denotes the smallest even integer greater than or equal to \(a\).
\end{thm}

The $(\GL_n^2, \GL_n)$ case is an incarnation of the well-known charge statistic stated in terms of rational $\GL_n$-tableaux. As discussed in the introduction, for the split family $(\GL_n,\Or_n)$ these polynomials give the coefficients of the Hodge $K$-character of the tempered spherical principal series representation of $\GL(n,\R)$ with real infinitesimal character. The $\varphi_i$ can be computed combinatorially via the signature rule as in Examples \ref{ex:sig rule 1}, \ref{ex:sig rule 2}, \ref{ex:sig rule 3}, and \ref{ex:sig rule 4}.

Section \ref{section:O_4-examples} presents explicit computations of a number of examples for $(\GL_4, \Or_4)$, with the goal of recomputing the table of $q$-multiplicities found at \cite[p.~80]{Willenbring2000}. Section \ref{sect:degree-one-k-type} determines the graded multiplicity of the unique irreducible representation occurring in degree one of the harmonics. Section \ref{sec:proof-of-graded-mults} contains the proof of Theorem \ref{thm:graded mults}.

\subsection{\texorpdfstring{$\Or_4$}{O4} examples} \label{section:O_4-examples}
We begin with several examples for $\Or_4$ and tie them to the results found in \cite[p.~80]{Willenbring2000}. The results in \cite{Willenbring2000, WallachWillenbring2000} are stated for $(\SL_4, \Sor_4)$ while the results in this paper are for $(\GL_4, \Or_4)$. It is easy to translate our results to the $(\SL_n, \Sor_n)$ setting by remembering how $\Or_n$-representations restrict to $\Sor_n$. This is an advantage of working with $\Or_n$ instead of $\Sor_n$. In general, it is less clear how to translate results for $(\SL_n, \Sor_n)$ to results about $(\GL_n, \Or_n)$. We also note that the $q$-multiplicities seem to have nicer properties in the $(\GL_n, \Or_n)$ setting. For example, in all cases used below to recompute the table in \cite[p.~80]{Willenbring2000}, the $(\GL_4,\Or_4)$ $q$-multiplicities are palindromic, whereas the corresponding $(\SL_4,\Sor_4)$ table contains non-palindromic $q$-multiplicities.

We begin by recalling how irreducible representations of $\Or_4$ restrict to $\Sor_4$, see \cite[Theorem 5.5.24]{GoodmanWallach2009} or \cite[Theorem 19.22]{FultonHarris1991}. Write $\Or_4 \cong \{1 , g_0 \} \ltimes \Sor_4$ as in \cite[p.~275]{GoodmanWallach2009}. Let $P_{++}(\Sor_4) = \{\mu_1 \epsilon_1 + \mu_2 \epsilon_2 : \mu_1 \geq |\mu_2| \}$ be the dominant weights of $\Sor_4$. We think of dominant weights as tuples via the identification $\mu_1 \epsilon_1 + \mu_2 \epsilon_2 \longleftrightarrow (\mu_1, \mu_2)$. Let $\varpi_1 = \frac{1}{2}(\epsilon_1 - \epsilon_2)$ and $\varpi_2 = \frac{1}{2}(\epsilon_1 + \epsilon_2)$ be the fundamental weights for $\sor_4$. The element $g_0$ acts on the dominant weight lattice by $g_0 \cdot \epsilon_1 = \epsilon_1$ and $g_0 \cdot \epsilon_2 = -\epsilon_2$. If $\nu \in \widehat{\Or}_4$,  its \textit{associated partition} $\overline{\nu}$ is the partition obtained from $\nu$ by changing the length of the first column from $\ell(\nu)$ to $4 - \ell(\nu)$ and leaving all other columns unchanged.

\begin{thm} \label{thm:O4-irreps}
    Let $\nu \in \widehat{\Or}_4$.
    \begin{enumerate}
        \item If $\nu \neq \overline{\nu}$, then $\pi^\nu_{\Or_4}|_{\Sor_4}$ is irreducible and $\pi^\nu_{\Or_4}|_{\Sor_4} \cong \pi^{\overline{\nu}}_{\Or_4}|_{\Sor_4}$.
        \item If $\nu = \overline{\nu}$, then $\pi^\nu_{\Or_4}|_{\Sor_4} \cong \pi^\nu_{\Sor_4} \oplus \pi^{g_0 \cdot \nu}_{\Sor_4}$.
    \end{enumerate}
\end{thm}

Case (1) says that the hook shape $\nu = (k, 1,1,0)$ restricts to the same $\Sor_4$-representation as $\overline{\nu} = (k, 0, 0,0)$, namely the representation with highest weight $k \epsilon_1 = (k, 0)$. Furthermore, the trivial representation $(0,0,0,0)$ restricts to the same representation as the determinant $(1,1,1,1)$, the trivial representation $(0,0)$ of $\Sor_4$. Case (2) says that $\pi^{(k_1, k_2, 0,0)}_{\Or_4}$ restricts to  $\pi^{(k_1, k_2)}_{\Sor_4} \oplus  \pi^{(k_1, -k_2)}_{\Sor_4}$.

Observe that $\nu$ can only have $M$-weight zero tableaux if $|\nu|$ is even since the filling must have an even number of $i$ for all $i \in [4]$. We start with the simplest nontrivial example.

\begin{example}
    Consider $\nu = (2, 0,0,0) = \ydiagram{2}$. The zero $M$-weight tableaux are 
    $$(\mathcal{T}^{\nu}_{\Or_4})_0 = \{ \ytableaushort{2 2}, \ytableaushort{3 3}, \ytableaushort{4 4} \}.$$

    Next, we compute the grading. By the signature rule, 
    $$\varphi_i(\ytableaushort{2 2}) = \begin{cases}
        2 & \text{if } i=2,\\
        0 & \text{otherwise.}
    \end{cases}$$

    So,
\begin{align*}
    d(\ytableaushort{2 2}) &= \frac{1}{2}(|\ydiagram{2}| + \sum_{i=1}^3 (4-i) \cdot \lceil \varphi_i(\ytableaushort{2 2}) \rceil_{\mathrm{ev}} )\\
    &= \frac{1}{2}(2 + (4-2) \cdot 2)\\
    &=3.
\end{align*}
    Similarly, we find $d(\ytableaushort{3 3}) = 2$ and $d(\ytableaushort{4 4}) = 1$. Thus, the graded multiplicity is $q^1 + q^2 + q^3$.

    Observe also that $\overline{\nu} = (2,1,1,0) = \ydiagram{2,1,1}$ has no $M$-weight zero tableaux, so the $\Sor_4$-representation indexed by $(2,0) = 2 \epsilon_1 = 2\varpi_1 + 2 \varpi_2$ has the same graded multiplicity $q^1+q^2+q^3$. This matches the entry in row 1, column 1 of the table at \cite[p.~80]{Willenbring2000}. We denote the entry in row $i$, column $j$ of the table by $T_{i,j}$. Note the table in \cite{Willenbring2000} has rows indexing even powers of $\varpi_1$ and columns indexing even powers of $\varpi_2$, so $T_{1,1}$ corresponds to $2 \varpi_1 + 2 \varpi_2 = 2 \epsilon_1 = (2,0)$. 
\end{example}

\begin{example}
    Consider $\nu = (2,2,0,0) = \ydiagram{2,2}$. The zero $M$-weight tableaux are 
    $$(\mathcal{T}^{\nu}_{\Or_4})_0 = \{ \ytableaushort{2 2, 4 4}, \ytableaushort{3 3,4 4} \}.$$

    Next, we compute the grading. 
    $$\varphi_i(\ytableaushort{2 2, 4 4}) = \begin{cases}
        2 & \text{if } i=2,\\
        0 & \text{otherwise.}
    \end{cases}$$

    So,
\begin{align*}
    d(\ytableaushort{2 2, 4 4}) &= \frac{1}{2}(|\ydiagram{2,2}| + \sum_{i=1}^3 (4-i) \cdot \lceil \varphi_i(\ytableaushort{2 2, 4 4}) \rceil_{\mathrm{ev}} )\\
    &= \frac{1}{2}(4 + (4-2) \cdot 2)\\
    &=4.
\end{align*}
Similarly, we find $d(\ytableaushort{3 3, 4 4}) = 2$. Thus, the graded multiplicity is $q^2 + q^4$. Since $\nu = \overline{\nu}$, $\pi^\nu_{\Or_4}|_{\Sor_4} \cong \pi^{(2,2)}_{\Sor_4} \oplus \pi^{(2, -2)}_{\Sor_4}$. Hence, we see this graded multiplicity in entries $T_{0,2}$ and $T_{2,0}$ of \cite[p.~80]{Willenbring2000}.
    
\end{example}

\begin{example}
    The final $\Or_4$ example we work out in detail is $\nu = (4, 0,0,0) = \ydiagram{4}$. The zero $M$-weight tableaux are 
    \begin{align*}
        (\mathcal{T}^{\nu}_{\Or_4})_0 = \{ \ytableaushort{2 2 2 2}, \ytableaushort{2 2 3 3 },& \ytableaushort{2 2 4 4},\\
        \ytableaushort{3 3 3 3},& \ytableaushort{3 3 4 4 },\\
        & \ytableaushort{4 4 4 4 } \}.
    \end{align*}

Computing the degrees as above gives the graded multiplicity
$q^2 + q^3 + 2q^4 +q^5 + q^6$. Notice, the associated partition $\overline{\nu} = (4, 1, 1, 0) = \ydiagram{4, 1, 1}$ does have one zero $M$-weight tableau $T$,
$$(\mathcal{T}^{\overline{\nu}}_{\Or_4})_0 = \{ \ytableaushort{2 2 3 4, 3, 4} \}.$$
Recalling that $\lceil \cdot\rceil_{\mathrm{ev}}$ means round up to the nearest even integer,
\begin{align*}
    d(T) &= \frac{1}{2}(|\overline{\nu}| + (4-2)\lceil \varphi_2(T) \rceil_{\mathrm{ev}} + (4-3) \lceil \varphi_3(T) \rceil_{\mathrm{ev}} )\\
    &= \frac{1}{2}(6 + 2 \cdot\lceil 1 \rceil_{\mathrm{ev}} + 1 \cdot \lceil 1 \rceil_{\mathrm{ev}} )\\
    &= \frac{1}{2}(6 + 2 \cdot 2 + 1 \cdot 2)\\
    & = 6.
\end{align*}

As discussed above, the $\Or_4$-representation indexed by this associated partition also restricts to the $\Sor_4$-representation with highest weight $(4,0)$. Hence, the graded multiplicity for the $\Sor_4$-representation $(4,0) = 4 \varpi_1 + 4 \varpi_2$ is equal to the sum of the graded multiplicities of the $\Or_4$ irreducible representations indexed by $\nu$ and $\overline{\nu}$: 
$$(q^2 + q^3 + 2q^4 +q^5 + q^6) + (q^6) = q^2 + q^3 + 2q^4 +q^5 + 2q^6.$$

Again, this matches entry $T_{2,2}$ in \cite[p.~80]{Willenbring2000}. This also provides an illustration of how restricting to $\Sor_4$ can break palindromicity. 
\end{example}

Next, we present tables for the cases considered in \cite[p.~80]{Willenbring2000}. We computed the $q$-multiplicities in these tables by hand, but it would be easy to program as well. The tables are organized as in \cite[p.~80]{Willenbring2000}, by even multiples of the fundamental weights of $\sor_4$, so entry $(i,j)$ corresponds to highest weight $2i \varpi_1 + 2j \varpi_2$. 

The first table below gives the $q$-multiplicities of the $\Or_4$ irreducible representations, indexed by partitions with at most two rows, whose restrictions contain the indicated irreducible representation of $\Sor_4$. The entries strictly above the main diagonal correspond to two-row shapes $\nu=(k_1,k_2,0,0)$, and hence satisfy $\nu=\overline{\nu}$. By Theorem \ref{thm:O4-irreps}, each such $\Or_4$ irreducible representation restricts to a direct sum of two irreducible $\Sor_4$-representations. In the table below, we record the graded multiplicity of each such $\Or_4$-representation only once; the second summand appears later as the reflection across the main diagonal. The main diagonal corresponds to one-row shapes. Among the shapes recorded in this table, these are precisely the shapes satisfying $\nu\neq \overline{\nu}$.

\begin{center}
    \begin{tabular}{c|c|c|c|c}
         & 0 & 1 & 2&3 \\
         \hline
        0 & 1 & 0 & $q^2+q^4$ & $q^6$\\
        1 & - & $q + q^2 + q^3$ & $q^3+q^4+q^5$ & $q^3+q^4+ 2q^5 + q^6 + q^7$\\
        2 & - & - & $q^2 + q^3 + 2 q^4 + q^5 + q^6$ & $q^4 + 2q^5 + 2q^6 + 2 q^7 + q^8$\\
        3 & - & - & - & $q^3 + q^4 + 2 q^5 + 2q^6 +2q^7 +q^8 + q^9$
        
    \end{tabular}
\end{center}

For the one-row shapes $\nu$ appearing on the main diagonal of the previous table, the next table gives the $q$-multiplicities of the associated partitions $\overline{\nu}$.
\begin{center}
    \begin{tabular}{c|c|c|c|c}
         & 0 & 1 & 2&3 \\
         \hline
        0 & -&-&-&-\\
        1 & - & $0$ & - & -\\
        2 & - & - & $q^6$ & -\\
        3 & - & - & - & $q^7 + q^8 + q^9$
        
    \end{tabular}
\end{center}

Restricting the $\Or_4$-representations to $\Sor_4$ has the effect of reflecting the first table across the main diagonal, producing a symmetric matrix; this is Theorem \ref{thm:O4-irreps}, part (2). We must also add the contributions from the second table to the diagonal entries, by Theorem \ref{thm:O4-irreps}, part (1). These additional contributions destroy palindromicity for some of the $q$-multiplicities along the main diagonal, corresponding to the one-row shapes. The resulting $\Sor_4$ $q$-multiplicities are displayed below and agree exactly with the table in \cite[p.~80]{Willenbring2000}.

\begin{center}
\small
    \begin{tabular}{c|c|c|c|c}
         & 0 & 1 & 2&3 \\
         \hline
        0 & 1 & 0 & $q^2+q^4$ & $q^6$\\
        1 & 0 & $q + q^2 + q^3$ & $q^3+q^4+q^5$ & $q^3+q^4+ 2q^5 + q^6 + q^7$\\
        2 & $q^2 + q^4$ & $q^3+q^4+q^5$ & $q^2 + q^3 + 2 q^4 + q^5 + 2q^6$ & $q^4 + 2q^5 + 2q^6 + 2 q^7 + q^8$\\
        3 & $q^6$ & $q^3+q^4+ 2q^5 + q^6 + q^7$ & $q^4 + 2q^5 + 2q^6 + 2 q^7 + q^8$ & $q^3 + q^4 + 2 q^5 + 2q^6 +3q^7 +2q^8 + 2q^9$
    \end{tabular}
\end{center}

\subsection{The degree-one \texorpdfstring{$K$}{K}-type} \label{sect:degree-one-k-type}
Next, we compute the graded multiplicity of the unique irreducible representation appearing in degree one of $\Ha$. Since the grading comes from $\C[\p]$, we see the irreducible representation is just the non-invariant piece of $\p$. 

\begin{example}
    For $\GL_n$, $\p = \gl_n = \C \oplus \sll_n$ under the adjoint action of $\GL_n$. The non-invariant piece, $\sll_n$, has highest weight the highest root $\nu = \epsilon_1 - \epsilon_n$ which corresponds to the staircase 
    $$\begin{ytableau}
    \none & \none[1] &\\
    \none & \none[2]\\
    \none & \none[\vdots]\\
    \none & \none[n-1]\\
     & \none[n]\\
\end{ytableau}$$

The $M$-weight zero tableaux are
$$(\mathcal{T}^\nu_{\GL_n})_0 =\begin{ytableau}
    \none & \none[1] & 2\\
    \none & \none[2]\\
    \none & \none[\vdots]\\
    \none & \none[n-1]\\
     \overline{2} & \none[n]\\
\end{ytableau}, \begin{ytableau}
    \none & \none[1] & 3\\
    \none & \none[2]\\
    \none & \none[\vdots]\\
    \none & \none[n-1]\\
     \overline{3} & \none[n]\\
\end{ytableau}, \dots \begin{ytableau}
    \none & \none[1] & n\\
    \none & \none[2]\\
    \none & \none[\vdots]\\
    \none & \none[n-1]\\
     \overline{n} & \none[n]\\
\end{ytableau}.$$

\noindent The filling by $\overline{1}, 1$ is not a $\GL_n$-tableau. Let $T_k$ denote the tableau above with entries $\overline{k}$ and $k$. $T_k$ has word $\overline{k} k$. It is easy to check using the signature rule that 
$$\varphi_i(T_k) = \begin{cases}
    1 & \text{if } i = k-1,\\
    0 & \text{otherwise.}
\end{cases}$$ 

\noindent It follows that $T_k$ contributes to degree $n+1-k$ and the graded multiplicity is $q^1 + q^2 + \cdots + q^{n-1}$.
\end{example}

\begin{example}
    Similarly, we compute the graded multiplicity of $\nu = \ydiagram{2}$ for $\Or_n$. The relevant $M$-weight zero tableaux are
    $$(\mathcal{T}^\nu_{\Or_n})_0 = \{ \ytableaushort{ 2 2}, \ytableaushort{3 3}, \cdots , \ytableaushort{n n} \}. $$
    \noindent and for $1 \leq i \leq n-1$
    $$\varphi_i(\ytableaushort{k k}) = \begin{cases}
    2 & \text{if } i = k,\\
    0 & \text{otherwise.}
\end{cases}$$ 

It follows that the graded multiplicity is $q^1 + q^2 + \cdots + q^{n-1}$.
\end{example}

\begin{example}
    Finally, we compute the graded multiplicity of $\nu = \ydiagram{1,1}$ for $\Sp_{2n}$. The relevant $M$-weight zero tableaux are
    $$(\mathcal{T}^\nu_{\Sp_{2n}})_0 = \{ \ytableaushort{3,4}, \ytableaushort{5,6}, \dots , \ytableaushort{{\scalebox{0.5}{$2n-1$}},{2n}}  \}.$$
    From this, we compute the graded multiplicity to be $q^1 + q^2 + \cdots + q^{n-1}$.
\end{example}

\subsection{Proof of Theorem \ref{thm:graded mults}}\label{sec:proof-of-graded-mults}
We start with
$$\mult_q(\pi^\nu_K, \C[\p]) = \sum_{\lambda \in \Spec(\C[\p])} q^{d_G^\lambda}b^\lambda_\nu$$
\noindent from the introduction and apply the relevant branching rules in Theorem \ref{thm:main_branching} and Proposition \ref{prop:GLn branching rule} to the multiplicity-free decompositions in Theorem \ref{thm:mult-free spaces}.

This gives the following table

\begin{center}
\begin{tabular}{c|l}
\multicolumn{1}{c}{$(G,K)$} & \multicolumn{1}{c}{$\mult_q(\pi^\nu_K, \C[\p])$}\\
    \hline
    $(\GL_n, \Or_n)$ & $\displaystyle \sum_{\lambda, \mu \in \Par^{(2)}_n} |\LC^\lambda_{\mu \nu}(\Or_n)| \cdot q^{|\lambda|/2} $\\
    $(\GL_n^2, \GL_n) $ & $\displaystyle \sum_{\lambda, \mu \in \Par_n, \lambda = \mu} |\LC^\lambda_{\mu \nu}(\GL_n)| \cdot q^{|\lambda|}$\\
    $(\GL_{2n}, \Sp_{2n})$ & $\displaystyle \sum_{\lambda, \mu \in \Par^{(1,1)}_{2n}} |\LC^\lambda_{\mu \nu}(\Sp_{2n})|\cdot q^{|\lambda|/2}$\\
\end{tabular}
\end{center}

For the second row, we use the symmetry $c^\nu_{\lambda^* \lambda} = c^\lambda_{\lambda \nu}$ from Section \ref{subsec:LR-symmetries} to write $\sum_{\lambda \in \Par_n} c^\nu_{\lambda^* \lambda} \cdot q^{|\lambda|} = \sum_{\lambda \in \Par_n} c^\lambda_{\lambda \nu} \cdot q^{|\lambda|}$. Then we set $\mu = \lambda$ so our notation is uniform with the other two cases. The remaining two rows are immediate. Let $Q$ denote the relevant indexing set of partitions in the summation, that is $Q = \Par^{(2)}_n, \Par_n$, or  $\Par^{(1,1)}_{2n}$. Set $N = n$ for $(\GL_n, \Or_n)$, $(\GL_n^2, \GL_n)$ and $N = 2n$ for $(\GL_{2n}, \Sp_{2n})$. Let \(I_K\) be the indexing set for the sums above
\[
I_K =
\begin{cases}
\{(\lambda,\mu):\lambda,\mu\in\Par_n^{(2)}\}, &(\GL_n,\Or_n),\\
\{(\lambda,\mu):\lambda=\mu\in\Par_n\}, & (\GL_n^2,\GL_n),\\
\{(\lambda,\mu):\lambda,\mu\in\Par_{2n}^{(1,1)}\}, &(\GL_{2n},\Sp_{2n}).
\end{cases}
\]

The key is that $\nu$ is fixed in the sums, so each $\LC_{\mu \nu}^\lambda(K)$ appearing may be viewed as a subset of $\mathcal{T}^\nu_{K}$. We can rewrite the sums above over the set 
$$\Mult^\nu_K = \bigsqcup_{(\lambda,\mu)\in I_K} \LC^\lambda_{\mu \nu}(K).$$
We write an element of $\Mult^\nu_K$ as $(\lambda, \mu, T)$, where $(\lambda, \mu ,T) \in \Mult^\nu_K$ if and only if $T \in \LC^\lambda_{\mu \nu}(K)$. Given $(\lambda, \mu ,T) \in \Mult^\nu_K$, $\lambda$ is determined by $T$ and $\mu$. For $(\GL_n^2, \GL_n)$, $\lambda= \mu$ so this is immediate. For $(\GL_n, \Or_n)$ and $(\GL_{2n}, \Sp_{2n})$, $T \in \LC^\lambda_{\mu \nu}(K)$, so $L^\mu_N \otimes T \equiv L^\lambda_N$ and $\lambda = \mu + w_0(\wt(T))$ where $w_0$ is the longest Weyl group element for $\GL_N$. In either case, write this uniquely determined $\lambda$ as $\lambda(\mu, T)$. Thus, an arbitrary element of $\Mult^\nu_K$ can be written as $(\lambda(\mu, T), \mu ,T)$ and the pair $(\mu, T)$ uniquely determines the element of the disjoint union. Going forward, we will think of the elements of $\Mult^\nu_K$ as tuples $(\mu, T)$. With this, we can rewrite the sum determining $\mult_q(\pi^\nu_K, \C[\p])$ as

\begin{center}
\begin{tabular}{c|l}
\multicolumn{1}{c}{$(G,K)$} & \multicolumn{1}{c}{$\mult_q(\pi^\nu_K, \C[\p])$}\\
    \hline
    $(\GL_n, \Or_n)$ & $\displaystyle \sum_{(\mu, T) \in \Mult^\nu_{\Or_n}} q^{|\lambda(\mu, T)|/2} $\\
    $(\GL_n^2, \GL_n) $ & $\displaystyle \sum_{(\mu, T) \in \Mult^\nu_{\GL_n}} q^{|\lambda(\mu, T)|}$\\
    $(\GL_{2n}, \Sp_{2n})$ & $\displaystyle \sum_{(\mu, T) \in \Mult^\nu_{\Sp_{2n}}} q^{|\lambda(\mu, T)|/2}$\\
\end{tabular}
\end{center}

Now, we want to factor out the invariants from this sum. To do this, we consider the map $p: \Mult^\nu_K \to \mathcal{T}^\nu_K$ projecting onto the $T$ factor, i.e. $p(\mu, T) = T$. First, we characterize the image of this map.

\begin{prop} \label{prop:image of p}
    The image of $p$ is $(\mathcal{T}^\nu_K)_0$, the set of $M$-weight zero elements in $\mathcal{T}^\nu_K$.
\end{prop}

\begin{proof}
    Although the three cases are very similar, we proceed case by case. 
    
    $(\GL_n, \Or_n)$: First suppose $T$ is in the image of $p$. Then $T \in \LC^\lambda_{\mu \nu}(\Or_n)$ for some $\lambda, \mu \in \Par^{(2)}_n$. In terms of the $\GL_n$ weight lattice basis $\epsilon_1, \dots , \epsilon_n$, this means $\mu = \mu_1 \epsilon_1 + \cdots + \mu_n \epsilon_n$ and $\lambda = \lambda_1 \epsilon_1 + \cdots + \lambda_n \epsilon_n$ where all $\mu_i, \lambda_j$ are even nonnegative integers. The same holds for the lowest weights $w_0(\mu)$ and $w_0(\lambda)$, since $w_0$ just permutes the $\epsilon_i$. By Proposition \ref{prop:LC_equivalent_sets},
    \begin{align*}
        \wt(T) &= \wt(L^\lambda_n) - \wt(L^\mu_n)\\
        &= w_0(\lambda) - w_0(\mu)\\
        & = (\lambda_1 - \mu_1)\epsilon_n + (\lambda_2-\mu_2)\epsilon_{n-1} + \cdots + (\lambda_n - \mu_n) \epsilon_1
    \end{align*}

    \noindent must also have even coefficients in terms of the $\epsilon_i$ basis. Hence, by Definition \ref{def:K-tableaux}, $\mwt_{\Or_n}(T) = 0$.

    Now suppose $T \in (\mathcal{T}^\nu_{\Or_n})_0$, i.e. $\mwt_{\Or_n}(T) = 0$. We need to find $\mu \in \Par_n^{(2)}$ such that $T \in \LC^{\lambda(\mu,T)}_{\mu \nu}(\Or_n)$ and $\lambda(\mu, T) \in \Par_n^{(2)}$. We begin by finding a $\mu \in \Par_n^{(2)}$ such that $T \in \LC^{\lambda(\mu,T)}_{\mu \nu}(\Or_n)$. Since $\LC^{\lambda(\mu,T)}_{\mu \nu}(\Or_n) = \LC^{\lambda(\mu,T)}_{\mu \nu}(\GL_n) \cap \mathcal{T}^\nu_{\Or_n}$ and we know $T \in \mathcal{T}^\nu_{\Or_n}$, we just need to find a $\mu \in \Par^{(2)}_n$ such that $T \in \LC^{\lambda(\mu,T)}_{\mu \nu}(\GL_n)$. By Proposition \ref{prop:LC_equivalent_sets} (2), it suffices to find $\mu$ such that $\mu_0 \geq \varphi(T)^* = \varphi_1(T) \varpi_{n-1} + \varphi_2(T) \varpi_{n-2} + \cdots + \varphi_{n-1}(T) \varpi_1$. Taking any $\mu = m_1 \gamma_1 + \cdots + m_{n} \gamma_n$ with $m_i$ even and $m_i \geq \varphi_{n-i}(T)$ for all $i \in [n-1]$ and setting $m_n =0$ works. (Recall from Section \ref{sec:conventions for gln}: $\gamma_i = \epsilon_1 + \cdots + \epsilon_i$).
    
    Now we show that for any such $\mu$, $\lambda(\mu, T) \in \Par_n^{(2)}$. Since $\mwt_{\Or_n}(T) = 0$, there must be an even number of $i$ in $T$ for all $i \in [n]$ by Definition \ref{def:K-tableaux}. Since $\mu \in \Par^{(2)}_n$, $\lambda(\mu, T) = \mu + w_0(\wt(T))$ has even rows, and hence lies in $\Par_n^{(2)}$. 

    $(\GL_n^{2}, \GL_n)$: First suppose $T$ is in the image of $p$. Then $T \in \LC^\mu_{\mu\nu}(\GL_n)$ for some $\mu \in \Par_n$. Hence $\mwt_{\GL_n}(T) = \wt(T) = \wt(L^{\mu}_n) - \wt(L^{\mu}_n) = 0$ by Proposition \ref{prop:LC rational equiv set}.

    Now suppose $T \in (\mathcal{T}^\nu_{\GL_n})_0$, i.e. $\mwt_{\GL_n}(T) = 0$. We need to find $\mu \in \Par_n$ such that $T \in \LC^\mu_{\mu \nu}(\GL_n)$. By Proposition \ref{prop:LC rational equiv set}, it suffices to find $\mu$ such that $\mu_0 \geq \varphi(T)^* = \varphi_1(T) \varpi_{n-1} + \varphi_2(T) \varpi_{n-2} + \cdots + \varphi_{n-1}(T) \varpi_1$. Taking any $\mu = m_1 \gamma_1 + \cdots + m_{n-1} \gamma_{n-1}$ with integers $m_i \geq \varphi_{n-i}(T)$ for all $i \in [n-1]$ works.

    $(\GL_{2n}, \Sp_{2n})$: First suppose $T$ is in the image of $p$. Then $T \in \LC^\lambda_{\mu \nu}(\Sp_{2n})$ for some $\lambda, \mu \in \Par_{2n}^{(1,1)}$ so $\wt(T) = \wt(L^{\lambda}_{2n}) - \wt(L^{\mu}_{2n})$ by Proposition \ref{prop:LC_equivalent_sets}. Now for $\lambda = \sum_{i=1}^n \lambda_{2i} \gamma_{2i} \in \Par^{(1,1)}_{2n}$, $\wt(L^{\lambda}_{2n}) = w_0(\lambda) =\sum_{i=1}^n \lambda_{2i} w_0(\gamma_{2i})$ where $w_0(\gamma_{2i}) = \gamma_{2n} - \gamma_{2n-2i} = \epsilon_{2n} + \epsilon_{2n-1} + \cdots +  \epsilon_{2n-2i+1}$, see the discussion above Theorem \ref{thm:GLn_highest_lowest}. Hence, the coefficients of $\epsilon_{2i}$ and $\epsilon_{2i-1}$ are equal in $\wt(L^\lambda_{2n})$ for all $i \in [n]$. Similarly for $\wt(L^{\mu}_{2n})$. This implies that in the filling of $T$, the number of entries equal to $2i$ is equal to the number of entries equal to $2i-1$ for all $i \in [n]$ so, by Definition \ref{def:K-tableaux}, $\mwt_{\Sp_{2n}}(T) = 0$. 
    
    It remains to show that the ballot condition is satisfied. Suppose $T$ is not a ballot tableau and let $w(T) = x_1 x_2 \cdots x_k$ be the word of $T$. Then, there exists an $i \in [n]$ and a position $j$ such that the final segment $x_j x_{j+1} \cdots x_k$ has $\# \{ 2i \} > \#\{2i-1 \}$. By the signature rule, $\varepsilon_{2i-1}(T) > 0$. Recall from the definition of a crystal, Definition \ref{def:crystals},
    $$\varphi_{2i-1}(T) = \langle \wt(T), h_{2i-1} \rangle + \varepsilon_{2i-1}(T).$$
    Since $\#\{2i\} = \# \{2i-1\}$ in $T$, $\langle \wt(T), h_{2i-1} \rangle =0$, and we see $\varphi_{2i-1}(T) = \varepsilon_{2i-1}(T) > 0$.
    Hence $\varphi(T)^*$ has a positive coefficient in front of the odd fundamental weight $\varpi_{2n-2i+1}$. No $\mu \in \Par_{2n}^{(1,1)}$ can have $\mu_0 \geq \varphi(T)^*$ since any such $\mu$ must have $\mu_0$ with only even-indexed fundamental weights. So by Proposition \ref{prop:LC_equivalent_sets} (2), there is no $\LC^\lambda_{\mu \nu}(\Sp_{2n})$ with $\mu \in \Par^{(1,1)}_{2n}$ containing $T$, contradicting the assumption that $T$ is in the image of $p$. It follows that the ballot condition must be satisfied.

    Now suppose $T \in (\mathcal{T}^\nu_{\Sp_{2n}})_0$, i.e. $\mwt_{\Sp_{2n}}(T) = 0$ and the ballot condition is satisfied. We need to find $\mu \in \Par_{2n}^{(1,1)}$ such that $T \in \LC^{\lambda(\mu,T)}_{\mu \nu}(\Sp_{2n})$ and $\lambda(\mu, T) \in \Par_{2n}^{(1,1)}$. We begin by finding a $\mu \in \Par_{2n}^{(1,1)}$ such that $T \in \LC^{\lambda(\mu,T)}_{\mu \nu}(\Sp_{2n})$. Since $\LC^{\lambda(\mu,T)}_{\mu \nu}(\Sp_{2n}) = \LC^{\lambda(\mu,T)}_{\mu \nu}(\GL_{2n}) \cap \mathcal{T}^{\nu,H}_{\Sp_{2n}}$ and we have assumed $T \in (\mathcal{T}^\nu_{\Sp_{2n}})_0 \subset \mathcal{T}^{\nu,H}_{\Sp_{2n}}$, we just need to find a $\mu \in \Par^{(1,1)}_{2n}$ such that $T \in \LC^{\lambda(\mu,T)}_{\mu \nu}(\GL_{2n})$. As we saw above, since $T$ has zero $M$-weight and satisfies the ballot condition, the odd $\varphi_{2i-1}(T)$ vanish, so $\varphi(T)^* = \sum_{i =1}^{n-1} \varphi_{2i}(T) \varpi_{2n -2i}$, i.e. no odd fundamental weights show up. Hence, any $\mu =\sum_{i=1}^{n-1} m_{2i} \gamma_{2i} \in \Par^{(1,1)}_{2n}$ with $m_{2i} \geq \varphi_{2n-2i}(T)$ for $i \in [n-1]$ satisfies $\mu_0 \geq \varphi(T)^*$, so $T \in \LC^{\lambda(\mu,T)}_{\mu \nu}(\Sp_{2n})$.

    Now we show that for any such $\mu$, $\lambda(\mu, T) \in \Par_{2n}^{(1,1)}$. Since $\mwt_{\Sp_{2n}}(T) =0$, the coefficients of $\epsilon_{2i}$ and $\epsilon_{2i-1}$ are equal for all $i \in [n]$. This continues to hold for $w_0(\wt(T))$. Since $\mu \in \Par^{(1,1)}_{2n}$, $\lambda(\mu, T) = \mu + w_0(\wt(T))$ has equal coefficients of $\epsilon_{2i}$ and $\epsilon_{2i-1}$ for all $i \in [n]$, so $\lambda(\mu, T) \in \Par^{(1,1)}_{2n}$.
\end{proof}

Let $F_T = \{ \mu : (\mu, T) \in \Mult^\nu_K \}$ be the fiber above $T$. Then we can rewrite the sums again as follows:

\begin{center}
\begin{tabular}{c|l}
\multicolumn{1}{c}{$(G,K)$} & \multicolumn{1}{c}{$\mult_q(\pi^\nu_K, \C[\p])$}\\
    \hline
    $(\GL_n, \Or_n)$ & $\displaystyle \sum_{T \in (\mathcal{T}^\nu_{\Or_n})_0} \sum_{\mu \in F_T} q^{|\lambda(\mu, T)|/2} $\\
    $(\GL_n^2, \GL_n) $ & $\displaystyle \sum_{T \in (\mathcal{T}^\nu_{\GL_n})_0} \sum_{\mu \in F_T} q^{|\lambda(\mu, T)|}$\\
    $(\GL_{2n}, \Sp_{2n})$ & $\displaystyle \sum_{T \in (\mathcal{T}^\nu_{\Sp_{2n}})_0} \sum_{\mu \in F_T} q^{|\lambda(\mu, T)|/2}$\\
\end{tabular}
\end{center}

Next, we establish the structure of the fibers $F_T$. From the proof of Proposition \ref{prop:image of p}, $\mu \in F_T$ if and only if $\mu_0 \geq \varphi(T)^*$ and $\mu \in Q$. In each of the three cases, there is a unique minimal $\mu \in Q$ such that $\mu_0 \geq \varphi(T)^*$, which we denote $\mu_{\min}(T)$. Write $\varphi(T)^* = \varphi_{N-1}(T) \varpi_1 + \cdots + \varphi_1(T) \varpi_{N-1}$ in terms of the fundamental weights. Recall that $N = n$ in the first two cases and $N=2n$ in the symplectic case. The $\mu_{\min}(T)$ for the cases we consider are below.

\begin{center}
\begin{tabular}{c|l}
\multicolumn{1}{c}{$(G,K)$} & \multicolumn{1}{c}{$\mu_{\min}(T)$}\\
    \hline
    $(\GL_n, \Or_n)$ & $\displaystyle \sum_{i=1}^{n-1} \lceil \varphi_{i}(T) \rceil_{\mathrm{ev}} \gamma_{n-i}$\\
    $(\GL_n^2, \GL_n)$ &$\displaystyle \sum_{i=1}^{n-1} \varphi_{i}(T) \gamma_{n-i}$ \\
    $(\GL_{2n}, \Sp_{2n})$ & $\displaystyle \sum_{i=1}^{n-1} \varphi_{2i}(T) \gamma_{2n-2i}$\\
\end{tabular}
\end{center}

 Only the $(\GL_n, \Or_n)$ case requires the coefficients of $\varphi(T)^*$ to be modified. Here \(\lceil a\rceil_{\mathrm{ev}}\) denotes the smallest even integer greater than or equal to \(a\). This ensures that $\mu_{\min}(T)$ is in $\Par^{(2)}_n$.

\begin{lemma}
    $F_T = \mu_{\min}(T) +Q.$
\end{lemma}

\begin{proof}
    If $\mu \in F_T$, then $\mu_0 \geq \varphi(T)^*$ and $\mu \in Q$. It follows that $\mu - \mu_{\min}(T) \in Q$ by the explicit description of $\mu_{\min}(T)$ above and the definition of $\geq$ (discussed near the end of Section \ref{sec:reductive Lie algebras} and Section \ref{sec:conventions for gln}). Hence, $\mu = \mu_{\min}(T) + \delta$ for some $\delta \in Q$. 
    
    For $\delta \in Q$, any $\mu = \delta + \mu_{\min}(T)$ is in $F_T$. To see this, notice $\mu_0 \geq \varphi(T)^*$ so $T \in \LC_{\mu \nu}^{\lambda(\mu, T)}(\GL_N)$. Also, $\mu \in Q$ since $Q$ is an additive monoid and $\lambda(\mu,T)$ is in $Q$ by the proof of Proposition \ref{prop:image of p}.
\end{proof}

Since $\lambda$ is uniquely determined by $(\mu, T)$, there is a minimum $\lambda$ corresponding to $\mu_{\min}(T)$. Set $\lambda_{\min}(T) = \lambda(\mu_{\min}(T),T) = \mu_{\min}(T) + w_0(\wt(T))$. Then for every $\mu = \mu_{\min}(T) +\delta$ with $\delta \in Q$, we have $\lambda(\mu, T) = \lambda(\mu_{\min}(T) + \delta, T) = \lambda_{\min}(T) + \delta$. Since $q^{|\lambda(\mu, T)|} =q^{|\lambda_{\min}(T) + \delta|} = q^{|\lambda_{\min}(T)|} \cdot q^{|\delta|}$, we can rewrite the sums as 

\begin{center}
\begin{tabular}{c|l}
\multicolumn{1}{c}{$(G,K)$} & \multicolumn{1}{c}{$\mult_q(\pi^\nu_K, \C[\p])$}\\
    \hline
    $(\GL_n, \Or_n)$ & $\displaystyle \sum_{T \in (\mathcal{T}^\nu_{\Or_n})_0} q^{|\lambda_{\min}(T)|/2} \sum_{\delta \in \Par_n^{(2)}} q^{|\delta|/2} $\\
    $(\GL_n^2, \GL_n) $ & $\displaystyle  \sum_{T \in (\mathcal{T}^\nu_{\GL_n})_0} q^{|\lambda_{\min}(T)|} \sum_{\delta \in \Par_n} q^{|\delta|} $\\
    $(\GL_{2n}, \Sp_{2n})$ & $\displaystyle  \sum_{T \in (\mathcal{T}^\nu_{\Sp_{2n}})_0} q^{|\lambda_{\min}(T)|/2} \sum_{\delta \in \Par_{2n}^{(1,1)}} q^{|\delta|/2} $\\
\end{tabular}
\end{center}

We recognize the sum over $Q$ as the invariants and rewrite it as a rational function:
\begin{lemma}
    $$\sum_{\delta \in \Par_n^{(2)}} q^{|\delta|/2} = \sum_{\delta \in \Par_n} q^{|\delta|} = \sum_{\delta \in \Par_{2n}^{(1,1)}} q^{|\delta|/2} = \frac{1}{\prod_{i=1}^n (1-q^i)}.$$
\end{lemma}

\begin{proof}
    The additive monoids $\Par^{(2)}_n, \Par_n$, and $\Par^{(1,1)}_{2n}$ are freely generated by the following $n$ elements:
    \begin{center}
        \begin{tabular}{c|c}
           $\Par^{(2)}_n$  &  $(2,0,0, \dots ,0), (2,2,0, \dots ,0), \dots , (2,2,2, \dots ,2)$\\
           $\Par_n$  & $(1,0,0, \dots, 0), (1,1,0,\dots, 0), \dots , (1,1,1, \dots , 1)$\\
           $\Par^{(1,1)}_{2n}$ & $(1,1,0,\dots ,0), (1,1,1,1,0,\dots, 0), \dots , (1,1,1,\dots ,1)$
        \end{tabular}
    \end{center}

    These generators have degrees $1, 2, \dots n$ with respect to the degree function appearing in the lemma and the lemma follows.
\end{proof}

We arrive at a combinatorial analogue of the separation of variables in the Kostant-Rallis Theorem (Theorem \ref{thm:Kostant-Rallis}):

\begin{center}
\begin{tabular}{c|l}
\multicolumn{1}{c}{$(G,K)$} & \multicolumn{1}{c}{$\mult_q(\pi^\nu_K, \C[\p])$}\\
    \hline
    $(\GL_n, \Or_n)$ & $\displaystyle \frac{1}{\prod_{i=1}^n (1-q^i)}\sum_{T \in (\mathcal{T}^\nu_{\Or_n})_0} q^{|\lambda_{\min}(T)|/2}  $\\
    $(\GL_n^2, \GL_n) $ & $\displaystyle  \frac{1}{\prod_{i=1}^n (1-q^i)}\sum_{T \in (\mathcal{T}^\nu_{\GL_n})_0} q^{|\lambda_{\min}(T)|} $\\
    $(\GL_{2n}, \Sp_{2n})$ & $\displaystyle  \frac{1}{\prod_{i=1}^n (1-q^i)} \sum_{T \in (\mathcal{T}^\nu_{\Sp_{2n}})_0} q^{|\lambda_{\min}(T)|/2}$\\
\end{tabular}
\end{center}
Comparing this with $\C[\p]\cong \C[\p]^K\otimes \Ha$, the graded multiplicity in $\Ha$ is given by the numerator in each row.

From this, we recover the statistics on the zero $M$-weight tableaux presented in Theorem \ref{thm:graded mults} by expanding $|\lambda_{\min}(T)|$. 

\begin{center}
\begin{tabular}{c|l}
\multicolumn{1}{c}{$(G,K)$} & \multicolumn{1}{c}{$d(T)$}\\
    \hline
    $(\GL_n, \Or_n)$ & $\displaystyle \frac{1}{2}( |\nu| + \sum_{i = 1}^{n-1} (n-i) \cdot \lceil \varphi_{i}(T) \rceil_{\mathrm{ev}} )$\\
    $(\GL_n^2, \GL_n) $ & $\displaystyle  \sum_{i = 1}^{n-1} (n-i) \cdot \varphi_{i}(T)$\\
    $(\GL_{2n}, \Sp_{2n})$ & $\displaystyle \frac{1}{2}(|\nu|+\sum_{i = 1}^{n-1} (2n-2i) \cdot \varphi_{2i}(T))$\\
\end{tabular}
\end{center}
The $|\nu|$ term in the $(\GL_n, \Or_n)$ and $(\GL_{2n}, \Sp_{2n})$ cases comes from the contribution of $w_0(\wt(T))$ to $|\lambda_{\min}(T)|$. In the $(\GL_n^2, \GL_n)$ case, $\wt(T) = 0$, so no $|\nu|$ term appears.
\pagebreak 

During the review of this work, the author used ChatGPT to assist with language and copyediting. The author reviewed and edited the output as needed and takes full responsibility for the content of the published article. 
\raggedright 
\printbibliography

\end{document}